%\magnification1200
\let\ReallyTheEnd=\end
\def\end{\vfil\eject}
\def\IMSmarkvadjust{0 pt}
\def\IMSmarkhadjust{0 pt}
\def\IMSmarkhpadding{0 pt}
\def\IMSpubltext{Published in modified form:}
\def\SBIMSMark#1#2#3{
 \font\SBF=cmss10 at 10 true pt
 \font\SBI=cmssi10 at 10 true pt
 \setbox0=\hbox{\SBF \hbox to \IMSmarkhpadding{\relax}
                Stony Brook IMS Preprint \##1}
 \setbox2=\hbox to \wd0{\hfil \SBI #2}
 \setbox4=\hbox to \wd0{\hfil \SBI #3}
 \setbox6=\hbox to \wd0{\hss
             \vbox{\hsize=\wd0 \parskip=0pt \baselineskip=10 true pt
                   \copy0 \break%
                   \copy2 \break% 
                   \copy4 \break}}
 \dimen0=\ht6   \advance\dimen0 by \vsize \advance\dimen0 by 8 true pt
                \advance\dimen0 by -\pagetotal
	        \advance\dimen0 by \IMSmarkvadjust
 \dimen2=\hsize \advance\dimen2 by .25 true in
	        \advance\dimen2 by \IMSmarkhadjust

%
%   Check for publication info
%
%  \newread\jref
  \openin2=publishd.tex
  \ifeof2\setbox0=\hbox to 0pt{}
  \else 
     \setbox0=\hbox to 3.1 true in{
                \vbox to \ht6{\hsize=3 true in \parskip=0pt  \noindent  
                {\SBI \IMSpubltext}\hfil\break
                \input publishd.tex 
                \vfill}}
  \fi
  \closein2
  \ht0=0pt \dp0=0pt
 \ht6=0pt \dp6=0pt
 \setbox8=\vbox to \dimen0{\vfill \hbox to \dimen2{\copy0 \hss \copy6}}
 \ht8=0pt \dp8=0pt \wd8=0pt
 \copy8
 \message{*** Stony Brook IMS Preprint #1, #2. #3 ***}
}

\SBIMSMark{1992/11}{July 1992}{}
%\magnification=1200
\overfullrule=0pt
\def\oldQP{\narrower\smallskip\noindent}
\def\QP{\medskip\leftskip=.4in\rightskip=.4in\noindent}
\def\[{$\,}
\def\]{\,$}
\def \QED { \rlap{$\sqcup$}$\sqcap$\smallskip}
\def\mod{{\rm mod}}
\def\C{{\bf C}}
\def\inter{{\rm interior}}
\def\ref{\smallskip\hangindent=1pc \hangafter=1 \noindent}
\def\scd{{\rm scd}}
\font\bit=cmssi12 at 12truept
\font\tenmsy=msym10
\textfont8=\tenmsy
\mathchardef\ssm="7872
\abovedisplayskip=6pt
\belowdisplayskip=6pt
\input psfig
% This file is supposed to imitate the -ras.tex- file which contained
% the macros -insertRaster- and -RasterBox- for including SUN raster files
% in TeX.  But here we assume that the original raster files have been
% converted to PostScript, so we use -psfig- rather than -sunbitmap-.
% There was a bug in the way -ras.tex- scaled these bitmaps, and this file
% attempts to preserve this bug so that the resulting figures are identical
% to the raster figures in the original paper.

\input psfig
%%\magnification=1200

\newif\ifboxfigure      % set to true if you want the figure boxed
\boxfigurefalse

\def\BoxIt#1#2{%        % put a box around #1, leaving a gap of #2
	\vbox{\hrule
	\hbox{\vrule\kern#2\vbox{\kern#2#1\kern#2}\kern#2\vrule}
		   \hrule}}

\def\insertRaster #1 pixels #2  by #3 scaled #4 {
%  #1 is the raster file
%  #2 is the width of the picture in pixels
%  #3 is the height of the picture in pixels
%  #4 is scaling 500 means half size 2000 means double size 1000 is actual  
			\medskip
			 \hbox to \hsize{%

			 \hss
			 \RasterBox {#1} {#2} {#3} {#4}
			 \hss
			 }%
}

\def\RasterBox #1 #2 #3 #4{

% dimen0 and dimen1 are the dimensions of the hbox which are setup in
% the same way as -ras.tex- unfortunately psfig includes images in a 
% different manner from sunbitmap, so it must be scaled using different
% values, hence dimen2 and dimen3 and MAGIC number 933

\dimen0=#2pt   % #2 is in pixels, so must multiply by 1.1076 when
			  % used
\dimen1=#3pt
\divide\dimen0 by 11076 % really want 1.1076
\multiply\dimen0 by 10 % (because #4 is out of 1000 )
\divide\dimen1 by 11076
\multiply\dimen1 by 10
 \dimen2=#3pt
 \divide\dimen2 by 933
 \dimen3=#2pt
 \divide\dimen3 by 933
\setbox4=\hbox to #4\dimen0{
 \vbox to #4\dimen1{
 \vss
 \psfig{figure=#1,height=#4\dimen2,width=#4\dimen3}
 }
 \hss
 }
 \ifboxfigure\BoxIt{\box4}{0pt}
 \else\box4
 \fi

 }

% \insertRaster towrrab.pic pixels 500 by 500 scaled 400

%table of content format
\def\leaderfill{{\leaders\hbox to 1em{\hss.\hss}\hfill}}
\def\eleaderfill{\hfill\hbox to 1em{\hss\hss}\hfill}

\centerline{\bf Local Connectivity of Julia Sets: Expository Lectures.}
\smallskip

\centerline{J. Milnor}\smallskip

\centerline{Stony Brook, July 1992}\bigskip

\centerline{\bf Contents}\smallskip
{\QP\line{\S1. Local Connectivity of Quadratic Julia Sets\hfil}
{\indent (following Yoccoz)\leaderfill 2}
\line{\S2. Polynomials for which All But One of the Critical Orbits
Escape\hfil}
{\indent (following Branner and Hubbard)\leaderfill 17}
\line{\S3. An Infinitely Renormalizable Non Locally Connected Julia Set\hfil}
{\indent(following Douady and Hubbard)\leaderfill 30}\break
{Appendix: Length-Area-Modulus Inequalities\leaderfill 36}\break
{Errata\leaderfill 42}\break
{References\leaderfill 45}
\bigskip}

\centerline{\bf Introduction}\smallskip

The following notes provide an introduction to recent work of Branner,
Hubbard and Yoccoz on the geometry of polynomial Julia sets. They
are an expanded version of lectures
given in Stony Brook in Spring 1992. I am indebted to help
from the audience.\smallskip

Section 1 describes unpublished
work by J.-C. Yoccoz on local connectivity of quad\-ratic Julia sets.
It presents only the ``easy'' part of his work, in
the sense that it considers only non-renormalizable polynomials, and makes no
effort to describe the much more difficult arguments which are needed
to deal with local connectivity in parameter space. It is based on
second hand sources, namely Hubbard [Hu1] together with lectures by Branner
and Douady. Hence the presentation is surely quite different from
that of Yoccoz.\smallskip

Section 2 describes the analogous arguments
used by Branner and Hubbard [BH2] to study higher degree polynomials for which
all but one of the critical orbits escape to infinity. In this case, the
associated Julia set \[J\] is never locally connected. The basic problem
is rather to decide when \[J\] is
totally disconnected. This Branner-Hubbard
work came before Yoccoz, and its technical
details are not as difficult. However, in
these notes their work is presented simply as another application of the same
geometric ideas.\smallskip

Chapter 3 complements the Yoccoz results by describing a family of examples,
due to Douady and Hubbard (unpublished), showing that an
infinitely renormalizable quadratic polynomial may have
non-locally-connected Julia set. An Appendix describes needed tools
from complex analysis, including the Gr\"otzsch
inequality.\medskip

We will assume that the
reader is familiar with the basic properties of Julia sets and the Mandelbrot
set. (For general background, see for example [Be], [Br2], [D1], [D2], [EL],
[L1], as well as the brief outline in \S3.)
In particular, we will make use of {\bit external rays\/}
for a polynomial Julia set \[J(f)\subset\C\]. (Compare [DH1], [DH2], [M2],
[GM].)% For general background concerning Julia sets and the Mandelbrot set,
%see for example [Be], [Br2], [D1], [D2], [EL], [L1].

\vfil\eject %\bigskip

\centerline{\bf \S1. Local Connectivity of Quadratic Julia Sets (following
Yoccoz).} %\footnote{${}^1$
\medskip

This section will prove the following. % Let \[f_c(z)=z^2+c\].
\medskip

{\leftskip=.5in\rightskip=.8in\noindent{\bf
Theorem 1.} \it If \[f_c(z)=z^2+c\] is a quadratic polynomial such that:
%\noindent\hskip .1in

$(1)$ the Julia set \[J(f_c)\] is connected,

$(2)$ both fixed points are repelling, and

$(3)$ \[f_c\] is not renormalizable\footnote{$^1$}
{\rm For the definition of renormalizability, see Figure 10 and the
associated discussion.} %  $\,($see Figure 10 and its context$\,),$

\noindent then \[J(f_c)\] is locally connected.\medskip}

In terms of the familiar
parameter space picture for the family of quadratic maps \[f_c(z)=z^2+c\],
Condition (1)
says that the parameter value \[c\] belongs to the Mandelbrot
set \[M\], while (2)
says that \[c\] does not belong to the closure of the
central region bounded by the cardioid, and (3) says that
\[c\] does not belong to any one of the many small copies of \[M\]
which are scattered densely around the boundary of \[M\]. (Compare Figure 1.)
\bigskip

\midinsert
\hbox to \hsize{ 
\hfil 
\hbox{\RasterBox {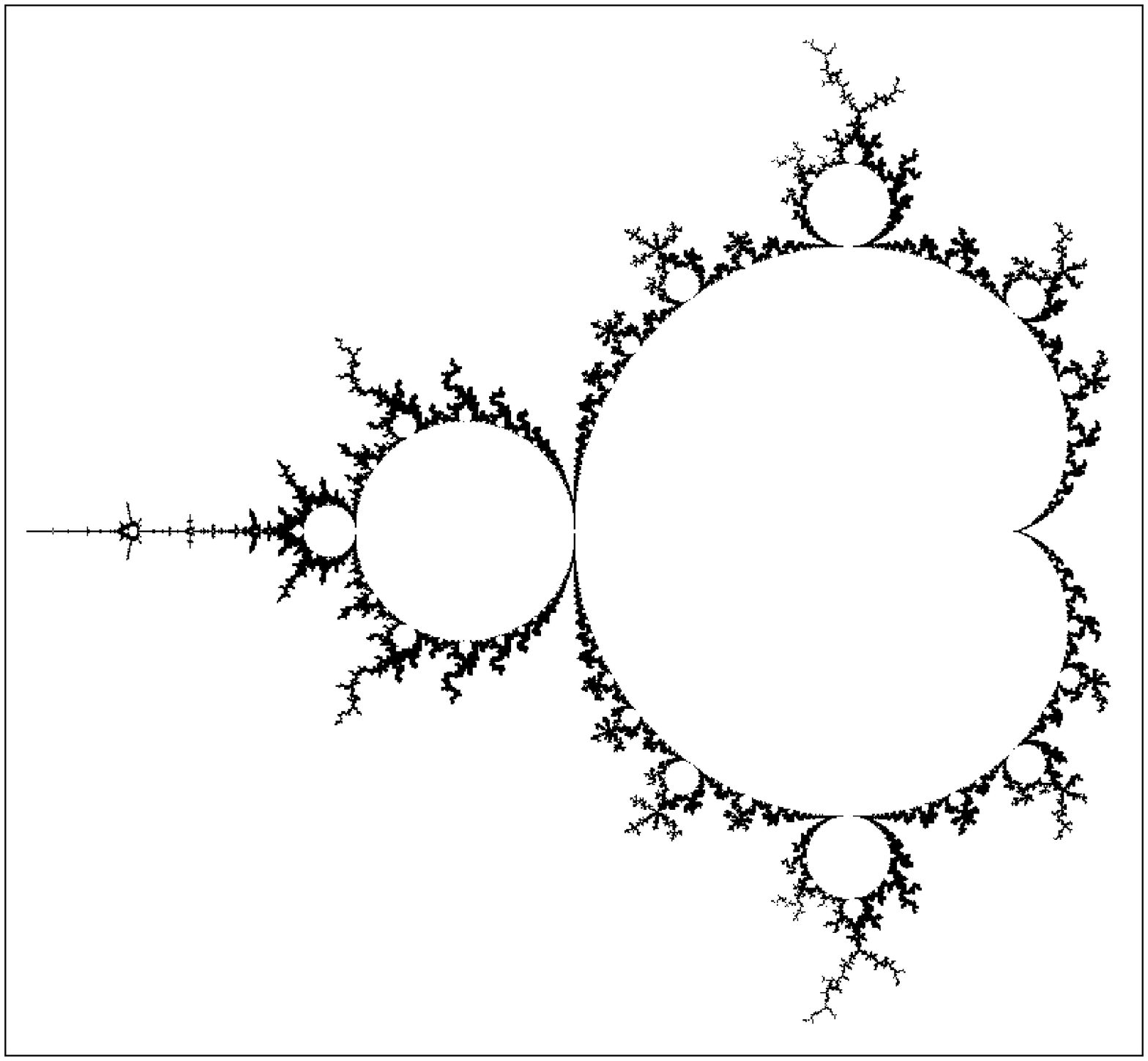} {832} {768} {225}}
\hskip .2in 
\hbox{\RasterBox {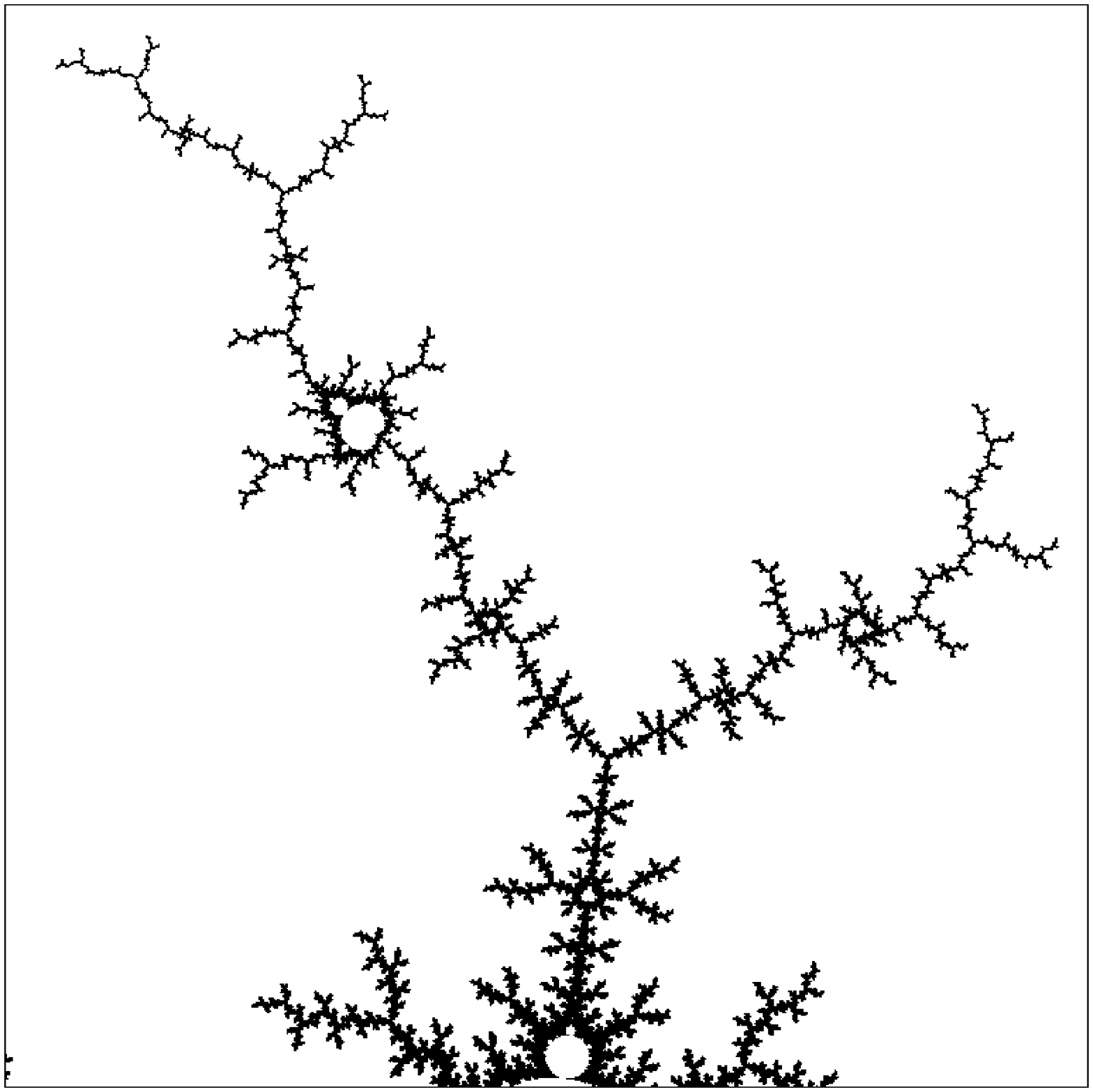} {768} {768} {225}}
\hfil 
}
\smallskip
{\QP \bit Figure 1. Boundary of the Mandelbrot
set \[M\], and a detail of the one-third limb
showing several small copies of \[M\].}
\endinsert

{\bf Remarks:} The proof will give much more, since it will
effectively describe the
Julia set by a new kind of symbolic dynamics. With a little more work,
the Yoccoz method can also deal with the finitely renormalizable case.
Since the case of a map with attracting or parabolic fixed point had been
understood much earlier [DH2], we see that
Conditions (2) and (3) can be actually
replaced by the following weaker pair of conditions:\smallskip

$(2')$ \[f\] has no irrationally indifferent periodic points
(Cremer or Siegel points), and\smallskip

$(3')$ \[f\] is not {\it infinitely\/} renormalizable.\smallskip

\noindent These modified Conditions (1), $(2')$ and $(3')$ are all essential. In
fact,
\eject

$(1)\!:$ For \[c\not\in M\] the Julia set \[J(f_c)\] is a Cantor set, which
is certainly not locally connected.
%Fatou and Julia showed that a non-connected Julia set has infinitely
%(in fact uncountably) many distinct connected components, so such a Julia set
%certainly cannot be locally connected. 

$(2')\!:$ Sullivan and Douady showed
that a polynomial Julia set with a Cremer point is never locally connected.
(Compare [Su], [D1], [M2]. For a more explicit description of
non local connectivity, see [S\o].) Similarly,
Herman has constructed quadratic polynomials with a Siegel disk
having no critical point on its boundary. The corresponding Julia set cannot
be locally connected. (See [He,\S17.1], [D1,II,5], [D4].)

$(3'):$ In \S3,
following unpublished work by Douady and Hubbard, we will describe
infinitely renormalizable polynomials for which \[J\] is not locally connected.
\smallskip

Thus the sharpened version of Theorem 1 comes fairly close to deciding
exactly which quadratic polynomial Julia sets are locally connected.
Yoccoz has also proved a corresponding result in parameter space: For \[c\]
in the Mandelbrot set \[M\], if the polynomial \[f_c(z)=z^2+c\] is not
infinitely renormalizable, then
\[M\] is locally connected at \[c\]. For varied proofs
of this more difficult result, see [Hu2], [HF], [K].
%\bigskip
%\vfil\eject

\midinsert
\insertRaster 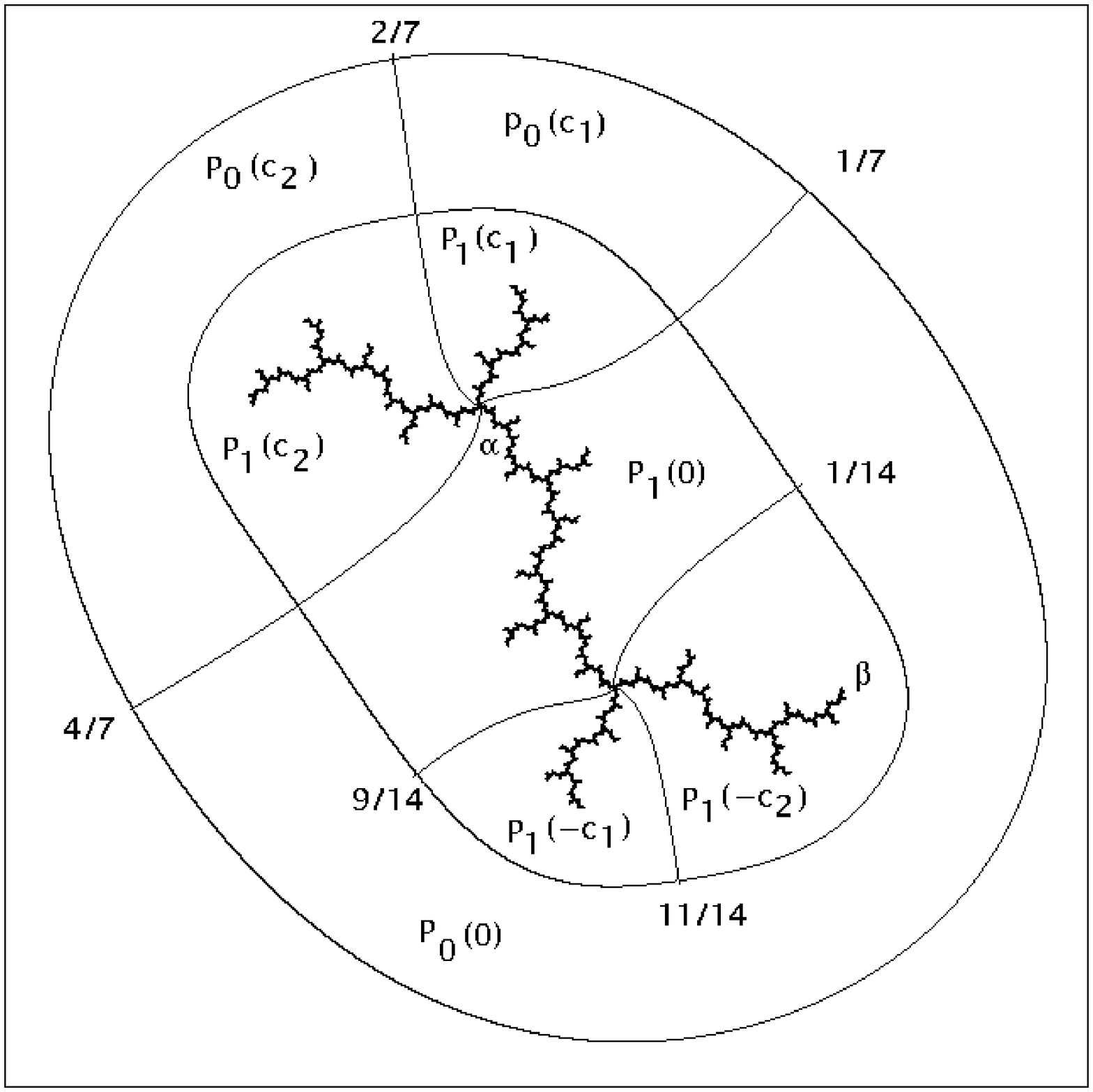 pixels 640 by 640 scaled 450
{\QP{\bit Figure 2.} \bit Julia set for \[f(z)=z^2+i\],
showing the Yoccoz puzzles of depth zero and one. Here \[q=3\].
The \[1/7\,,\;2/7\] and \[4/7\] rays land at the fixed point
\[\alpha\].\medskip}
\endinsert

\vskip -.2in
{\bf The Yoccoz jigsaw puzzle.} The proof of Theorem 1 begins as follows.
If a connected quadratic Julia set has
two repelling fixed points, then one fixed point
(to be called \[\beta\])
is the landing point of the zero external ray, and the other (called
\[\alpha\]) is the landing point of a cycle of \[q\]
external rays, where \[q\ge 2\]. (Compare [M2], [Pe].)
Let\break \[0=c_0\mapsto c_1\mapsto\cdots\]
be the critical orbit.
The {\bit Yoccoz puzzle\/} of {\bit depth zero\/} consists
of \[q\] pieces \[P_0(c_0)\,,\,P_0(c_1)\,,\,\ldots\,,\,P_0(c_{q-1})\]
which are obtained by
cutting the region \[G\le 1\] along the \[q\]
external rays landing at \[\alpha\]. Here  \[G\] is the canonical
potential function. The various pieces have been labeled so that each
\[P_0(c_i)\] contains the post-critical point \[c_i=f^{\circ i}(0)\].

{\bf Inductive Construction:}
If \[P_d^{(1)}\,,\,\ldots\,,\,P_d^{(m)}\]
are the puzzle pieces of depth \[d\], then the connected components of the
sets \[f^{-1}(P_d^{(i)})\] are the puzzle pieces \[P_{d+1}^{(j)}\] of
depth \[d+1\].
As an example, there are always \[2q-1\] pieces of depth one.
These consist of
\[q\] pieces \[P_1(c_0)\,,\,P_1(c_1)\,,\,\ldots\,,\,P_1(c_{q-1})\]
which touch at the fixed point \[\alpha\], together with \[q-1\] additional
pieces \[P_1(-c_1)\,,\,\ldots\,,\,P_1(-c_{q-1})\] which touch at the
pre-image \[-\alpha\].\smallskip

{\bf The Main Problem:} Let \[z\in J(z)\] be any point whose forward orbit
never hits\footnote{$^1$}
{The slight modifications of the argument needed to deal with the case
\[f^{\circ d}(z)=\alpha\] are described near the end of this section.}
the fixed point \[\alpha\], and
%Let \[P_d(z)\] be the\footnote{${}^1$}
%{When using this notation, we must always assume that the orbit of \[z\]
%is disjoint from the fixed point \[\alpha\], so that the puzzle pieces
%\[P_d(z)\] are uniquely defined.
%The case \[f^{\circ d}(z)=\alpha\] is briefly
%discussed at the end of this section.}
let \[P_d(z)\] be the unique puzzle piece of depth \[d\]
which contains this point \[z\], so that
$$ P_0(z)~\supset~ P_1(z)~\supset~ P_2(z)~\supset~\cdots~.$$
{\bit Does the intersection \[\bigcap_d P_d(z)\] consist of the single point
\[z\]\/}?
\medskip

{\bf The associated annuli.} Let \[P_d(z)\supset P_{d+1}(z)\] be the
puzzle pieces of two consecutive depths containing some given \[z\in J(f)\].
If we are lucky, the smaller puzzle piece \[P_{d+1}(z)\] will be contained
in the interior of \[P_d(z)\]. In this case the difference
$$	A_d(z)~=~\inter(P_d(z))\ssm P_{d+1}(z) $$
is an annulus, whose modulus \[~\mod\, A_d(z)~\] is a positive real number.
(Compare the Appendix.)
For example in Figure 2 the annulus \[A_0(-c_1)\] has positive modulus.
On the other hand, it may happen that \[P_{d+1}(z)\] intersects the
boundary of \[P_d(z)\]. In this case we describe \[A_d(z)\] as a
{\bit degenerate annulus\/}, and define its modulus to be
zero. For example, in Figure 2
the ``annulus'' \[A_0(0)\] around the critical point is degenerate.\smallskip

{\bf Modified Main Problem (Branner and Hubbard) :} Given a point
\[z\in J(f)\],\break is
the sum \[\sum_d
\mod\,A_d(z)\] infinite? If so, using the Gr\"otzsch inequality, it is
not difficult to prove that the intersection \[\bigcap P_d(z)\] consists
of the single point \[z\]. (See Appendix.)\smallskip

\midinsert
\hbox to \hsize{
\hfil
\hbox{\RasterBox {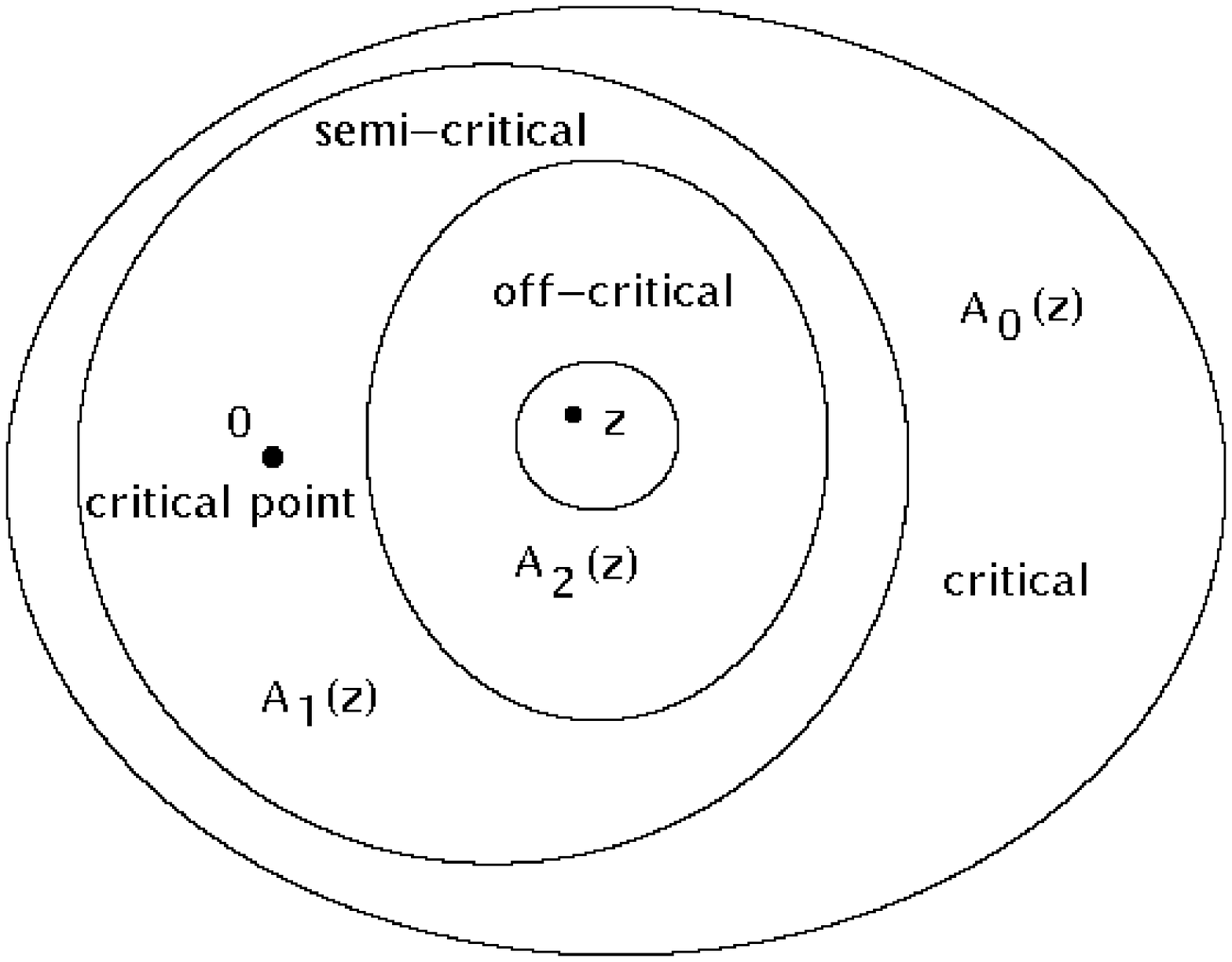} {640} {480} {350}}
\hskip -.5in
\hbox{\RasterBox {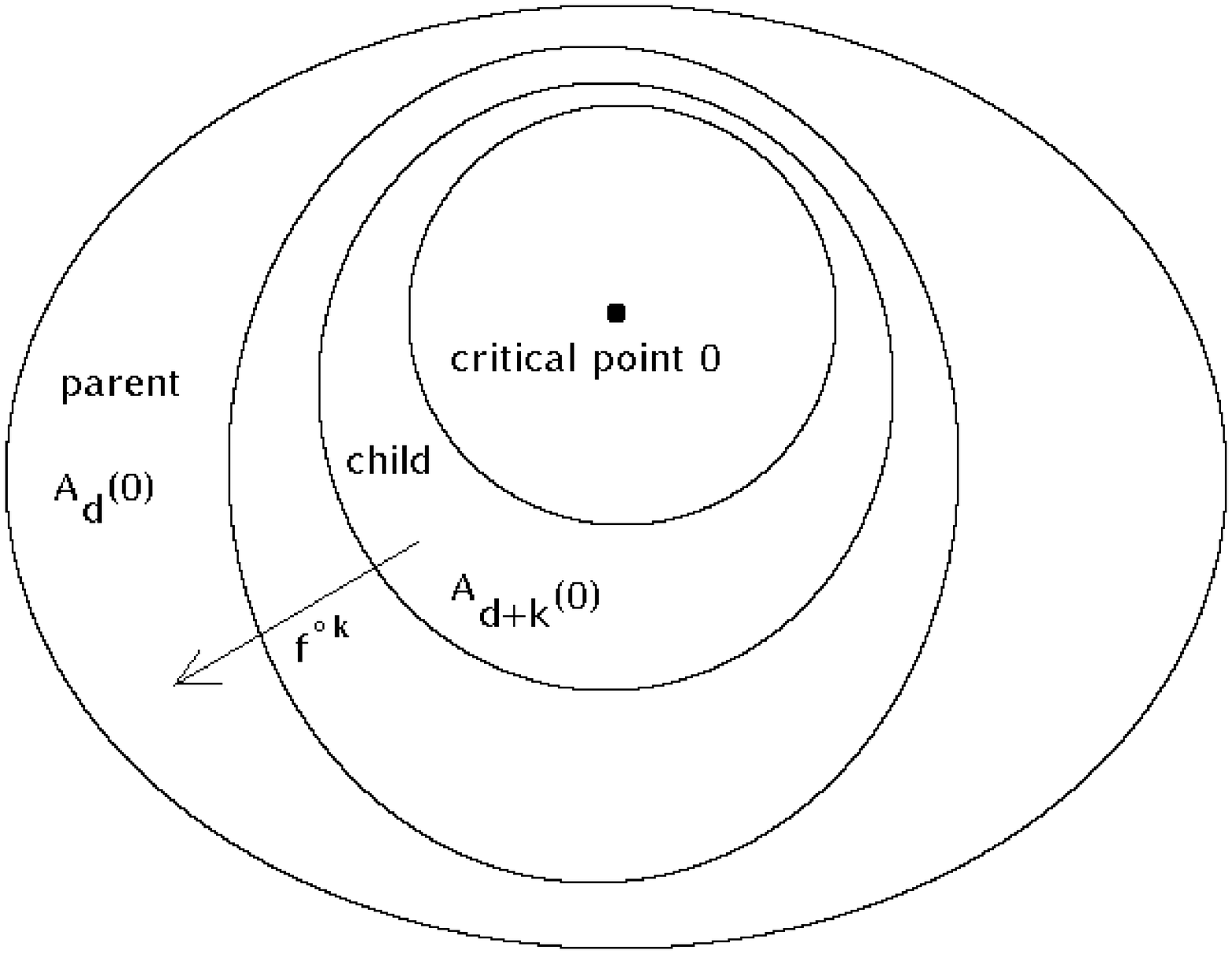} {640} {480} {350}}
\hfil
}
\centerline{\bit Figure 3. Critical, semi-critical and off-critical annuli
$($left$)$, and a ``child'' $($right$)$.} \endinsert

Note that \[f\] carries any puzzle piece \[P_d(z)\] of depth \[d>0\] onto
\[P_{d-1}(f(z))\]. This map \[~P_d(z)\to P_{d-1}(f(z))~\]
is either a conformal isomorphism or a two-fold
branched covering according as \[P_d(z)\] does or does not contain the
critical point.\smallskip

{\bf Definition.} We describe an annulus as being either {\bit
semi-critical\/} or
{\bit critical\/} or\break {\bit off-critical\/}
according as the critical point belongs to the annulus itself or to the
bounded  or the unbounded component of its complement.
Thus, in the schematic picture (Figure 3 left), the annulus
\[A_0(z)\] is critical, while \[A_1(z)\] is semi-critical,
and \[A_2(z)\] is off-critical.
Let \[A_d(z)\] be the annulus of depth
\[d>0\] in the Yoccoz puzzle surrounding a point \[z\in J(f)\].
If \[A_d(z)\] is critical or off-critical, then evidently \[f\] maps \[A_d(z)\]
onto \[A_{d-1}(f(z))\]:
\vskip -.1in

\hskip .7in by an unramified two-fold covering in the critical case,\smallskip

\hskip .7in by a conformal isomorphism in the off-critical case.\smallskip

%\hskip .25in{conformally isomorphic to \[A_{d-1}(f(z))\] whenever \[A_d(z)\]
%is off-critical, and}\smallskip

%\hskip .25in{an unramified two-fold covering of \[A_{d-1}(f(z))\] whenever
%\[A_d(z)\] is critical.}\smallskip

%\hskip .25in{a ramified two-fold covering when \[A_d(z)\]
%is semi-critical, and}\smallskip

\noindent
In the first case it follows that \[~\mod\,A_d(z)={1\over 2}\,
\mod\,A_{d-1}(f(z))\,\],
while in the second case \[~\mod\,A_d(z)=\mod\,A_{d-1}(f(z))\,\].
In the semi-critical case, the appropriate statement is that
\[~\mod\,A_d(z)>{1\over 2}\,\mod\,A_{d-1}(f(z))\,\].
% the modulus of \[A_d(z)\] is half
%the modulus of \[A_{d-1}(f(z))\], and in the second case
%these two moduli are equal. 
The proof is more complicated.
(Compare Problem 1-2 at the end of this section.) Even in the semi-critical case,
it is easy to check that:\break \[A_d(z)\] {\it has positive modulus if
and only if \[A_{d-1}(f(z))\] has positive modulus.}
\medskip

{\bf Definition (Figure 3 right).}
The critical annulus \[A_{d+k}(0)\] in
the Yoccoz puzzle will be called
a {\bit child\/}
\footnote{${}^1$}
{Our terminology is based on [Hu1], but with several modifications.
Thus in Hubbard's terminology, the level \[d+k\] is called
a ``legitimate child'' of the level \[d\]. Similarly, we have replaced
Hubbard's ``marked grid'' by {\bit tableau\/}, and his
``noble'' by {\bit excellent\/}.}
of the critical annulus \[A_d(0)\]
if and only if \[A_{d+k}(0)\] is an unramified two-fold covering of \[A_d(0)\]
under the map \[f^{\circ k}\].\smallskip

Note that the modulus of such a child is always
exactly half the modulus of the parent. Our strategy for solving the
Modified Main Problem can now be summarized as follows:

%\pageinsert

{\QP\bit Find a
critical annulus of positive modulus, and prove that it has so many
descendents, children and
grandchildren and so on, that the modulus sum is infinite.\medskip}

%\centerline{\psfig{figure=tabdef.ps,height=4.5in}}
%\smallskip
{\bf Definition.} The {\bit tableau}$^1$ associated with an orbit \[~z_0\mapsto
z_1\mapsto z_2\mapsto\cdots\]
in \[J(f)\ssm\{\alpha\}\] is an array with one column associated with
each \[z_i\] and one row associated with each depth in the Yoccoz puzzle.
We will draw a solid vertical line at depth \[d\] in the\break \[j$-th column to
indicate
that the annulus \[A_d(z_j)\] in the Yoccoz puzzle coincides with the
critical annulus \[A_d(0)\]. A double vertical line will indicate that \[A_d(z_j)\]
is semi-critical, and no line at all will indicate
that it is off-critical. % Thus a column
Diagonal arrows (pointing north-east)
correspond to iterates of the map \[f\]. Thus annuli along the same
diagonal line either all have zero modulus or all have non-zero modulus.
\smallskip
This concept of tableau, due to Branner and Hubbard, provides a language
which we will use to describe the Yoccoz proof. (It is not the
language which Yoccoz himself uses.)
\smallskip

\midinsert
\centerline{\psfig{figure=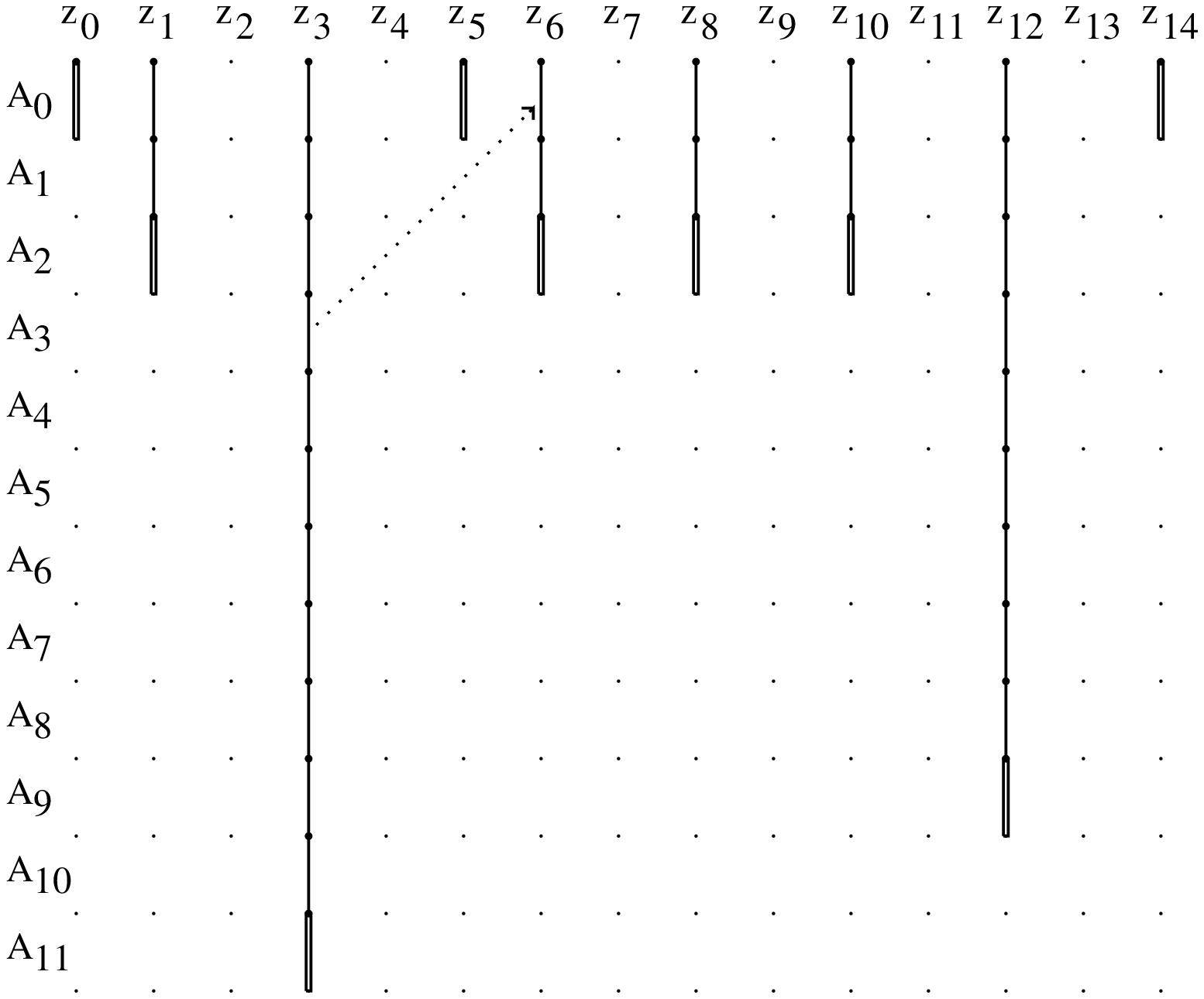,height=3.8in}}
\vskip -.05in
{\QP\bit Figure 4. An example: the tableau associated with the orbit of the point\break
\[z_0=1\] for the map \[f(z)=z^2-1.6\]. Here we see by following the
diagonal arrow
that the critical annulus \[A_3(z_3)=A_3(0)\]
is an unramified two-fold covering of the critical annulus\[A_0(z_6)=
A_0(0)\]. In particular, the critical annulus at depth \[3\]
is a child of the critical annulus at depth \[0\]. Similarly, the critical
annuli at depths \[4\,,\,6\,,\,8\,,\,10\] are children
of the critical annulus of depth \[1\].\bigskip}
\endinsert
\eject

We are particularly interested in the tableau of the critical orbit
$$0=c_0~\mapsto~ c_1~\mapsto~ c_2~\mapsto~\cdots~,$$ which has vertical
line segments going all the way down in column zero. {\bit To simplify
the discussion, we will always assume that
the critical orbit is disjoint from the
fixed point \[\alpha\], so that this critical tableau is well defined.\/}
(If the critical orbit ends at \[\alpha\], then we are in the post-critically
finite case, and local connectivity can be established by other methods.
Compare [DH2].)
%\endinsert

\midinsert
\centerline{\psfig{figure=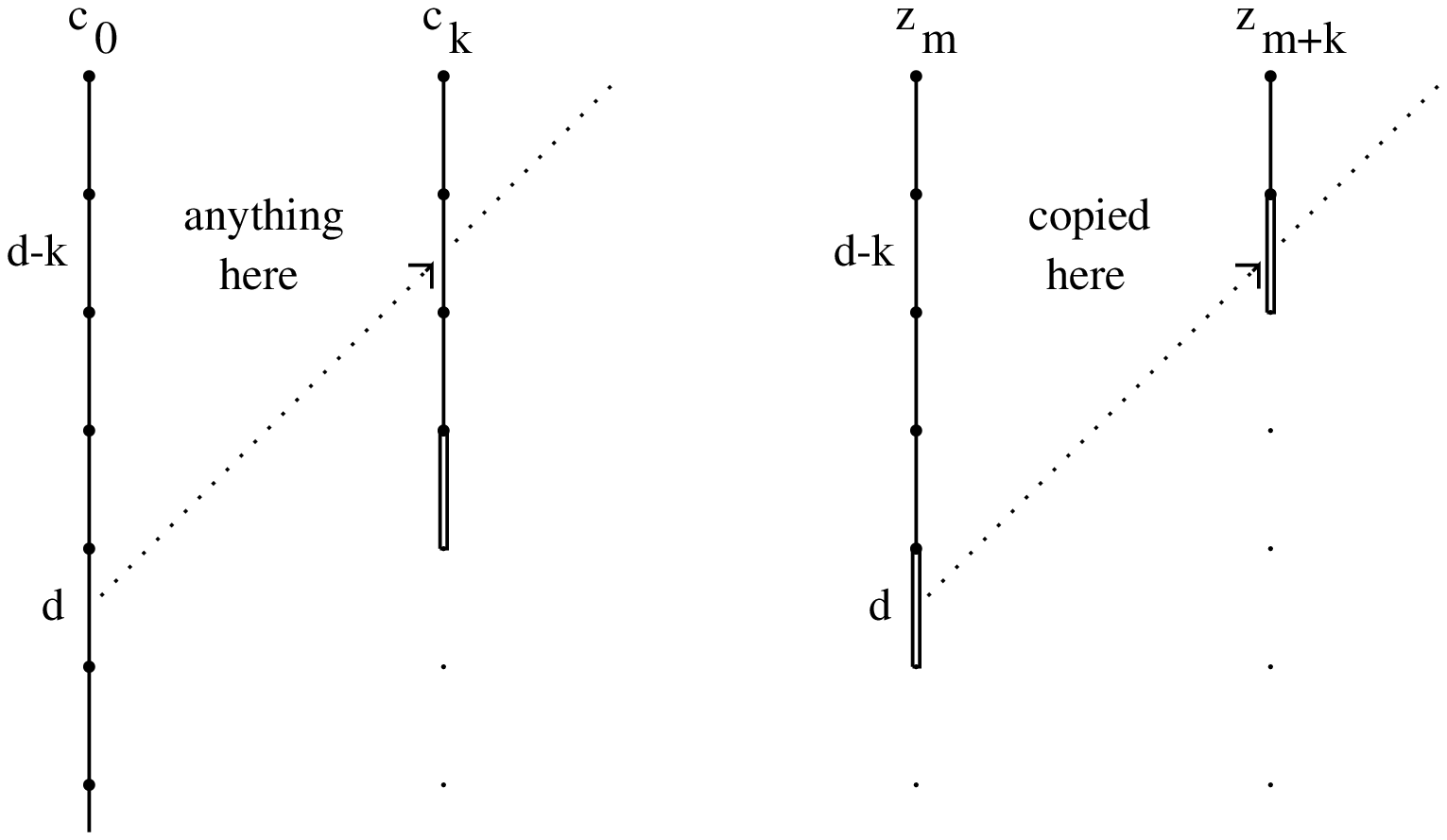,height=2.9in}}

{\QP\bit Figure 5. Illustration for the second and third tableau
rules, with the critical tableau on the left, and the tableau
for \[z_0\mapsto z_1\mapsto\cdots\] on the right.\par}\endinsert

{\QP {\bf First tableau rule:} \it Every column of  a tableau
is either all critical, or all off-critical, or has exactly one
semi-critical depth and is critical above and off-critical below.\medskip}

(Thus each column, corresponding to a point \[~z_j\in J(f)~\], can be
completely described by a single number, the ``{\bit semi-critical depth\/}''
\[~-1~\le~ \scd(z_j)~\le~\infty~,\]
defined by the condition that \[~~P_d(z_j)=P_d(0)~\Longleftrightarrow
~d\le\scd(z_j)~.~\] A large value of \[\scd(z_j)\] means
that \[z_j\,\] is very close to the critical point.)\bigskip
\vfil\eject

Now let us compare the tableau of the critical orbit
$$	0=c_0~\mapsto~ c_1~\mapsto~ c_2~\mapsto~\cdots $$
with the tableau of some given orbit \[z_0\mapsto z_1\mapsto\cdots\] in
\[J(f)\]. (The case \[z_0=c_0\] is not excluded.)
If the tableau of \[\{z_j\}\] is critical or semi-critical at depth \[d\]
in column
\[m\], draw a line ``north-east'' from this critical or semi-critical annulus,
as indicated by the dotted line on the right, and draw a corresponding
line north-east from depth \[d\] of column zero in the critical tableau.

{\QP{\bf Second
tableau rule:\/} {\it Everything strictly above the diagonal line
on the left must be copied above the diagonal line on the right.}\medskip}

The proofs of these two rules are easily supplied.\QED

Now suppose that the critical annulus of depth \[d\] is a child of
the critical annulus of depth
\[d-k\], as indicated in the figure, and suppose that the tableau of
\[\{z_j\}\] is semi-critical at depth \[d\] of column \[m\].

{\QP{\bf Third tableau rule:\/} \it Following
the diagonal arrow from this semi-critical annulus of depth \[d\]
in the tableau of \[\{z_j\}\], we must reach a
semi-critical annulus at depth \[d-k\], as illustrated.\medskip}

{\bf Proof.} According to the hypothesis,
\[f^{\circ k}\] maps \[A_{d}(0)\] onto \[A_{d-k}(0)\], where the
point \[z_m\] is an element of this annulus
\[A_{d}(0)\]. Therefore \[f^{\circ k}(z_m)=z_{m+k}\] must belong to
\[A_{d-k}(0)\].\QED

{\bf Definition.} We will say that a critical annulus \[A_d(0)\] is
{\bit excellent\/} if it contains no post-critical points, or equivalently
if there is no semi-critical annulus in the
\[d$-th row of the critical tableau.\smallskip

\vfil
%\endinsert

\centerline{\psfig{figure=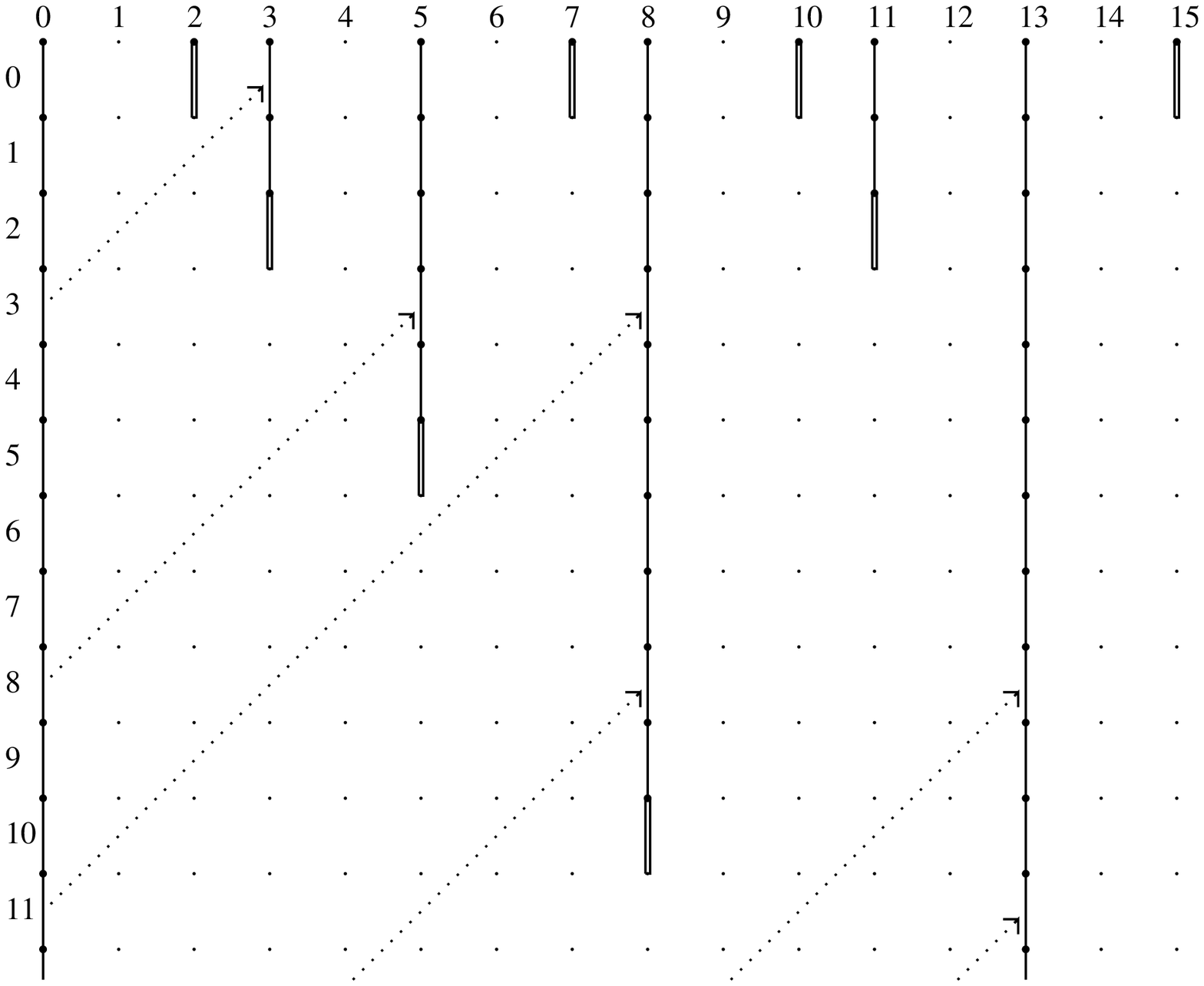,height=4.2in}}

{\oldQP\bit Figure 6. Another example: the Fibonacci tableau.
(Compare [BH2].) Here the closest recurrences
of the critical orbit come after a Fibonacci number of iterations:
$$	1\,,~2\,,~3\,,~5\,,~8\,,~13\,,~21\,,~34\,,~~55\,,~\ldots~. $$
If \[u_n\] is the \[n$-th Fibonacci
number, then the semi-critical annulus for column number \[u_n\]
occurs at depth \[u_{n+1}-3\].
The diagonal dotted lines illustrate the genealogies
\vskip -.15in
% $$	0~\leftarrow~ 3~\leftarrow~ 8~~{\rm and}~~11~\leftarrow~\cdots~, $$
$$	\matrix{
	\qquad \qquad 0 & \leftarrow & 3 & \leftarrow & 8 & \leftarrow & \cdots\cr
	  &            &   & \nwarrow \cr
	  &            &   &   &   11 & \leftarrow & \cdots} %\rm
\eqno ({\bit arrow~from~child~to~parent}).\qquad\qquad $$
\medskip}

In Figure 6,
note that the critical annuli at depths
\[1\,,\;3\,,\;4\,,\;6\,,\;7\,,\;8\,,\; 9\,,\;\ldots\] are excellent.
Each one has exactly two
children, which are also excellent. On the other hand, \[0\,,\;2\,,\;5\,,\;
10\,,\; 18\] (Fibonacci numbers minus three) are not excellent. Each of these
has only one child; however this child is excellent.\smallskip

To see the force of the three tableau rules, suppose that we try to modify
this Fibonacci
tableau by changing just one column. For example, suppose that we place the
semi-critical annulus for column \[5\] at depth \[3\] or \[4\] or at depth
\[\ge 6\], instead of at depth \[5\], without changing columns 0 through 4.
{\bf Exercise:} Show that the tableau rules would then force column 8 to end already
at the semi-critical depth \[0\] or \[1\] or \[2\] respectively.

\eject

\centerline{\psfig{figure=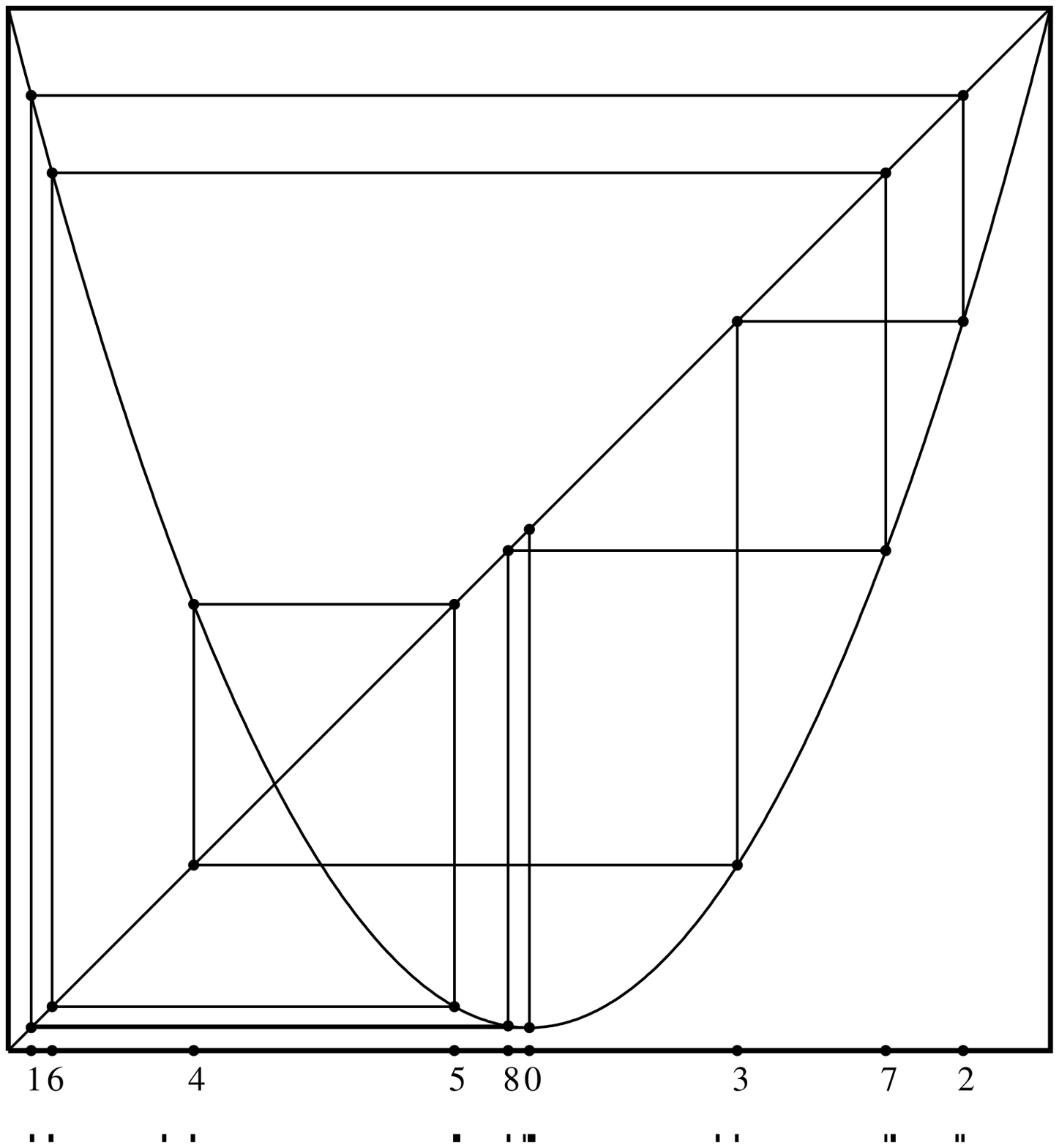,height=3.82in}}
\smallskip
{\oldQP{\bit Figure 7. Graph of the unimodal map \[x\mapsto x^2-1.8705286\cdots\]
which realizes the Fibonacci tableau. (Compare [ML].) The first eight
points on the critical orbit are marked. The critical orbit closure is
a rather thin Cantor set, which is plotted underneath the graph. (Compare
Problem 1-7.)} % By way of contrast for the map \[x\mapsto x^2-1.8\] of Figure 4,
%the critical orbit appears to be dense in the inteval \[[f(0)\,,\,f(f(0))]\].)}
\bigskip}

\insertRaster 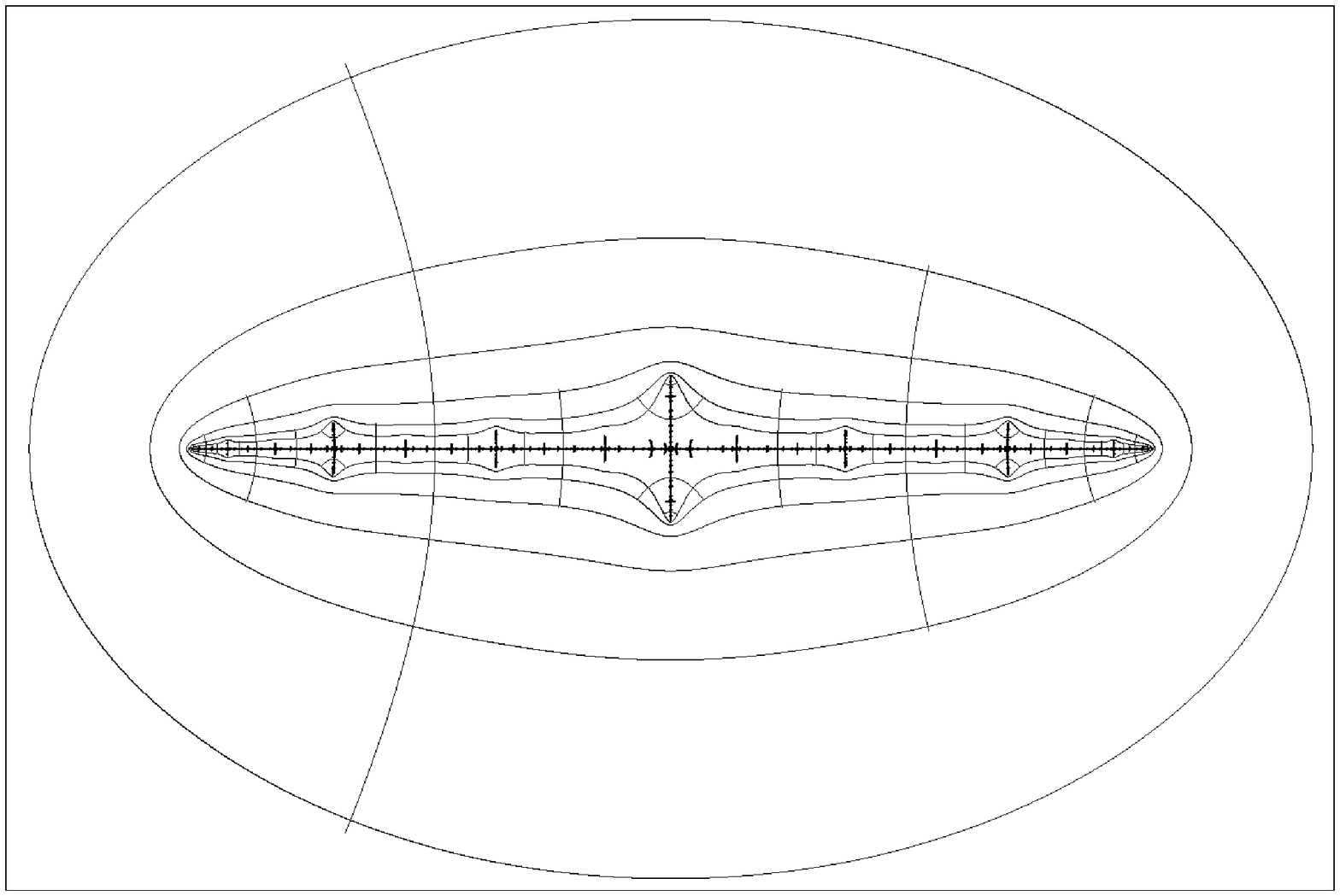 pixels 1440 by 960 scaled 225
\medskip
{\oldQP{\bit Figure 8. The Yoccoz puzzle at depths zero through five
for this quadratic %\break
Fibonacci map, drawn to the same scale. Note that the critical pieces at
any depth are the biggest ones.}
\par}
\eject

\centerline{\psfig{figure=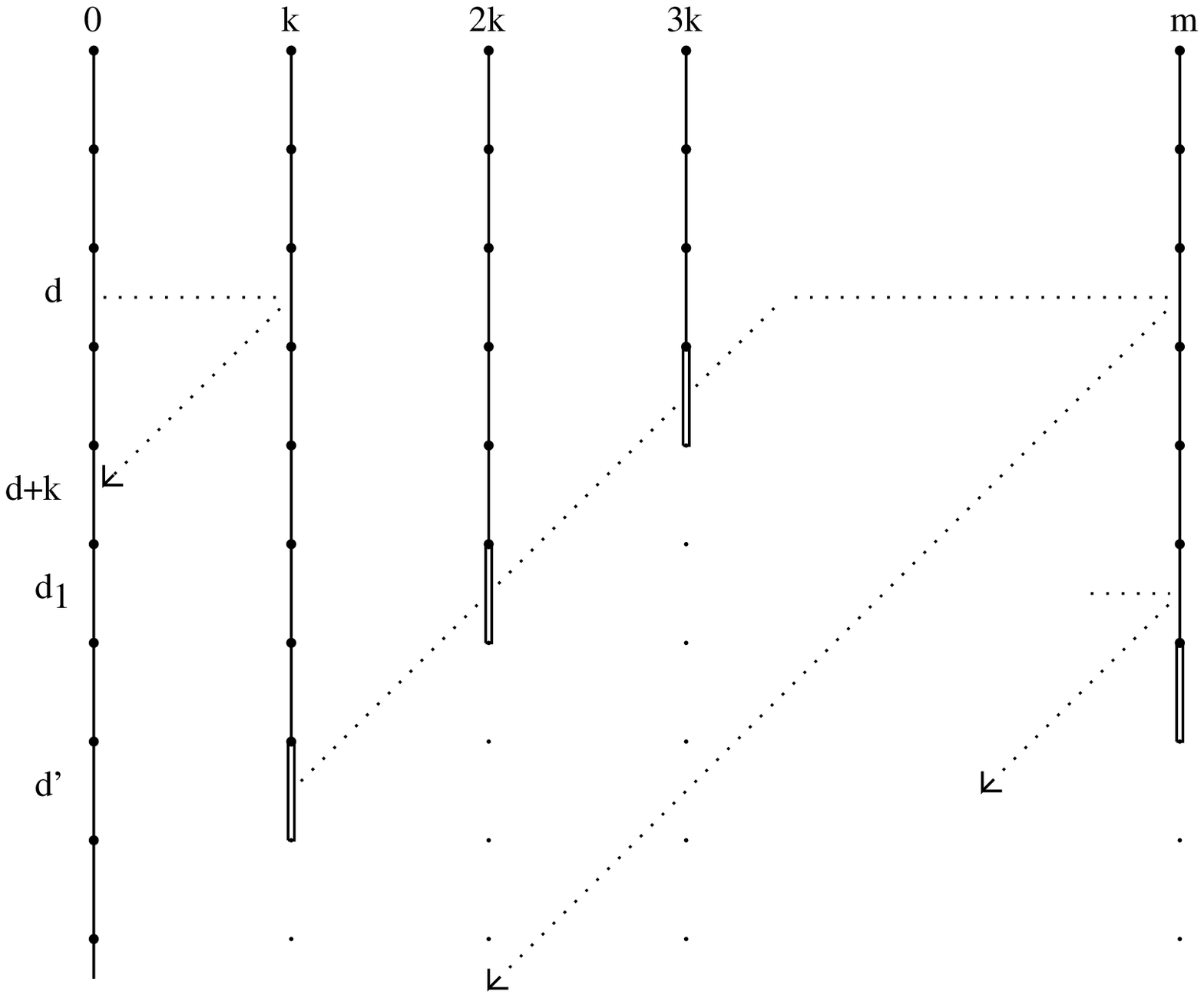,height=3.3in}}
\centerline{\bit Figure 9. Finding Children.}\medskip

Recall that a
critical annulus \[A_d(0)\] is {\bit excellent\/} if the corresponding
tableau row contains no semi-critical entries, or equivalently if the
annulus \[A_d(0)\] contains no post-critical points.
{\bf Definition:} We will say that the critical
tableau is {\bit recurrent\/} if there are columns
with \[k>0\] which go arbitrarily far down; and that it is
{\bit periodic\/} if at least one
such column goes all the way down, so that \[P_d(0)=P_d(c_k)\]
for all depths \[d\]. (Compare Lemma 2 below.)
% \vskip -.05in

{\QP{\bf Lemma 1:} {\it If the critical tableau is recurrent but not periodic,
then:

$(a)$ Every critical annulus has at least one child.

$(b)$ Every excellent critical annulus has at least two children.

$(c)$ Every child of an excellent parent is excellent.

$(d)$ Every only child is excellent.}\medskip}

 {\bf Proof of (a).} Start on the left at depth \[d\],
and march to the right
until we first meet another critical annulus (\[=\]solid vertical line),
say at column \[k\]. Now marching diagonally ``south-west'', we find the
first child at depth \[d+k\].

{\bf Proof of (b).} By hypothesis, the \[k$-th column
cannot be critical all the way down.
Hence it must be semi-critical at some depth \[d'>n\]. Starting at column \[k\]
and depth \[d'\], proceed diagonally right (north-east). By the tableau rules,
column \[2k\] must be semi-critical at depth \[d'-k\]. Similarly
column \[3k\] must be 
semi-critical at depth \[d'-3k\], and so on, until we again reach depth \[d\].
Furthermore, as we follow this diagonal, we do not meet any other
critical or semi-critical annuli. In particular, at the point where
we reach depth \[d\], there cannot be a critical or semi-critical annulus.
(The hypothesis that \[d\] is excellent comes in at this point. Compare
Figure 6.) Now let us 
again march to the right at depth \[d\] until we reach a critical
annulus, say at column \[m\]. Again turning \[135^\circ\]
and proceeding diagonally south-west, we cannot
meet any critical or semi-critical annulus until we
are back at column \[0\]. In this
way we prove that the critical annulus of depth \[m+d\] is also a child
of \[d\].

{\bf Proof of (c).} This is clear, since if
\[f(A_{d+k}(0))=A_d(0)\] where \[A_{d+k}(0)\] contains
a post-critical point, then \[A_d(0)\]
does also.

{\bf Proof of (d), by contradiction.} Consider a child
\[d'\] which is not excellent, and let \[d_1=n'-k\] be the parent. Then
for some \[k'\ge k\] the \[k'$-th column
has semi-critical annulus at depth \[d'\]. (The case \[k'=k\] is illustrated.)
Following the diagonal up from column \[k'\] and depth \[d'\], we must meet
a semi-critical
annulus at depth \[d_1\] by the third tableau rule.
(Thus the parent is not excellent.)
Now proceed to the right at depth \[d_1\]
until we meet a critical annulus, say at column
\[m\]. Then it follows as above that \[d_1+m\] is a second child; hence
\[d'\] was not an only child.\QED

\insertRaster 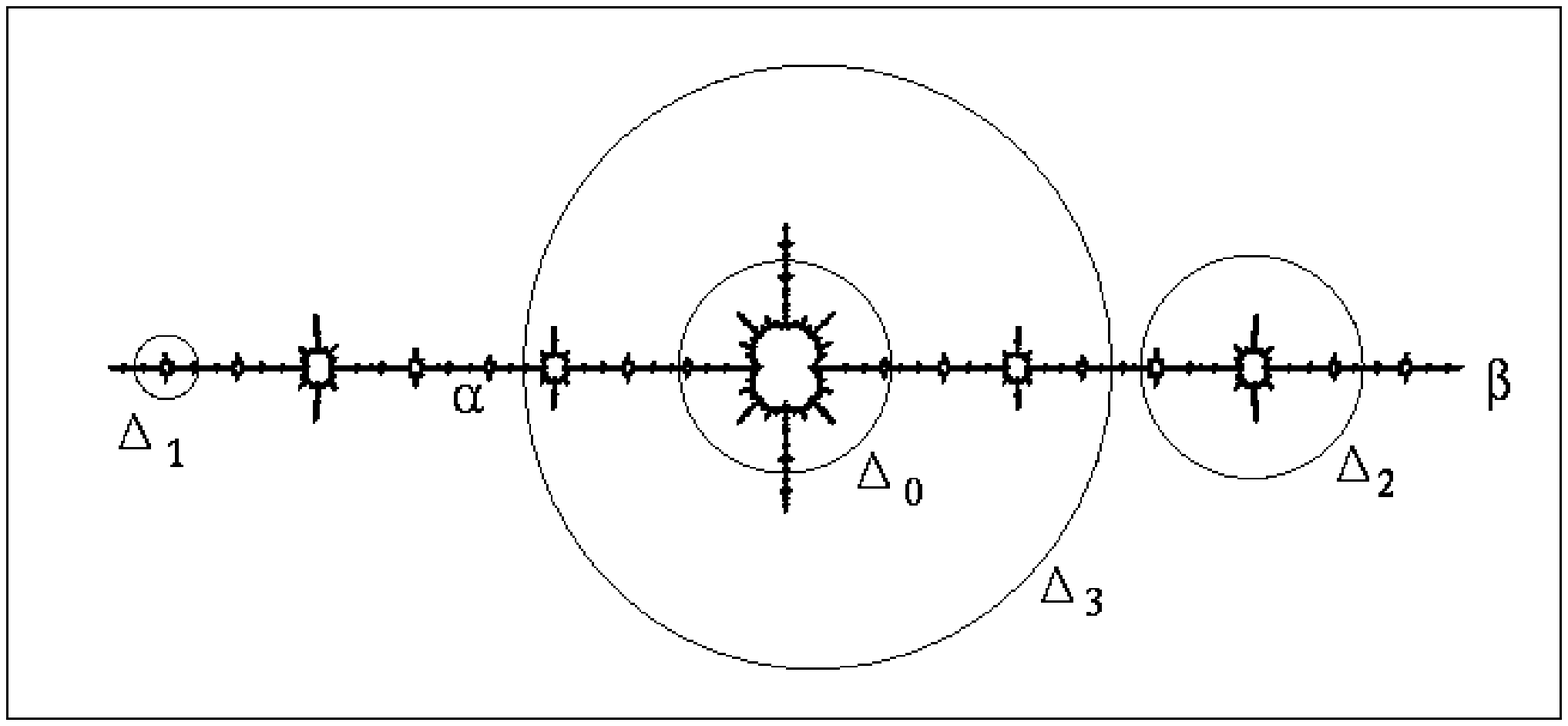 pixels 720 by 330 scaled 450

{\QP\bit Figure 10. Julia set for
\[f(z)=z^2-1.75\], with renormalization period
\[p=3\]. (The parameter value \[c=-1.75\]
is the root point of a small copy of the
Mandelbrot set.) The filled Julia set \[K(f^{\circ 3}|\Delta_0)\]
is a topological disk bounded by
a ``cauliflower''.\medskip}

{\bf Definition.} (Compare [DH3].)
A quadratic polynomial \[f\] is
{\bit renormalizable\/} if there exists an integer \[p\ge 2\] and a closed
topological disk \[\Delta_0\] around the critical point
with the following properties:

\vskip -.05in
{\QP\ref (1) \[\Delta_0\] should be centrally symmetric about the critical point
so that its image \[\Delta_1=f(\Delta_0)\] is also a closed topological disk.

\ref (2) \[f\] should induce conformal isomorphisms
$$	\Delta_1~\buildrel\cong\over\longrightarrow~\Delta_2~
	\buildrel\cong\over\longrightarrow~\cdots~
	\buildrel\cong\over\longrightarrow~\Delta_p~,$$
where \[\Delta_i=f^{\circ i}(\Delta_0)\]. In particular, the critical point
should be disjoint from \[\Delta_1\cup\cdots\cup\Delta_{p-1}\].

\ref (3) However, the disk \[\Delta_p\] should contain \[\Delta_0\]
in its interior.

\ref (4) Finally, the entire orbit of the critical point under
the map \[f^{\circ p}\] should be contained in the original disk \[\Delta_0\].
\smallskip}

\noindent (For further information, see the discussion of ``polynomial-like
mappings'' in \S2, and the discussion of ``tuning'' in \S3.)
We will call \[p\] the {\bit renormalization period\/}.
If these conditions are satisfied, then the
{\bit filled Julia set\/} \[K(f^{\circ p}|\Delta_0)\] can be defined
as the compact connected set consisting
of points whose orbits remain in
\[\Delta_0\] under all iterations of \[f^{\circ p}\]. This is a proper subset
of the filled Julia set \[K(f)\] of the original map \[f\].
We can now state two
basic lemmas.

{\QP{\bf Lemma 2.} \it If the critical tableau associated with
\[f\] is periodic, then \[f\] is renormalizable. More precisely, if
\[P_d(c_p)=P_d(0)\] for all depths \[d\],
and therefore \[P_d(c_{i+p})=P_d(c_i)\] for all \[i\] and \[d\] by the
Second Tableau Rule, then \[f\] is renormalizable of period \[p\].\medskip}

\noindent The converse is also true, but will not be proved here.\smallskip

{\QP{\bf Lemma 3.} \it If
the critical orbit lies completely within the union $$ P_1(c_0)~\cup~
P_1(c_1)~\cup~\cdots~\cup~ P_1(c_{q-1}) $$ of those puzzle pieces of depth one
which touch the fixed point \[\alpha\], then \[f\] is renormalizable
of period \[q\].\medskip}

\hbox to \hsize{
\hskip -.8in
\hbox{\RasterBox {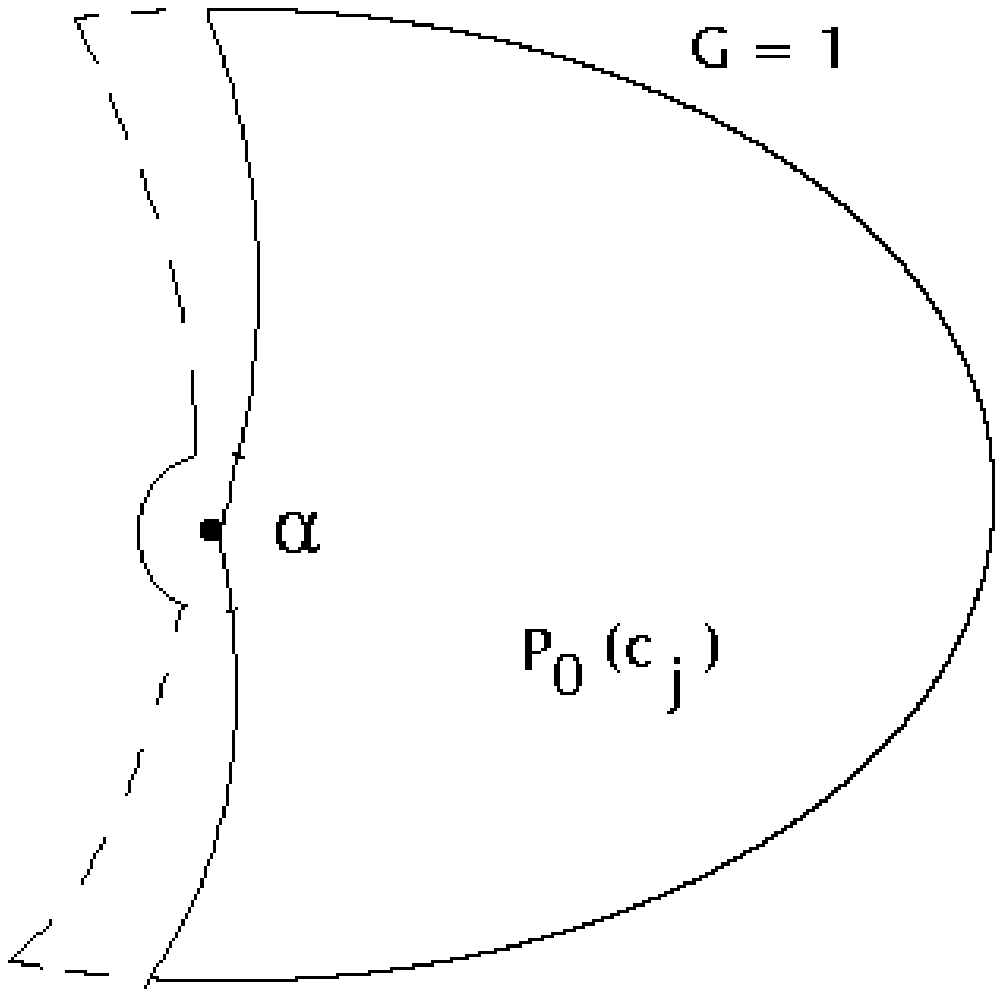} {640} {480} {400}}
\hskip -.3in
\hbox{\RasterBox {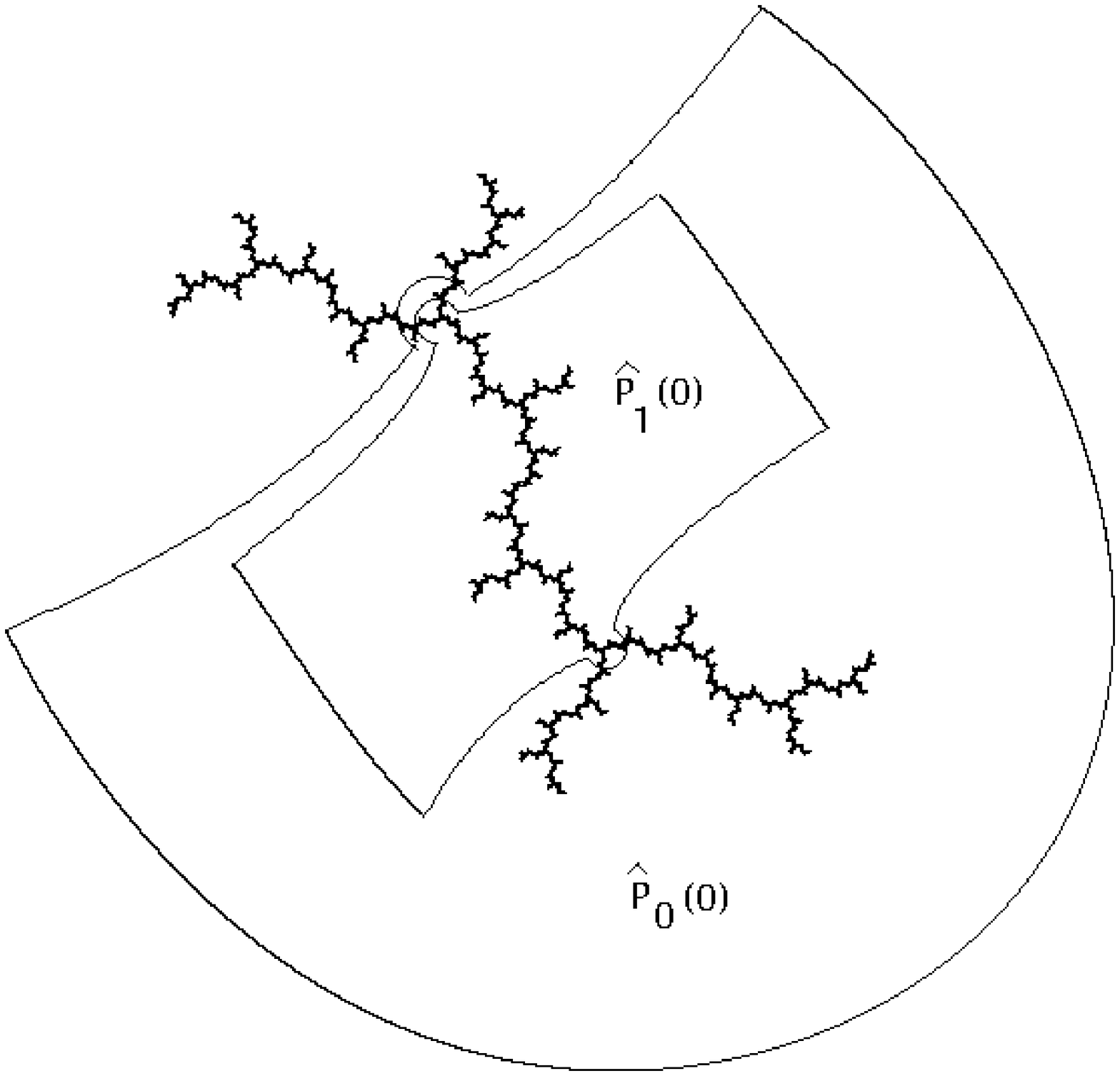} {640} {640} {300}}
\hfil
}
\smallskip
\centerline{\bit Figure 11. Left: A
Puzzle Piece and Thickened Piece of Level \[0\].}

\centerline{\bit Right: Thickened Pieces of Level \[0\] and \[1\]
for \[z\mapsto z^2+i\].}\medskip

The following construction
will be needed to prove Lemmas 2 and 3. Recall that each of the
puzzle pieces \[P_0(c_i)\] of depth zero
consists of points in some closed subset of the filled Julia set \[K(f)\],
together with points outside of \[K(f)\] which have potential \[G\]
and external angle \[t\] satisfying inequalities of the form
$$	0~<~G~\le~1\,,\qquad t_i~\le~t~\le~t'_i~~. $$
(Here \[0\le i<q\].) We construct a
{\bit thickened puzzle piece\/} \[\widehat P_0(c_i)\supset P_0(c_i)\] in two
steps, as follows. First choose a small disk \[D_\epsilon(\alpha)\]
about the fixed point \[\alpha\].
Second, choose \[\eta>0\] so small that every external
ray whose angle differs from \[t_i\] or \[t'_i\] by at most \[\eta\]
intersects this disk. Now let \[\widehat P_0(c_i)\]
consist of the disk \[D_\epsilon(\alpha)\],
together with the region bounded by:

\ref $~~(1)$ the segment \[t_i-\eta\le t\le t'_i+\eta\] of
the equipotential \[G=1\],

\ref $~~(2)$ the external ray segments of angle \[t_i-\eta\] 
and \[t'_i+\eta\] which extend from this equipotential \[G=1\] to their
first intersection with \[D_\epsilon(\alpha)\], and

\ref $~~(3)$ an arc of the boundary of \[D_\epsilon(\alpha)\].\smallskip

\noindent Evidently this thickened puzzle piece \[\widehat P_0(c_i)\]
contains the original \[P_0(c_i)\].
We now construct thickened puzzle pieces of greater depth by the usual
inductive procedure: If \[\widehat P_d^{(j)}\] is a thickened puzzle
piece of depth \[d\], then each component of \[f^{-1}(\widehat
P_d^{(j)})\]
is a thickened puzzle piece of depth \[d+1\].\smallskip

The virtue of these thickened pieces is the following statement, which is
easily proved by induction: {\it If a puzzle piece \[P_d^{(j)}\]
contains \[P_{d+1}^{(k)}\], then the corresponding thickened piece
\[\widehat P_d^{(j)}\] contains \[\widehat P_{d+1}^{(k)}\] in
its interior.\/} In other words, this construction replaces all of our
annuli by non-degenerate annuli.\medskip
%\eject

{\bf Caution:} For this construction to be useful, we also need the following
condition:\break
{\it A thickened puzzle piece \[\widehat P_d(z)\] contains the critical
point only if \[P_d(z)\] already contains the critical point.\/} However, for
a fixed choice of \[\epsilon\] and \[\eta\], this condition will only
be satisfied
for appropriately bounded values of \[d\]. For it may well happen
that the fixed point
\[\alpha\] is an accumulation point of the critical orbit. Thus
we cannot avoid having
points \[c_d\] of the critical orbit within
\[\widehat P_0(c_i)\ssm P_0(c_i)\], hence we cannot avoid having thickened
puzzle pieces \[\widehat P_d(z)\] with \[0\in\widehat P_d(z)\ssm
P_d(z)\], when \[d\] is large.\smallskip

However, we always avoid complication by assuming that the
critical orbit does not actually hit the point \[\alpha\]. This does not
seriously limit the scope of the Theorem, since if \[c_i=\alpha\] then
\[f\] is post-critically finite and
one knows already by [DH2] that \[J(f)\] is locally connected. With this
hypothesis, we can always
choose \[\epsilon\] and \[\eta\] small enough so that the above condition
is satisfied for any specified \[d\].
\medskip

{\bf Proof of Lemma 2.} Choose \[d\ge p\] so large that the critical annulus
\[A_d(0)\] is a child of \[A_{d-p}(0)\], and
let \[\Delta_0=\widehat P_d(0)\] be the thickened
critical puzzle piece of depth \[d\]. Then \[\Delta_p=f^{\circ p}
(\Delta_0)\] is equal to \[\widehat P_{d-p}(c_p)=\widehat P_{d-p}(0)\],
which contains \[\Delta_0\]
in its interior. The hypothesis guarantees that the successive points
\[c_p\,,\;c_{2p}\,,\;c_{3p}\,,\;\ldots\] all belong to \[\Delta_0\].
Further details will be left to the reader.\QED\smallskip

{\bf Proof of Lemma 3} (suggested by M. Lyubich). Let \[\Delta_0=\widehat
P_1(0)\]. Then it is easy to check that the successive images\[\Delta_i
=f^{\circ i}(\Delta_0)\]
are disjoint from the critical point for \[0<i<q\], and that \[\Delta_q\]
contains \[\widehat P_0(0)\supset\Delta_0\] in its interior. %(So far, we have
%not made any use of the hypothesis.) 
We will prove inductively that
\[c_{qi}\in P_1(0)\subset\Delta_0\] for
every \[i\]. It certainly follows from this
inductive hypothesis that \[c_{qi+1}\in P_0(c_1)\cap K(f)\subset P_1(c_1)\],
and similarly that \[c_{qi+j}\in K(f)\cap P_0(c_j)\subset P_1(c_j)\]
for \[0<j<q\]. In particular, \[c_{qi+q-1}\] must belong to \[P_1(c_{q-1})\],
hence \[c_{qi+q}\in P_0(c_q)=P_0(0)\]. By the hypothesis of Lemma 3,
the point \[c_{qi+q}\] is known to lie in one of the
puzzle pieces \[P_1(c_h)\] which lie around \[\alpha\]. Evidently
it can only lie in \[P_1(0)\], as required.\QED\smallskip

{\QP{\bf Corollary.} \it If \[f\] is not renormalizable, then there exists a
critical annulus of positive modulus.\medskip}

{\bf Proof.} According to Lemma 3, the critical orbit must visit one of the
puzzle pieces \[P_1(-c_1)\,,\;\ldots\,,\;P_1(-c_{q-1})\] which surround
the pre fixed point \[-\alpha\]. Suppose for example that
\[c_d\in P_1(-c_i)\]. It is easy
to check that the corresponding annulus \[A_0(-c_i)\] has positive
modulus. (Compare Figure 2.) Pulling this annulus back inductively along
the critical orbit, it follows that \[A_d(0)\] also has positive modulus.
\QED\smallskip

Here is another application of thickened puzzle pieces. As usual, we
assume that the orbit of \[z_0\] does not hit the fixed point \[\alpha\].
\smallskip

{\QP{\bf Lemma 4.} \it Suppose that some orbit \[z_0\mapsto z_1\mapsto\cdots\]
in the Julia set never reaches the neighborhood \[P_N(0)\] of the critical
point. Then the intersection \[\bigcap P_d(z_0)\] of the puzzle pieces
containing \[z_0\] reduces to the single point \[z_0\].\medskip}

{\bf Proof.} Thickening the puzzle pieces
very slightly, we may assume that the orbit of \[z_0\]
never reaches \[\widehat P_N(0)\]. Let \[U_0\,,\;U_1\,,\;\ldots U_m\]
be the interiors of the various thickened puzzle pieces of depth \[N-1\],
numbered so that the critical value \[c_1=f(0)\] belongs to \[U_0\].
We will make use of the Poincar\'e metric on the \[U_i\]. For \[i>0\],
note that there are exactly two branches of \[f^{-1}\] on \[U_i\], call them
\[g_1\] and \[g_2\]. For \[i>0\], each of these maps \[g_k\] on \[U_i\], is a
holomorphic map
which carries \[U_i\] into a proper subset of some \[U_j\]. Hence it
strictly reduces the Poincar\'e metric. For each puzzle piece of depth \[N\]
contained in the \[U_i\], it follows that \[g_k\]
shrinks distances by some definite
factor \[\lambda<1\]. (We need the thickening at this point, to insure that
these puzzle pieces are compactly contained in \[U_i\].) Let \[C\] be the
maximum of the Poincar\'e
diameters of these puzzle pieces of depth \[N\]. Then for a puzzle piece
\[P_{N+h}(z_0)\], since the successive images \[P_{N+h-i}(z_i)\] for \[0\le i
\le h\] never meet the critical value region \[U_0\], it follows inductively
that the Poincar\'e diameter of \[P_{N+h}(z_0)\] is at most \[\lambda^h\,C\],
which tends to zero as \[h\to\infty\].\QED \smallskip

In order to deal with the possibility that \[f^{\circ d}(z_0)=\alpha\],
so that the puzzle piece \[P_d(z_0)\] is not uniquely defined,
we will need the following.\smallskip

{\bf Definition.} For any point \[z\] in the Julia set
let \[P^*_d(z)\] be the union of the puzzle pieces of depth \[d\]
containing \[z\]. (In most cases, \[P^*_d(z)\] is equal to the unique
depth \[d\] puzzle piece \[P_d(z)\] which contains \[z\]. However, if
\[f^{\circ d}(z)=\alpha\] then \[P^*_d(z)\] is a union of \[q\] distinct
puzzle pieces.)\smallskip

We can now state and prove the principal result of this
section.

{\QP{\bf Theorem 2.} \it Suppose as usual that \[f\] is quadratic, with
connected Julia set, with both fixed points repelling, and not renormalizable.
Suppose further that the critical orbit is disjoint from the fixed point
\[\alpha\]. Then for any\break \[z\in J(f)\] we have \[~\bigcap_d
P_d^*(z)=\{z\}\].\medskip}

{\bf Proof.} Since we assume that
\[f\] is not renormalizable, we know from the Corollary above that there
exists some critical annulus \[A_m(0)\] which has positive modulus. We will
first prove that \[\bigcap P_d(0)=\{0\}\], then prove that \[\bigcap P_d(z)
=\{z\}\] for any \[z\in J(f)\] which is not an iterated pre-image of
\[\alpha\], and finally prove the corresponding result when
\[f^{\circ n}(z)=\alpha\].
\smallskip

{\bf Critically Recurrent Case.} Suppose
that the critical orbit is recurrent, so that Lemma
1 applies. First consider
the puzzle pieces \[P_d(0)\] around the critical point. If the non-degenerate
annulus \[A_m(0)\]
has at least \[2^k\] descendents in the \[k$-th generation for each
\[k\], then each of these contributes exactly \[~\mod\,A_m(0)/2^k~\] to the sum
\[~\sum_d\mod\,A_d(0)~\]. Hence this sum is infinite, as required. On the other
hand, if there are fewer descendents in some generation, then one of them
must be an only child, hence excellent by Lemma 1d. Using Lemma 1b and 1c,
we again see that \[~\sum_d\mod\,A_d(0)~\] is infinite. Therefore the intersection \[\bigcap
\,P_d(0)\] reduces to the single point \[0\].
\smallskip

Now consider a point \[z_0\ne 0\] of the Julia set. We assume that the
orbit\break
\[z_0\mapsto z_1\mapsto\cdots\] is disjoint from \[\alpha\], so that
the puzzle pieces \[P_d(z_0)\] are well defined.
If the orbit of \[z_0\] does not accumulate
at zero, then we have \[\bigcap P_d(z_0)=\{z_0\}\] by Lemma 4. Suppose, on
the other hand,
that the origin is an accumulation point of \[\{z_n\}\].
In other words, suppose that the tableau of the point \[z_0\] has critical
annuli reaching down to all depths. For each depth \[d\], let us start at
column zero (corresponding to the point \[z_0\] itself) and advance to the
right until we first hit a critical annulus at column \[n\], then proceed
diagonal back
down until we reach column zero at depth \[n+d\]. It follows from this
construction that the annulus \[A_{n+d}(z_0)\] is conformally isomorphic
to \[A_d(0)\]. Furthermore, distinct values of \[d\] must correspond to
distinct
values of \[n+d\]. Thus the sum \[~\sum\mod\,A_d(z_0)~\] is also infinite,
hence \[~\bigcap P_d(z_0)=\{z_0\}\], as required.\smallskip

{\bf Critically Non-recurrent Case.} Now
suppose that the critical orbit is not recurrent. Then the critical
value \[f(0)\] satisfies the hypothesis of Lemma 4. Hence\break \[\bigcap
P_d(f(0))=\{f(0)\}\], and it follows easily that \[\bigcap P_d(0)=\{0\}\].
Next consider a point \[z_0\ne 0\]. Again we may assume that
the orbit \[\{z_n\}\] accumulates at zero, since otherwise the conclusion
would follow from Lemma 4. Again the Corollary above tells us that there
exists one depth
\[m\] such that \[A_m(0)\] has positive modulus. The corresponding
depth \[m\] for the tableau of \[z_0\] must have infinitely many columns
\[k\] which are critical. For each of these, let us proceed diagonally back
down in the tableau of \[z_0\], ignoring whatever critical or semi-critical
annuli we may meet, until we reach column zero at depth \[m+k\]. Each time
we meet a critical or semi-critical annulus, we lose up to half of the modulus.
(Problem 1-2 below.)
However, the hypothesis that the critical orbit is non-recurrent guarantees
that such losses will only occur a bounded number of times. Thus \[~\sum \mod\,
A_d(z_0)~\] has infinitely many summands which are bounded away from zero.
Hence this sum is infinite, and we have proved that
\[\bigcap P_d(z_0)=\{z_0\}\] in this case also.
%all cases. Together with Problem 1-1 below,
%this completes the proof of the Theorem.\QED
\smallskip

{\bf Iterated Pre-images of \[\alpha\].} If
some forward image \[z_n=f^{\circ n}(z_0)\] is equal to the
fixed point \[\alpha\], then the above arguments do not make sense as stated.
In this case, there are \[q\] distinct puzzle pieces \[P_n^{(i)}\] of depth
\[n\] which have \[z_0\] as common boundary point.
Each of these is contained in a
unique sequence of nested puzzle pieces \[~P_n^{(i)}
\supset P_{n+1}^{(i)}\supset\cdots~\] which have \[z_0\] as common boundary
point. {\bf Assertion:} {\it For
each one of these \[q\] nested sequences the intersection \[\bigcap_d
 P_d^{(i)}\] reduces to the single
point \[z_0\].} In fact, the proof of Lemma 4 applies equally well to this
situation. Evidently the statement that \[\bigcap_d
P^*_d(z_0)=\{z_0\}\] follows immediately. %\smallskip

Thus we have proved Theorem 2: \[\bigcap_dP^*_d(z)=\{z\}\] in all cases.
Theorem 1, as stated on page 2, is a
straightforward consequence.
(Compare Problem 1-1 below.) %\QED\hskip -.1in\QED
\rlap{$\sqcup$}$\sqcap$\hskip -.01in\QED
\centerline{--------------------------------}\smallskip

Here are some problems for the reader.\smallskip

{\bf Problem 1-1. Local Connectivity.} Prove that the intersection of
\[J(f)\] with each puzzle piece is
connected. Conclude that \[J(f)\] is locally
connected at \[z\] whenever \[\bigcap P^*_d(z)=\{z\}\].

{\bf Problem 1-2. Semi-Critical Annuli.} If \[A_d(z)\] is a non-degenerate
semi-critical annulus of depth \[d>0\], show that \[A_d(z)\] is the union of

(1) a ramified two-fold covering of \[A_{d-1}(f(z))\], and

(2) a conformal copy of \[P_d(f(z))\].

\noindent Using the Gr\"otzsch inequality, prove that
\[~~\mod \,A_d(z)\,>\,{\textstyle{1\over 2}}\,\mod\,A_{d-1}(f(z))~.\]

{\bf Problem 1-3. Non-degenerate Annuli.} Show that an annulus \[A_d(z_m)\]
is non-degenerate if and only if the corresponding annulus
\[A_0(z_{d+m})\] of depth zero is semi-critical.

{\bf Problem 1-4. Further Tableau Rules.} Let \[q\ge 2\] be
the number of external rays landing at the fixed point \[\alpha\].
Show that the semi-critical depth
of a tableau column can never take the values \[1\,,\,\ldots\,,\,q-1\].
Show that at most \[q-1\] consecutive
columns can be completely off-critical (semi-critical depth \[\scd(z_i)=-1\]),
and show that \[\scd(z_{i})=-1\] for \[m<i<m+q\]
if and only if the \[m$-th column has semi-critical depth \[\scd(z_m)\ge q\].

%If the semi-critical depth \[\scd(z_m)\]
%is \[\ge q\], show that the following \[q-1\] columns are completely
%off-critical: \[\scd(z_{m+i})=-1\] for \[0<i<q\].
%always satisfies either \[\scd(z_m)=-1\] or \[\scd(z_m)= 0\] or
%\[scd(z_m)\ge q\], where \[q\ge 2\] is the number of external rays landing
%
%at \[\alpha\]. If \[scd(z_m)\ge q\], show that
%$$	\scd(z_{m+1})\;=\;\cdots\;=\;\scd(z_{m+q-1})\;=\;-1~,\qquad\scd(
%z_{m+q})\ge 0~. $$
%On the other hand, if
%\[\scd(z_m)\le 0\], show that \[\scd(z_{m+i})\ge 0\] for some \[i\]
%with \[0<i<q\]. 

{\bf Problem 1-5. The Critical Orbit is Generically Dense.} It is convenient
to say that a property of certain points in a compact set is {\bit
generically\/} true if it is true throughout a countable intersection of
dense open subsets. For example, one can show that for generic \[c\] in the
boundary of the Mandelbrot set, the map \[f_c\] is
non-renormalizable\footnote{$^1$}{Proof outline: With notation as in \S3, the
set of
renormalizable points in the Mandelbrot set (together with associated root
points) forms a countable union of compact subsets \[H*M\]. This union is
nowhere dense in the boundary \[\partial M\]. In fact the set of \[c\] such
that the critical orbit of \[f_c\] eventually lands on the fixed point
\[\beta\] is everywhere dense in \[\partial M\]. (Compare [Br2] or [M1,
App. A].) Such a map \[f_c\] cannot be renormalizable [D3].},
with both fixed points repelling. Let \[U_d\subset
\partial M\] be the set of parameter values \[c\] in the boundary of the
Mandelbrot set such that every puzzle piece of depth \[d\] for \[f_c\]
contains a
post-critical point \[c_i=f_c^{\circ i}(0)\] in its interior. Show that \[U_d\]
contains a dense relatively
open subset of \[\partial M\]. (To prove density, use the fact
that periodic points are dense in \[J(f_c)\], and use Montel's Theorem.)
{\it For a generic parameter value \[c\in\partial M\], conclude that
the closure of the critical orbit is the entire Julia set \[J(f_c)\].\/}
Conclude also that no non-degenerate annulus can be excellent
in the generic case.

{\bf Problem 1-6. The Yoccoz \[\tau$-function.}
For each depth \[d>0\],
if there exists an integer \[0\le\tau<d\] so that \[f^{\circ d-\tau}\]
maps \[P_d(0)\] onto \[P_\tau(0)\], then we define \[\tau(d)\] to be the
%let \[\tau=\tau(d)\] be the largest integer \[\tau<d\]
%with \[f^{\circ \tau-d}(P_d(0))= P_{\tau}(0)\], let \[\tau(d)\] be the
largest such integer, and describe the critical puzzle piece \[P_{\tau(d)}(0)\]
as the ``parent'' of \[P_d(0)\]. Otherwise, if no such integer exists,
we set \[\tau(d)=-1\].
Thus \[-1\le\tau(d)<d\] in all cases. Show that \[\tau(d+1)\le \tau(d)+1\].
Show that the annulus \[A_d(0)\]
is a child of \[A_{n}(0)\] if and only if
% \[\tau(d)=n\ge 0\] and \[\tau(d+1)=n+1\].
\[~~n=\tau(d)=\tau(d+1)-1\ge 0\;\].

{\bf Problem 1-7. Persistent Recurrence.} 
The critical orbit is said to be {\bit persistently recurrent\/} if
it is non-renormalizable with \[\lim_{d\to\infty}\,\tau(d)=+\infty\].
Show that the Fibonacci tableau is persistently recurrent.
In the persistently recurrent case, show that the critical orbit is bounded
away from the fixed point \[\beta\]. (Otherwise, for any depth \[d\] we could
choose a post-critical point \[c_n\in P_d(\beta)\], then choose
the smallest \[k\]
with \[0\in P_{d+k}(c_{n-k})\], and conclude
that \[\tau(d+k)=0\].) Show that
the critical orbit is bounded away from every iterated pre-image of \[\beta\],
and hence that its closure is nowhere dense in the Julia set. More
generally, show that the critical orbit is bounded away from any periodic
point. Using the fact that it is bounded away from \[\alpha\], show that
the critical orbit closure is a Cantor set. (For further information,
see [L2].)\medskip
%\end

% written in plaintex   Standard TEX Header FILE

%\magnification=1200

\input psfig
\def\QP{\narrower\smallskip\noindent}
\def\ref{\hangindent=1pc \hangafter=1 \noindent}
\def\QED{  \rlap{$\sqcup$}$\sqcap$ \smallskip}
\def\mod{{\rm mod\;}}
\def\[{$\,}
\def\]{\,$}
\def\C{{\bf C}}
\def\R{{\bf R}}
\def\Z{{\bf Z}}

\def\={\;=\; }
\def\lln{\smallskip\centerline{------------------------------------------------------}\smallskip}
\font\bit=cmssi12 at 12truept
\font\tenmsy=msym10

\textfont8=\tenmsy
\mathchardef\ssm="7872
\mathsurround = 1pt
\abovedisplayskip=6pt
\belowdisplayskip=6pt
\parskip=2pt

\def\area{{\bf area}}
\def\mod{{\bf mod\,}}
\def\QP{\smallskip\leftskip=.4in\rightskip=.4in\noindent}
\def\oldQP{\narrower\noindent}
\parskip=4pt
%\pageno=17

\centerline{\bf \S2. Polynomials for which all but one
 of the critical orbits escape}

\centerline{\bf (following Branner and Hubbard).}\medskip

%\centerline{(Class Notes, J. Milnor)}\bigskip

Let \[f:\C\to\C\] be a polynomial of degree \[d\ge 2\], with filled Julia set
\[K\] and with Julia set \[J=\partial K\]. The object of this
section is to study the case where only one critical point has bounded orbit.
However, we begin with the simpler case where no critical point has
bounded orbit. The following is well known.

{\QP{\bf Theorem 3.} \it If all critical orbits of \[f\] escape to infinity,
then \[J=K\] is a Cantor set of measure zero\footnote{$^1$}
{\rm Shishikura has shown that the Hausdorff dimension of this Cantor set
can be arbitrarily close to two.}.
Furthermore,
the dynamical system \[~(J\,,\,f|J)\]\break is homeomorphic to the one-sided
shift on \[d\] symbols\footnote{$^2$}
{\rm Przytycki and Makienko have both announced the sharper result
that any rational Julia set which is totally disconnected and contains no
critical point must be isomorphic to a one-sided shift.}.\smallskip}

{\bf Proof.} Let \[\omega_1\,,\,\ldots\,,\,\omega_{d-1}\] be the (not
necessarily distinct) critical points of \[f\], and let
\[G:\C\to\R_+\] be the canonical potential function ($=$
Green's function), which satisfies \[G(f(z))=d\,G(z)\] and
vanishes precisely on the filled Julia set \[K\].
Thus \[G(\omega_i)>0\] by hypothesis. The critical points of \[G\] in
\[\C\ssm K\] are the points
\[\omega_i\] and also all of
their preimages under iterates of \[f\]. Hence the
critical values of \[G\] (other than
zero) are the numbers of the form \[G(\omega_i)/d^k\] with \[k\ge 0\].

The {\bit Branner-Hubbard puzzle\/} of \[f\] is constructed as follows.
Choose a number \[G_0\], not of the form \[G(\omega_i)/d^k\], so that
\[~0<G_0<{\rm Min}\{G(f(\omega_i))\}~.\]
Then the region \[G^{-1}(0\,,\,G_0]\] contains no critical values of \[f\],
and is bounded by smooth curves since \[G_0\] is not a critical value
of \[G\]. Similarly, each locus
$$	G^{-1}\Big[0\,,\,{G_0\over d^k}\Big]~=~\Big\{z\in\C~:~G(z)\le {G_0\over
d^k}\Big\}  \eqno (*_k) $$
is bounded by smooth curves.
Note that the complementary region \[G^{-1}(G_0/d^k\,,
\,\infty)\] cannot have any bounded component, since the harmonic function
\[G\] cannot have a local maximum. It follows that the locus \[(*_k)\]
is the disjoint union of a finite number of closed topological disks.
By definition, each of these closed disks
will be called a {\bit puzzle piece\/} \[P_k\] of depth \[k\]. Since these
puzzle pieces contain no critical value of \[f\], it follows easily that
\[f\] maps each \[P_k\] of depth \[k>0\] by a conformal isomorphism onto
some puzzle piece \[f(P_k)\] of depth \[k-1\].
%\vskip -.2in
\midinsert
\insertRaster 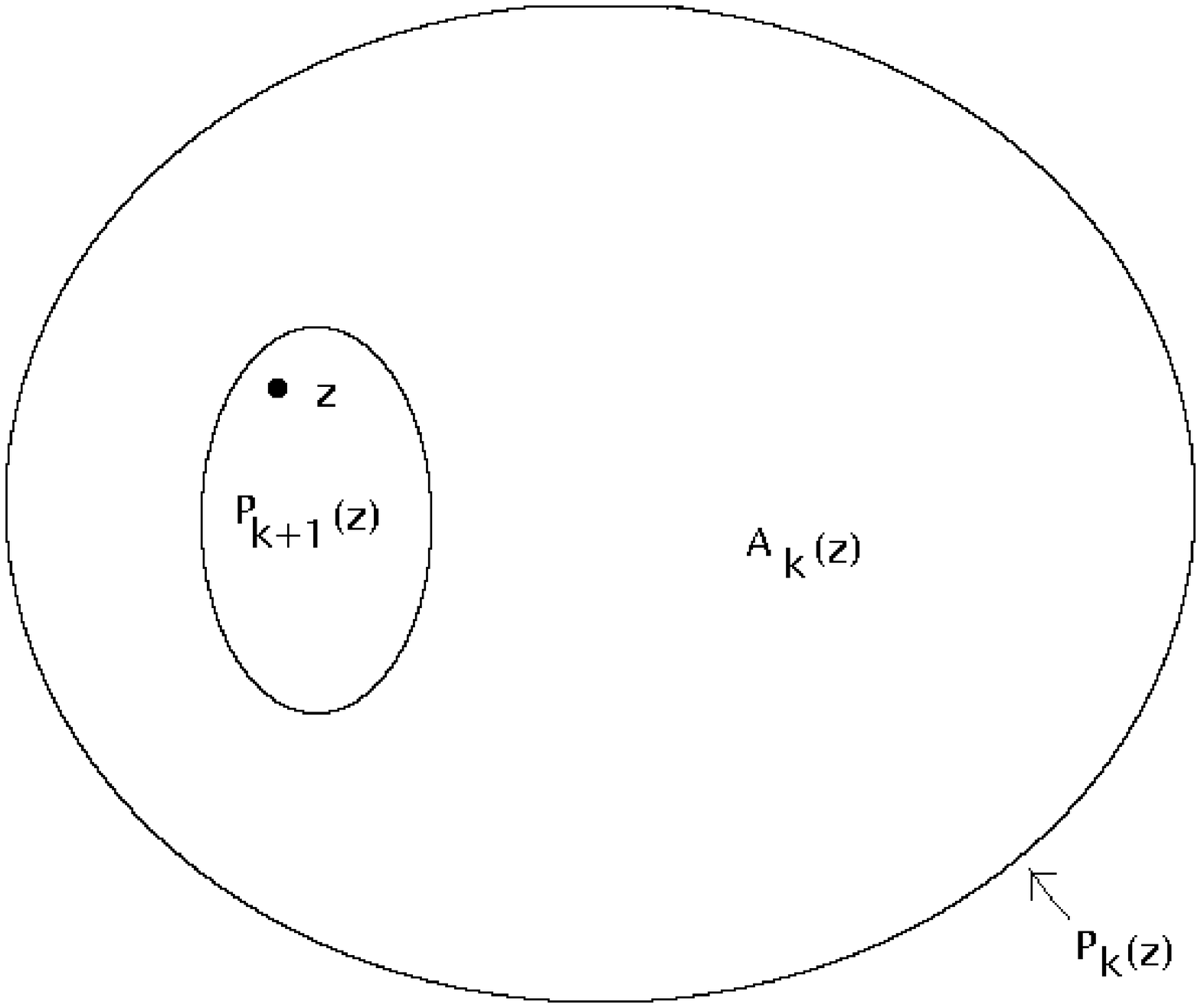 pixels 640 by 480 scaled 380
\centerline{\bit Figure 12. Nested puzzle pieces, and the annulus \[A_k(z)\].}
\endinsert

If \[P_k\supset P_{k+1}\] are nested puzzle pieces of depths \[k\] and \[k+1\],
then the set
$$	A_k~=~{\rm Interior}(P_k)\,\ssm\, P_{k+1}  $$
is a well defined annulus of strictly positive modulus. (Figure 12.)
We will call such an
\[A_k\] an {\bit annulus of depth\/} \[k\] in the Branner-Hubbard puzzle. 
Evidently \[f^{\circ k}\] maps
each annulus of depth \[k\] by a conformal isomorphism onto an annulus of
depth zero. {\it Hence the moduli of all annuli constructed in this way
are uniformly bounded away from zero.}
\vfil\eject

For any point \[z\] in the filled
Julia set we can form the nested sequence of\break
puzzle pieces \[P_0(z)\supset P_1(z)\supset\cdots\], all containing \[z\].
Since the moduli of the annuli \[({\rm Interior}\,P_k(z))\ssm P_{k+1}(z)\] are
bounded away from zero, it follows that the intersection \[\bigcap
P_k(z)\] reduces to the single point \[z\]. Since each boundary circle
\[\partial P_k\] is disjoint from \[K\], this implies that
\[K\] is totally disconnected.\smallskip

The proof that \[J\] has measure zero will be based on the {\bit
McMullen inequality}
$$	\area(P_{k+1})~\le~{\area(P_{k})\over 1+4\pi\,\mod(A_k)} ~, $$
with \[A_k={\rm Interior}(P_k)\ssm P_{k+1}\]
as above. (Compare the Appendix.)
However, to apply this inequality in a useful way, we
will need to construct annuli somewhat differently.
\vskip -.2in
\midinsert
\insertRaster 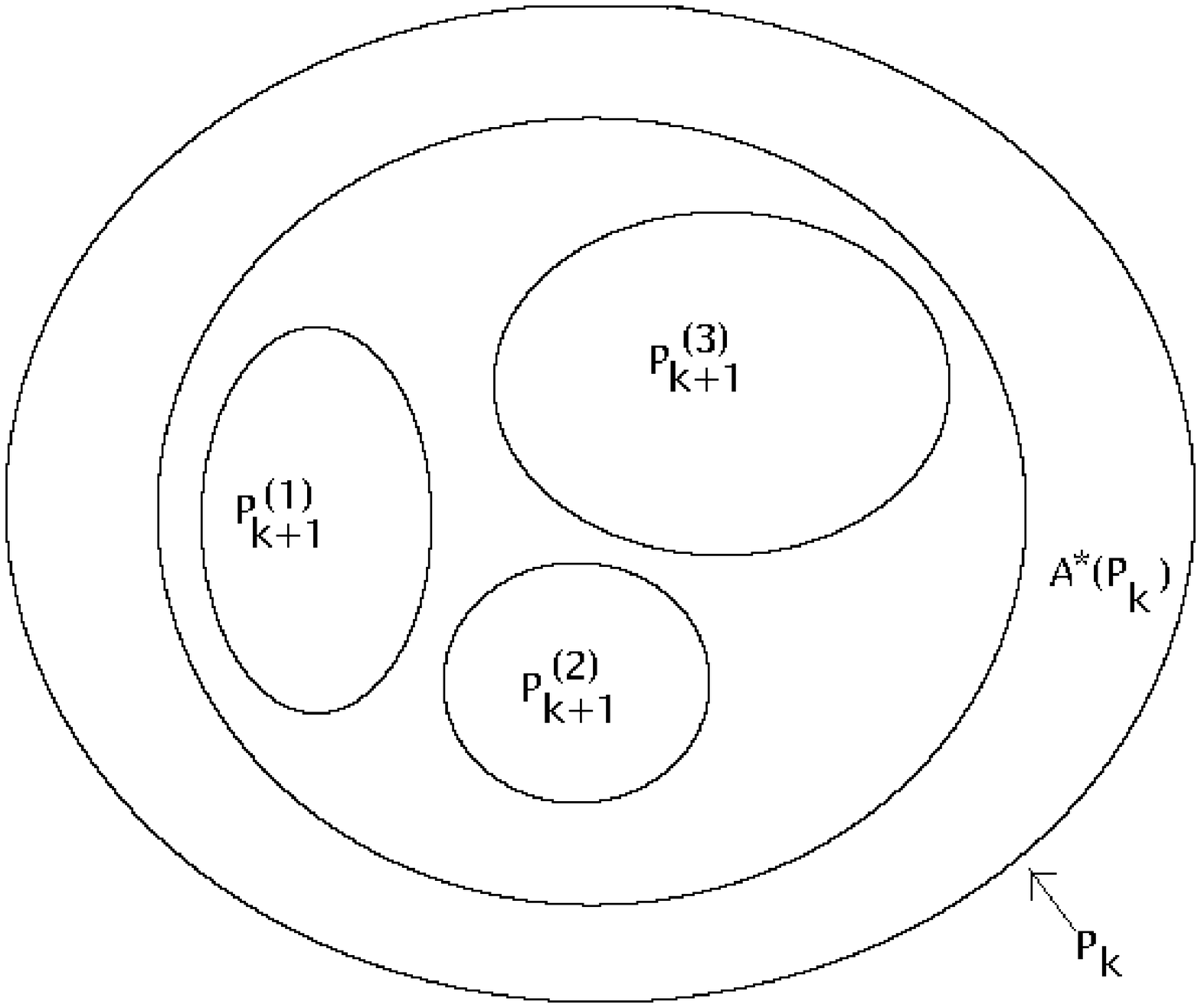 pixels 640 by 480 scaled 380
\centerline{\bit Figure 13. The ``thin annulus'' \[A^*(P_k)\subset A_k\subset
P_k\].}
%(Compare Figure 1.)}
\endinsert
\vfil\eject
Choose a number \[\epsilon>0\] which is small enough so that there are
no critical values of \[G\] within the closed interval
\[[G_0-\epsilon\,,\,G_0]\]. Then each connected component of\break
\[G^{-1}[(G_0-\epsilon)/d^k\,,\,G_0/d^k]\] is an annulus, and there is one
such ``thin annulus''
$$	A^*(P_k)~~=~~P_k\,\cap\, G^{-1}\Big[{G_0-\epsilon\over d^k}\,,\,{G_0\over d^k}
\Big] $$
within each connected component \[P_k\] of \[G^{-1}[0\,,\,G_0/d^k]\].
The construction is such that:

{\QP (a) all of these annuli \[A^*(P_k)\] have modulus
strictly bounded away from
zero, say \[~\mod A^*(P_k)\ge c>0\,\],~ and\smallskip

\noindent
(b) every puzzle piece \[P_{k+1}\] which is contained in \[P_k\] is actually
contained in the smaller disk \[~P_k\ssm A^*(P_k)~\]
(the bounded component of \[~\C\ssm A^*(P_k)~\]).\smallskip}

\noindent 
Thus the McMullen inequality takes the following sharper form. For each fixed
puzzle piece \[P_k\],
%$$	\sum\{\area(P_{k+1}) : P_{k+1}\subset P_k\}
$$	\sum_{P_{k+1}\subset P_k}\area(P_{k+1})
~~\le~~{\area(P_k)\over
1+4\pi\,\mod(A^*(P_k))}~~\le ~~{\area(P_k)\over 1 +4\pi\,c}~, $$
to be summed over those puzzle pieces of depth \[k+1\] which are contained
in the given \[P_k\]. It now follows inductively
that the sum of the areas of {\it all\/}
puzzle pieces of depth \[k\] satisfies
$$	\sum\area(P_k)~\le~\sum\area(P_0)/(1+4\pi\,c)^k~. $$
Since this tends to zero as \[k\to\infty\], and since \[J\subset\bigcup P_k\],
it follows that \[J\] has area zero.\smallskip

To prove that \[(J\,,\,f|J)\] is isomorphic to the one-sided \[d$-shift, we
proceed as follows. We will construct a closed subset \[\Gamma\subset \C\],
consisting of paths leading out to infinity, such that:

(1) \[\Gamma\] contains all critical values of \[f\],

(2) \[f(\Gamma)\subset\Gamma\], and

(3) the complement \[U=\C\ssm\Gamma\] is simply-connected, and contains the
Julia set.

\noindent
In fact, starting from each critical value \[f(\omega_i)\], we can follow the
gradient directions of \[G\] (the path of steepest ascent) until we meet a
critical point of \[G\], or all the way out to infinity if we never meet a
critical point. At each critical point of \[G\] we must make a choice.
From a critical point of multiplicity \[\mu\], there are \[\mu+1\] distinct
directions in which we can continue along a path of steepest ascent. Choose
one of these \[\mu+1\] directions for each critical point of \[G\].
However, the choice must be consistent
in the following sense: If \[f(z)=z'\], where \[z\] and \[z'\]
are both critical points of \[G\], then \[f\] must carry the preferred path
from \[z\] to the preferred path from \[z'\]. It is not difficult to make
such a consistent choice, simply by ordering the critical points of \[f\]
so that \[G(\omega_1)\le\cdots\le G(\omega_{d-1})\], and making a choice
of ascending path first at \[\omega_1\], and its iterated pre-images, then at
\[\omega_2\] and its iterated pre-images, and so on. Now we can follow
the chosen paths from all of
the post-critical points \[f^{\circ k}(\omega_i)\] out to infinity. The paths
from different post-critical points may come together, but once they come
together they must stay together, out to infinity. The union \[\Gamma\]
of these preferred paths will be a locally finite union of disjoint
topological trees with the required properties.

Since \[f(\Gamma)\subset\Gamma\], it follows that \[f^{-1}(U)\subset U\].
Furthermore, since \[U\] is simply connected, and contains no critical
values of \[f\], it follows that every one of the \[d\] branches of \[f^{-1}\]
near a point of \[U\] extends uniquely to a holomorphic map
\[	g_i:U\to U\]. The images
$$	g_1(U)\,,\,\ldots\,,\,g_d(U)~\subset~ U $$
are disjoint open sets covering the Julia set. We will show that
the intersections
$$	J_i~=~J\cap g_i(U) $$
are disjoint compact sets which form the
required {\bit Bernoulli partition\/} \[J=J_1\cup\cdots\cup J_d\] of the
Julia set. That is:

{\QP\it For each sequence of integers \[i_0\,,\,i_1\,,\,\ldots\] between
\[1\] and \[d\] there exists one and only one orbit \[z_0\mapsto z_1\mapsto
z_2\mapsto\cdots\] in the Julia set with \[z_n\in J_{i_n}\] for every \[n\].
\smallskip}

\noindent
In fact \[z_0\] can be described as the intersection of the nested sequence
of sets\break
\[g_{i_0}\circ g_{i_1}\circ\cdots\circ g_{i_n}(J)\]. To prove this
statement, we use the Poincar\'e metric on \[U\]. Since each \[g_i\]
restricted to the compact set \[J\subset U\] shrinks Poincar\'e distances
by a factor bounded away from one, it follows that each such intersection
\[\bigcap_n g_{i_0}\circ g_{i_1}\circ\cdots\circ g_{i_n}(J)\]
consists of a single point. This proves Theorem 3.\QED\bigskip%\bigskip
%\eject

\centerline{\bf Maps with exactly one bounded critical orbit.}\smallskip

This will be an exposition of results due to Branner and Hubbard [BH3]. We
now suppose that exactly one of the \[d-1\] critical points of \[f\] has
bounded orbit, while the orbits of the remaining \[d-2\]
critical points, counted with multiplicity, escape to infinity. (Thus
we exclude examples such as \[z\mapsto z^3+i\], for which
a double critical point has bounded orbit; however, a double critical
point escaping to infinity is fine.)
Furthermore, we assume that
\[d\ge 3\], so that at least one critical orbit does escape. Then, according
to Fatou and Julia, the Julia set is disconnected, with uncountably
many connected components. Let
$$	c_0~\mapsto~ c_1~\mapsto~ c_2~\mapsto~\cdots $$
be the unique bounded critical orbit.

As in the proof of Theorem 3, choose a number \[G_0>0\] which is not a critical
value of \[G\], and so that the region \[G^{-1}(0\,,\,G_0]\] contains no
critical value of \[f\]. Again we define the puzzle pieces of
depth \[k\] to be the connected components \[P_k\] of the locus
$$	\bigcup P_k~=~G^{-1}\Big[0\,,\,{G_0\over d^k}\Big]~. $$
Thus
each point \[z\in K\] determines a nested sequence
\[P_0(z)\supset P_1(z)\supset\cdots\], and the central
problem is to decide whether or not \[\bigcap_k P_k(z)=\{z\}\]. Again
we look at the intermediate annuli
$$	A_k(z)~=~{\rm Interior}\,P_k(z)\ssm P_{k+1}(z)~. $$
As in the Yoccoz proof, such an annulus is said to be
semi-critical, critical, or off-critical
according as the critical point \[c_0\] belongs to the annulus itself,
or to the bounded or the unbounded component of its complement. (For this
purpose, we ignore the other critical points, whose orbits escape to
infinity.) This Branner-Hubbard puzzle is easier to deal
with than the Yoccoz puzzle for three reasons:

(a) All of the annuli \[A_k(z)\] are non-degenerate, with strictly positive
modulus.

(b) The various puzzle pieces of depth \[k\] are
pair-wise disjoint.

{\oldQP (c) For each \[z\in K\], the intersection \[\bigcap_k P_k(z)\] of the puzzle
pieces containing \[z\] has an immediate
topological description: It is equal to the connected component of the filled
Julia set \[K\] which contains the given point \[z\]. For this intersection
is clearly connected, and no larger subset can be connected since each boundary
\[\partial P_k(z)\] is disjoint from \[K\].\smallskip}

\noindent The {\bit
tableau\/} of an orbit \[z_0\mapsto z_1\mapsto\cdots\] can be described as
a record of
exactly which of the annuli \[A_k(z_i)\] are critical or semi-critical or
off-critical. 
First suppose that the tableau of the critical orbit is not periodic. %, or
%equivalently suppose that the post-critical points \[c_1\,,\,c_2\,,\,\dots\]
%are all disjoint from the critical component \[\bigcap_k P_k(c_0)\].
%. That is,
%suppose that for every \[i>0\] there exists a puzzle piece \[P_k(c_i)\] 
%which is distinct from the critical puzzle piece \[P_k(c_0)\]. 
We continue to
assume that \[d\ge 3\] and that exactly one of the \[d-1\] critical points
has bounded orbit.

{\QP{\bf Theorem 4.} \it If the critical tableau is not periodic, or
equivalently if the post-critical points \[c_1\,,\,c_2\,,\,\dots\]
are all disjoint from the critical component \[\bigcap_k P_k(c_0)\], then
for every point \[z_0\] of the filled Julia set \[K\] the sum\break \[~
\sum_k\mod A_k(z_0)~\] is infinite, hence \[\bigcap_k P_k(z_0)=
\{z_0\}\]. It follows that\break \[J=K\]
is a totally disconnected set of area zero.\smallskip}

\midinsert
\centerline{\psfig{figure=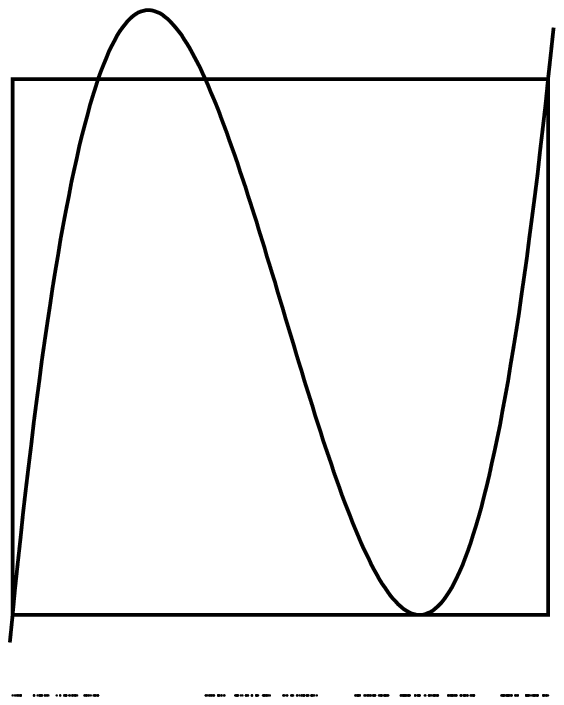,height=3in}}
\vskip -.1in
{\QP\bit Figure 14. An example:
Graph of the function \[f(x)=x(x-c_0)^2/(1-c_0)^2\] on the unit
interval \[I=[0,1]\], with \[c_0=.76\]. There is just
one bounded critical orbit \[c_0\mapsto 0\mapsto
0\mapsto\cdots\]. Since the iterated pre-images of any point of \[J\] are
dense in \[J\], and since
there are three distinct branches of \[f^{-1}\]
mapping \[I\] into itself, it follows that the Julia set \[J\]
is completely contained in the real interval \[I\].
This Julia set is  plotted underneath the graph. Since \[c_0\]
and \[0\] evidently belong to distinct connected components of \[J\], it
follows from Theorem 4 that \[J\] is totally disconnected.\par}\vskip -.1in
\endinsert
%\eject

Thus \[J\] is again homeomorphic to a Cantor set. However, in this case
\[(J\,,\,f|J)\]\break is not isomorphic to a shift, or even a sub-shift.
For there are critical points in \[J\], hence \[f\]
is not locally one-to-one on \[J\].\smallskip

The proof of this theorem is quite similar to the proof of the Yoccoz theorem.
However there are simplifications, leading to the sharper result which is
stated: \[\sum \mod\,A_k(z)=\infty\] for {\it all\/} \[z\in K\].
The main difference is the following statement, which is true whether or
not the critical tableau is periodic. Consider an
orbit \[z_0\mapsto z_1\mapsto\cdots\].

{\QP{\bf Lemma 5.} \it Suppose that the points \[z_1\,,\,z_2\,,\,\ldots\] are
all disjoint from some neighborhood \[P_N(c_0)\]
of the critical point \[c_0\]. Then the annuli
\[A_k(z_0)\] have modulus uniformly
bounded away from zero, hence \[\sum_k\mod A_k(z_0)=
\infty\].\smallskip}

{\bf Proof.} Each annulus \[A_k(z_i)\] of depth \[k>0\] has modulus at least
half of the modulus of \[A_{k-1}(z_{i+1})\]. In fact, if \[i>0\] and \[k>N\]
then these two annuli are conformally isomorphic. Thus
$$	\mod A_k(z_0)~\ge \mod A_0(z_k)/2^{N+1}~, $$
where the right side is bounded away from zero since there are only finitely
many annuli of depth zero.\QED
\vfil\eject

{\bf Proof of Theorem 4.} If the critical tableau is
recurrent, then the proof proceeds exactly as is the Yoccoz argument: Every
critical annulus either has at least \[2^n\] descendents in the \[n$-th
generation for every \[n\], or else has a descendent with this property.
Since all annuli are non-degenerate, it follows that
\[\sum\mod A_k(c_0)=\infty\]. On the other hand, if the critical tableau is
non-recurrent, then it follows from Lemma 5 that \[\sum\mod
A_k(c_0)=\infty\]. In particular, it follows that the collection of
puzzle pieces
\[\{P_k(c_0)\}\] forms a fundamental system of neighborhoods of the critical
point \[c_0\]. The corresponding statements for any iterated pre-image
of \[c_0\] follow immediately.

Now consider a point \[z_0\in K\], with orbit \[z_0\mapsto
z_1\mapsto\cdots\] which never meets the critical point \[c_0\].
If this orbit does not accumulate at \[c_0\], then the
statement that \[\sum_k\mod A_k(z_0)=\infty\] follows from Lemma 5.
On the other hand, if this orbit does accumulate at \[c_0\], then since
\[\sum\mod A_k(c_0)=\infty\], an easy tableau argument shows that
\[\sum_k\mod A_k(z_0)=\infty\] also.
Thus \[\bigcap_k P_k(z)=\{z\}\] in all cases.

Since each boundary \[\partial P_k(z)\] is disjoint from the filled Julia
set \[K\], it follows that \[K\] is totally disconnected, and hence that
\[J=K\]. To prove that this set has measure zero, we proceed as in the proof
of Theorem 3. Choose \[\epsilon>0\] so that the interval \[[G_0-\epsilon\,,\,
G_0]\] contains no critical values of the function \[G:\C\to\R_+\]. Then
each puzzle piece \[P_k(z)\] contains a unique component \[A^*(P_k(z))\] of the
set \[[G^{-1}[(G_0-\epsilon)/d^k\,,\,G_0/d^k]\]. These annuli \[A^*(P_k(z))
\subset
A_k(z)\] are also non-degenerate, and the proof above shows equally well that
\[\sum_k \mod A^*(P_k(z))=\infty\]. (The proof actually becomes a little
easier, since a
thin annulus can never contain the critical point \[c_0\].)
Hence, just as in the proof of Theorem 3,
for each fixed puzzle piece \[P_k\], the total
area of the puzzle pieces \[P_{k+1}\] of depth \[k+1\] which are contained in
\[P_k\] satisfies
$$	\sum_{P_{k+1}\subset P_k} \area\, P_{k+1}~~\le~~{\area\, P_k\over
	1+4\pi \,\mod A^*(P_k)}~.$$
Define the ratio \[\mu(P_k)\] by the formula
$$	\area\, P_k~~=~~\mu(P_k)\sum_{P_{k+1}\subset P_k} \area\, P_{k+1} $$
so that this McMullen inequality takes the form \[~ 1+4\pi\,\mod
A^*(P_k)~\le~\mu(P_k)\,\]. Then, %summing over {\it all\/} puzzle pieces of
%given depth, and 
substituting this formula inductively, we can write
$$	\sum_{ P_0}\area\, P_0~=~\sum_{P_0\supset P_1} \mu(P_0)\,\area(P_1)
	~=~\cdots~=~\sum_{P_0\supset\cdots\supset P_k}\mu(P_0)\,\cdots\,
\mu(P_{k-1})\,\area(P_k)~,$$
where the left hand expression is to be summed over all puzzle pieces of depth
zero, the next over all pairs \[P_0\supset P_1\], and so on. (If \[P_0\supset
\cdots\supset P_k\], note that \[P_0\,,\,\ldots\,,\, P_{k-1}\]\break
are uniquely determined by \[P_k\].)
Let \[\eta_k\] be the minimum value
of the product\break \[\mu(P_0)\cdots \mu(P_{k-1})\]
as \[P_0\,,\,\ldots\,,\,P_{k-1}\] varies over all sequences of nested
puzzle pieces \[P_0\supset P_1\supset\cdots\supset P_{k-1}\]. Then we see
from this last equality that
$$	\sum\area(P_0)~~\ge~~\eta_k\sum\area(P_k)~, $$
to be summed over all puzzle pieces of depth zero or \[k\] respectively.
Thus, if we can prove that \[\eta_k\to\infty\] as \[k\to\infty\], then it will
follow that
$$	\area(J)~\le~\sum\area(P_k)~\le ~\sum\area(P_0)/\eta_k~\to 0~, $$
hence \[\area(J)=0\] as required.

Clearly \[1<\eta_1<\eta_2<\cdots\]. If these numbers tended to a finite limit
\[L<\infty\], then for each \[k\] we could find puzzle pieces
\[P_0(k)\supset P_1(k)\supset\cdots\supset P_{k-1}(k)\] so that
\[\mu(P_0(k))\cdots \mu(P_{k-1}(k))\le L\]. Hence we could
choose a puzzle piece \[P_0\]
which occurs infinity often as \[P_0(k)\], then choose \[P_1\subset P_0\]
which occurs infinitely often as \[P_1(k)\], and so on. In this way, we could
find a sequence
$$	P_0~\supset~ P_1~\supset~P_2~\supset~\cdots $$
with \[\mu(P_0)\cdots \mu(P_{k-1})\le L<\infty\] for every \[k\]. Since
\[~ 1+4\pi\, \mod(A^*(P_i))\;\le\;\mu(P_i)\,\], this would imply that
$$	1+4\pi\sum\mod(A^*(P_i))~\le~ L<\infty~, $$
contradicting our statement that \[\sum\mod A_k(z)=\infty\] for all
\[z\in K\]. This completes the proof of Theorem 4.\QED
\bigskip

On the other hand, if the critical tableau is periodic, then we will prove the
following.
\vfil\eject

\pageinsert
\centerline{\psfig{figure=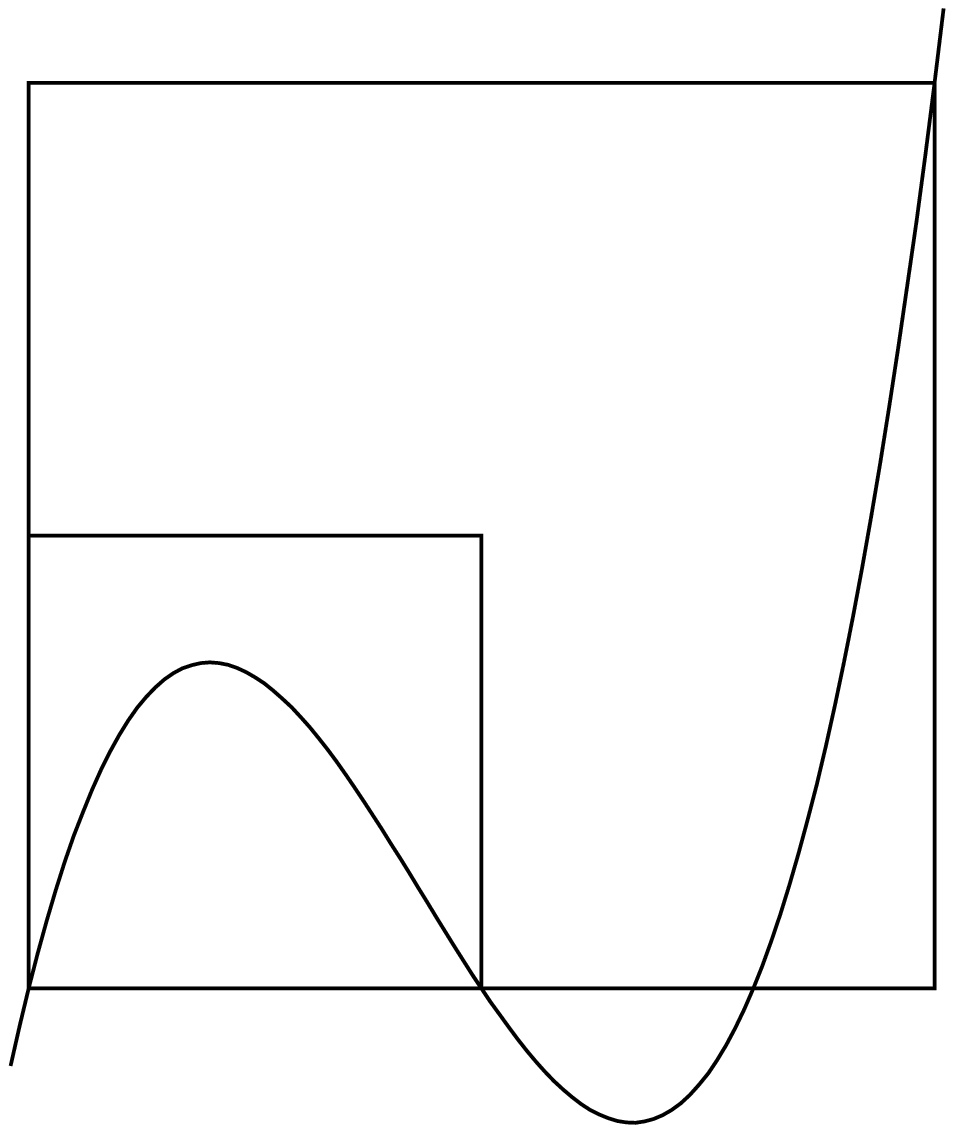,height=5in}}
\vskip -.2in
{\QP\bit Figure 15. Example for Theorem 5: Graph of the map \[f(x)=x(2x-1)
(5x-4)\] on the unit interval. In this case the connected interval
\[[0,{1\over
2}]\] is contained in the filled Julia set \[K\]. The orbit
of the critical point \[{2\over 3}\] escapes to \[-\infty\].\bigskip}
\insertRaster 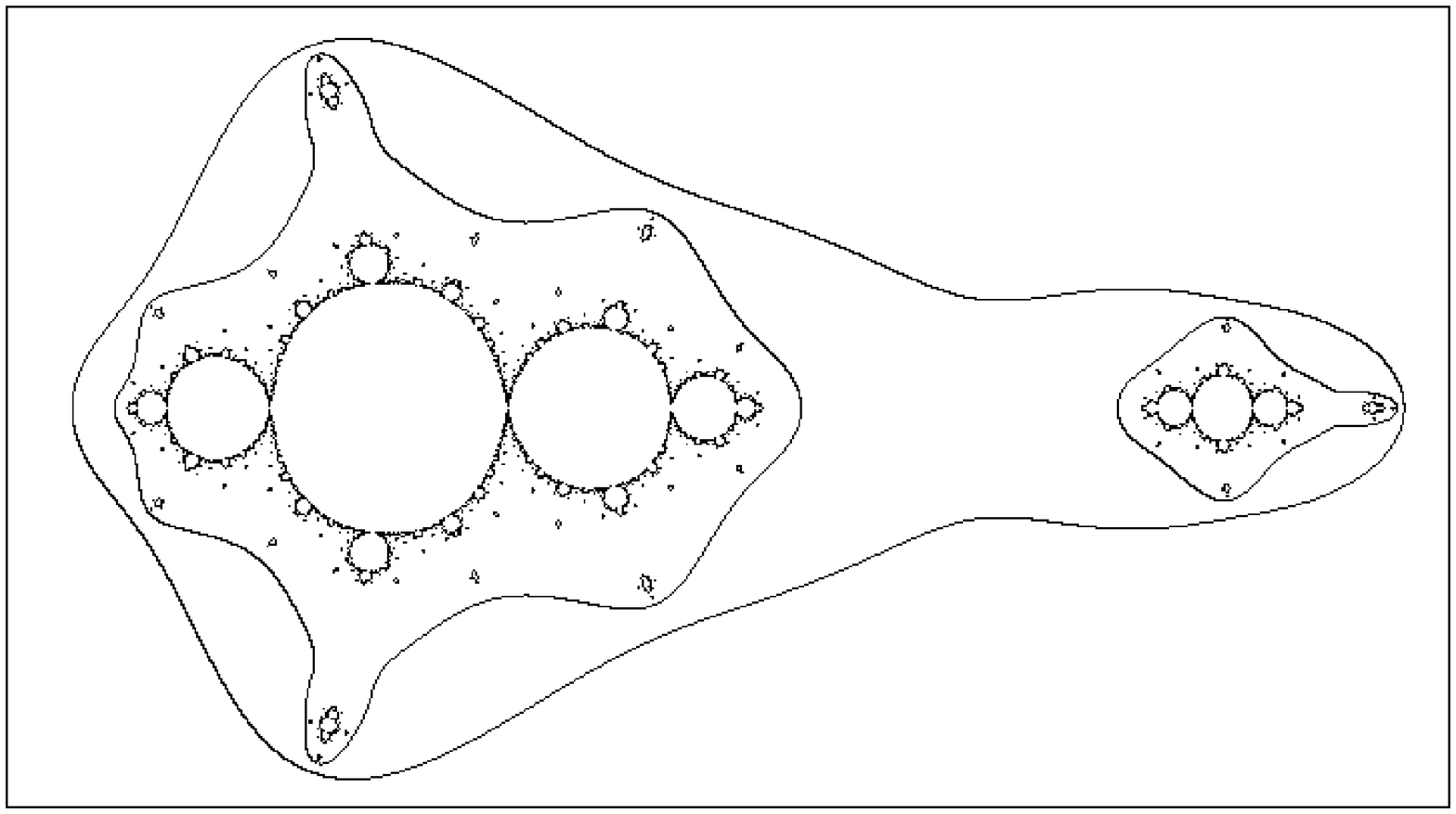 pixels 720 by 400 scaled 450
{\QP\bit Figure 16. Julia set for this map, drawn to the same scale,
showing
the puzzle pieces of level zero and one. Each non-trivial component of \[J\]
is homeomorphic to a certain quadratic Julia set.\par}
\endinsert

{\QP{\bf Theorem 5.} \it Still assuming that just one critical orbit
\[c_0\mapsto c_1\mapsto\cdots\] is
bounded, if the critical tableau is periodic of period \[p\ge 1\], so that\break
\[P_k(c_0)=P_k(c_p)\] for all depths \[k\], then the connected component of
the filled Julia set \[K=K(f)\] which contains \[c_0\] is non-trivial, that
is, consists of more than one point. In fact, a component of \[K\] is
non-trivial if and only if it contains some iterated pre-image of \[c_0\].
\smallskip}

\midinsert
\insertRaster 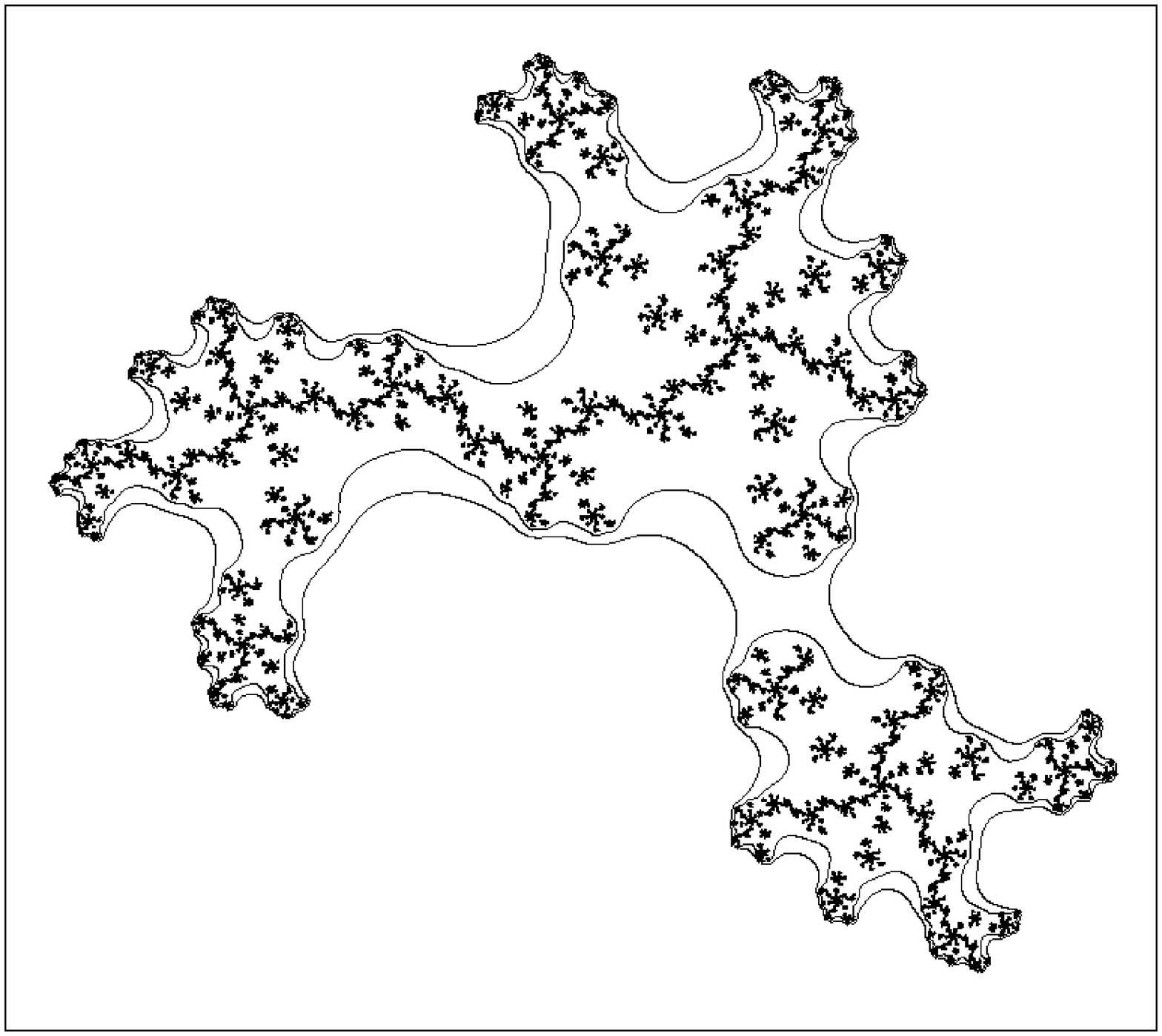 pixels 768 by 684 scaled 450
{\QP\bit Figure 17. Julia set for \[z\mapsto z^3+a\,z^2+1\],
with \[a=-1.10692+.63601\,i\], showing the puzzle pieces of level zero and one.
Each non-trivial component of \[J\] is homeomorphic to the
Julia set for the quadratic map \[z\mapsto z^2+i\].\par}
\endinsert

Thus there there are countably many non-trivial components
of \[K\]. These countably many components are everywhere dense in \[K\],
since the iterated pre-images of any point of the Julia set are dense in
the Julia set. Since a disconnected Julia set necessarily
has uncountably many components, it follows that there are uncountably many
single points components.
(The Julia set of a rational function may have uncountably
many non-trivial components. Compare [McM]. However, in the polynomial
case no such example is known.)
\vfil\eject

{\bf Proof of Theorem 5.}
%In all cases, note that the connected component of \[K\] containing a given
%point \[z\] coincides with the intersection of the puzzle pieces \[\bigcap_k
%P_k(z)\]. This intersection is certainly a connected subset of \[K\],
%and no larger subset
%of \[K\] can be connected since each \[\partial P_k(z)\] is disjoint from
%\[K\].
If \[~P_k(c_0)=P_k(c_k)~\] for all \[k\], then it follows from the tableau
rules that \[~P_k(c_i)=P_k(c_{p+i})~\]
for all \[i\] and \[k\]. Hence the entire orbit
$$	c_0~\mapsto~c_p~\mapsto c_{2p}~\mapsto~\cdots $$
of \[c_0\] under \[f^{\circ p}\] is contained in the critical component
\[\bigcap_k P_k(c_0)\subset K(f)\]. This\break intersection certainly
has more than one point. For either it contains \[c_0\ne c_p\],
or else \[c_0\] is a superattracting point, in which
case some entire neighborhood of \[c_0\] belongs to \[\bigcap_k P_k(c_0)\].

It follows easily that every pre-critical point in \[K(f)\] also belongs to
a non-trivial connected component. For if the orbit \[z_0\mapsto z_1\mapsto
\cdots\] intersects the critical component \[\bigcap_k P_k(c_0)\],
then we can choose the smallest \[\ell\ge 0\] for which \[z_\ell\] belongs to
this critical component. It follows easily that \[f^{
\circ\ell}\] maps the component \[\bigcap_k P_k(z_0)\] containing \[z_0\]
homeomorphically onto this critical component.

Now consider an orbit \[z_0\mapsto z_1\mapsto\cdots\] in \[K(f)\] which is
disjoint from this critical component. This means that the tableau of this
orbit has no columns which are completely critical. The proof is now divided
into two cases:

{\bf Case 1.} Suppose that there exists a fixed puzzle piece \[P_k(c_0)\]
which is disjoint from this orbit \[\{z_n\}\]. Then according to Lemma 5
we have \[\sum_k\mod A_k(z_0)=\infty\], hence \[\bigcap_k P_k(z_0)=\{z_0\}\].

\midinsert
\centerline{\psfig{figure=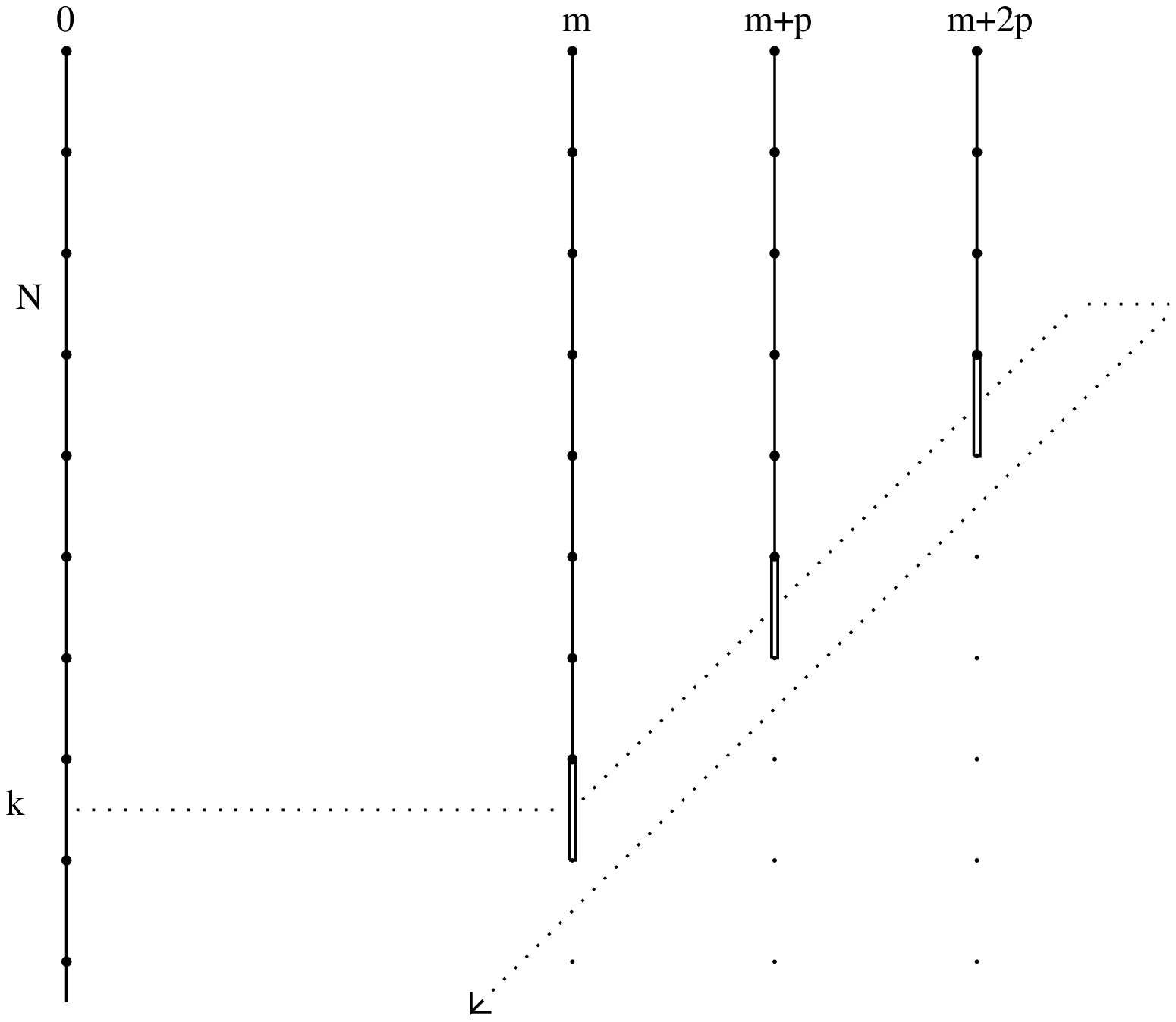,height=3in}}
\centerline{\bit Figure 18. Tableau for an orbit \[z_0\mapsto z_1\mapsto\cdots\]
which intersects every critical puzzle piece.}
\endinsert

{\bf Case 2.} If the orbit \[\{z_n\}\] intersects {\it every\/}
critical puzzle piece,
then we use a tableau argument as follows. Choose a depth \[N\] so that the
periodic critical tableau has no semi-critical annuli at depths \[\ge N\].
By hypothesis, there are infinitely many pairs \[(k,m)\], with \[k\ge N\],
so that the \[k$-th
row of the tableau for \[z_0\] is semi-critical in column \[m\] and
off-critical in earlier columns. (Compare Figure 18.)
Using the tableau rules, we can then compute the tableau in column \[m+i\]
and depth \[k-i\] for \[0<i<k-N\]. In fact the entries in column \[m+jp\]
and depth \[k-jp\] are semi-critical, and the others are off-critical.
It now follows that the annulus \[A_{k+m+1}(z_0)\] is conformally isomorphic to
an annulus of depth \[N\]. Hence this modulus is bounded away from zero.
It follows that \[\sum_\ell \mod A_\ell(z_0)=\infty\], which completes the
proof of Theorem 5.
\QED
\bigskip\bigskip
We can understand this proof better by introducing the following concepts,
which are due to Douady and Hubbard [DH3].

{\bf Definition.} By a {\bit polynomial-like map~}
is meant a pair \[(g\,,\,\Delta)\] where \[\Delta\subset\C\] is a closed
topological disk and \[g\] is a continuous mapping, holomorphic on the
interior of \[\Delta\],
which carries \[\Delta\] onto a closed topological
disk \[g(\Delta)\] which contains \[\Delta\] in its interior, such that
\[g\] maps boundary points of \[\Delta\] to boundary points of \[g(\Delta)\].
%\vskip -.1in

\midinsert
\insertRaster 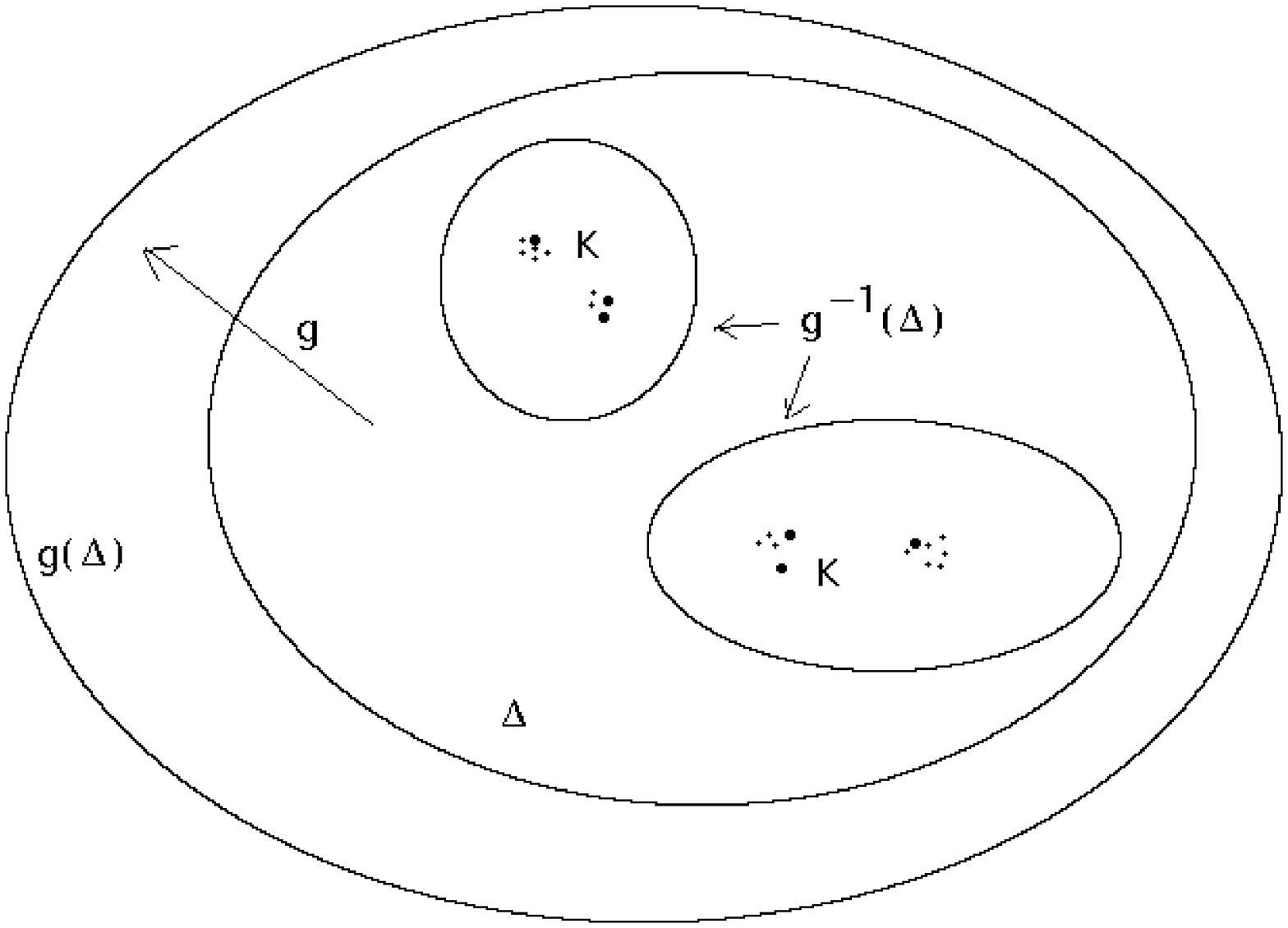 pixels 640 by 480 scaled 450
\centerline{\bit Figure 19. A polynomial-like mapping \[(g,\Delta)\]
with \[K=K(g,\Delta)\] totally disconnected.}
\endinsert

The {\bit degree\/} \[d\ge 1\] of such a polynomial-like mapping is a well
defined topological invariant. Note that almost every point of \[g(\Delta)\]
has precisely \[d\] pre-images in \[\Delta\].
The {\bit filled Julia set\/} \[K(g,\Delta)\] is defined
to be the compact set consisting of all \[z_0\in\Delta\] such that the entire
orbit \[z_0\mapsto z_1\mapsto\cdots\] of \[z_0\] under \[g\] is defined
and is contained in \[\Delta\].
%\eject

{\QP{\bf Lemma 6.} \it Such a polynomial-like map of degree \[d\] has \[d-1\]
critical points, counted with multiplicity, in the interior of \[\Delta\].
The filled Julia set
\[K(g,\Delta)\] is connected if and only if it contains all of these \[d-1\]
critical points.\smallskip}

{\bf Proof.} Consider the nested sequence of compact sets
$$	g(\Delta)~\supset~\Delta~\supset~ g^{-1}(\Delta)~\supset~
 g^{-2}(\Delta)	~\supset~\cdots $$
with intersection \[K(g,\Delta)\]. First suppose that the boundary
\[\partial\Delta\] contains no post-critical points, that is points
\[g^{\circ n}(\omega_i)\] with \[n>0\] where \[\omega_i\] is a critical point
of \[g\].
Then clearly each \[g^{-n}(\Delta)\] is a compact set bounded by one or more
closed curves. In fact, each \[g^{-n}(\Delta)\] is either a
closed topological disk or a finite union of closed topological disks,
each of which maps onto the entire disk \[\Delta\] under \[g^{\circ n}\].
To see this, note that for each component \[B\] of \[g^{-n}(\Delta)\] the
image \[g^{\circ n}(B)\] is compact, and that \[g^{\circ n}\] maps
the boundary \[\partial B\] into
\[\partial \Delta\] and maps the interior of \[B\] onto an open subset
of the interior of \[\Delta\]. Since the interior of \[\Delta\] is connected,
this implies that \[g^{\circ n}(B)=\Delta\]. If some such component \[B\] were
not simply-connected, then some component \[B'\] of \[\C\ssm
g^{-n}(\Delta)\] would be bounded. A similar argument would then show that
\[g^{\circ n}\] must map
\[B'\] onto the complementary disk \[\hat\C\ssm {\rm
Interior}(\Delta)\], which is impossible.

Applying the Riemann-Hurwitz formula
to the ramified covering \[\Delta\to g(\Delta)\], we see that
the number of critical points of \[g\] in the interior of \[\Delta\],
counted with multiplicity, is equal to
\[d\,\chi(g(\Delta))-\chi(\Delta)=d\cdot 1-1\]. (Here \[\chi\] is
the Euler characteristic.)
 Similarly, applying this formula to
\[g^{-1}(\Delta)\to \Delta\], we see that the number of critical points
in \[g^{-1}(\Delta)\] is equal to \[d\,\chi(g^{-1}(\Delta)-\chi(\Delta)\].
Thus all of the \[d-1\] critical points are contained in \[g^{-1}(\Delta)\]
if and only if \[\chi(g^{-1}(\Delta))=1\],
so that \[g^{-1}(\Delta)\] consists of a single topological disk. Similarly,
it follows by induction that all of the critical points
are contained in \[g^{-n}(\Delta)\] if and only if \[g^{-n}(\Delta)\] is
a topological disk. If every \[g^{-n}(\Delta)\] is a disk, then it follows that
the set
\[K=\bigcap g^{-n}(\Delta)\] is connected. On the other hand,
if some \[g^{-n}(\Delta)\]
consists of two or more disks, then each one of these disks must contain
a point of \[K\], which is therefore disconnected.

To complete the proof, we must allow for the possibility that \[\partial
\Delta\] may contain some post-critical point of \[g\]. Clearly there can be
at most \[d-1\] post-critical points in the annular region
\[g(\Delta)\ssm\Delta\]. Hence we can choose a disk \[\Delta_1\subset
g(\Delta)\] whose boundary avoids these post-critical points. If \[\Delta_1\]
contains \[\Delta\] in its interior, and also contains all critical values
\[g(\omega_i)\] in its interior, then it is easy to check that the pair
\[(g\,,\,g^{-1}(\Delta_1))\] is a polynomial-like map of the same degree, 
and with the same filled Julia set, but with no post-critical points in the
disk boundary. The proof then proceeds as above.\QED

{\bf Remark.} Douady and Hubbard prove much sharper statements: If \[(g\,,\,
\Delta)\] is polynomial-like of degree \[d\ge 2\], with \[K(g\,,\,
\Delta)\] connected, then there exists a polynomial map \[\psi\] of degree \[d\]
so that \[\psi\] on some neighborhood of \[K(\psi)\]
is quasi-conformally conjugate to \[g\] on a neighborhood of
\[K(g\,,\,\Delta)\]. Furthermore, this quasi-conformal conjugacy can be
chosen so as to satisfy the Cauchy-Riemann equations (in an appropriate sense)
on the compact set \[K(\psi)\].
The polynomial map \[\psi\] is then uniquely determined up to affine conjugacy.
In the case \[d=2\], one has the further statement that \[\psi\] depends
continuously on \[(g\,,\,\Delta)\].\medskip

Now let us return to the situation of Theorem 5.
\vfil\eject

{\QP{\bf Lemma 7.} \it
If the critical tableau is periodic of period \[p\ge 1\], then for any critical
puzzle piece \[P_r(c_0)\] with \[r\] sufficiently large,
the pair \[(f^{\circ p}\,,\,P_r(c_0))\]
is polynomial-like of degree two. Furthermore, the critical orbit
$$c_0~\mapsto~ c_p~\mapsto~ c_{2p}~\mapsto~\cdots$$ under \[f^{\circ p}\] is
completely contained in \[P_r(c_0)\], so that the filled Julia set
\[K(f^{\circ p}\,,\,P_r(c_0))\] is connected. In fact \[K(f^{\circ p}
\,,\,P_r(c_0))\] is equal to the intersection of the critical puzzle pieces
\[\bigcap_k P_k(c_0)\], and hence is precisely equal to the connected
component of \[K(f)\] which contains \[c_0\].\smallskip}

This is proved by
a straightforward tableau argument. Details will be left to the
reader.\QED

In this way, Branner and Hubbard show that each non-trivial component
of \[K(f)\] is homeomorphic to an appropriate quadratic Julia set.
As examples, Figures 16 and 17 illustrate the case \[p=k=1\]. 
%\vfil\eject
\end
%% 
%% \bigskip
%% 
%% 
%% \bigskip
%% 
%% \centerline{--------------------------------------------}\smallskip
%% %\eject
%% \bigskip\bigskip
%% 
%% \centerline{\bf ERRATA for ``The Yoccoz theorem on local connectivity of Julia
%% sets''}\smallskip
%% 
%% \noindent p. 11, Proof of Lemma 2: It may well be necessary to choose \[\Delta_0\]
%% as a puzzle piece of depth \[k>p\], in order to guarantee that
%% \[(f^{\circ p}\,,\,\Delta_0)\] is polynomial-like of degree two.\smallskip
%% 
%% \noindent p. 14, Problem 5: The statement should be that no {\it non-degenerate\/}
%% annulus can be\break excellent, for generic choice of \[f\]. (Otherwise, this
%% assertion would contradict\break Problem 4.)
%% \vfill
%% \rightline{Stony Brook, April 1992}
%% \end
%% 
%% 

% Written in plaintex.

%\magnification=1200

\def\QP{\medskip\leftskip=.4in\rightskip=.4in\noindent}
\def\ref{\hangindent=1pc \hangafter=1 \noindent}
\def\QED{  \rlap{$\sqcup$}$\sqcap$ \smallskip}
\def\mod{{\rm mod\;}}
\def\[{$\,}
\def\]{\,$}
\def\C{{\bf C}}
\def\R{{\bf R}}
\def\Z{{\bf Z}}

\def\car{_\heartsuit}
\def\={\;=\; }
\def\lln{\smallskip\centerline{------------------------------------------------------}\smallskip}
\font\bit=cmssi12 at 12truept
\font\tenmsy=msym10

\textfont8=\tenmsy
\mathchardef\ssm="7872
\mathsurround = 1pt
\abovedisplayskip=6pt
\belowdisplayskip=6pt
\parskip=2pt
\ifnum\pageno<2 \pageno=30 \fi

\centerline{\bf \S3. An infinitely renormalizable non locally connected Julia set.}
\medskip
%\centerline{(Class notes, J. Milnor, 3-31-92)}\bigskip

This section describes an unpublished example of Douady and Hubbard. It
begins with an outline, with few proofs, of results from
Douady and Hubbard [DH1], [DH2], [DH3], [D3].\medskip

\centerline{\bf Background Facts: Julia sets and the Mandelbrot set.}\smallskip

Let
\[	f_c(z)~=~z^2+c~. \]
%If \[f_c\] has a fixed point \[z_0\] of multiplier \[\lambda\], note that
%\[c=z_0-z_0^2\] and \[\lambda=2\,z_0\] hence
%$$	 c~=~c(\lambda)~=~\lambda/2-(\lambda/2)^2~. \eqno (1) $$
%As \[\lambda =e^{i\theta}\] traverses the unit circle, the corresponding
%point \[c(\lambda)\] traverses the {\bit cardioid\/}.\smallskip
By definition, the {\bit Mandelbrot set\/}
\[M\] is the compact set consisting
of all parameter values \[c\in\C\] such that
$${\rm the~Julia~set}~~	J(f_c)~~{\rm is~connected}\qquad\Longleftrightarrow\qquad 0~~{\rm
	has~bounded~orbit~.} $$
However, it is often convenient to identify \[M\] with the corresponding
set of polynomials \[f_c\]. A map \[f_c\in M\] is 
hyperbolic (on its Julia set) if and only if it has a necessarily
unique attracting periodic orbit. The hyperbolic maps in
\[M\] form an open subset of the plane, and each connected
component \[H\] in this open set is called a {\bit hyperbolic component\/}
in \[M\]. Let \[p\] be the period of the attracting orbit,
and let \[\lambda_p=\lambda_p(f_c)\in D\] be its multiplier. The basic facts
about hyperbolic components in \[M\] are as follows:

{\QP (1) {\it Any two maps in the same hyperbolic component
\[H\subset M\] have attracting orbits of the same period \[p\].}
We call \[p=p_H\] the {\bit period\/} of \[H\].\par}

{\QP (2) {\it Each hyperbolic
component \[H\] 
is conformally isomorphic to the open unit disk \[D\] under the
correspondence \[f_c\mapsto\lambda_p(f_c)\]. In fact this correspondence
extends uniquely to a homeomorphism between the closure \[\bar H\] and
the closed unit disk \[\bar D\].}\smallskip}

\noindent In particular, each \[H\]
has a unique {\bit center point\/} \[c_H\] which maps to \[\lambda_p(c_H)=0\],
and each boundary \[\partial H\] contains
a unique {\bit root point\/} \[r_H\in\partial H\]
which maps to \[\lambda_p(r_H)=1\].
If the map \[f_c\] has a superattracting periodic orbit, then evidently
\[c\] is the center point for one and only one hyperbolic component.

{\QP (3) {\it Similarly, if \[f_r\] has a parabolic periodic
orbit, then \[r\] is the root point for one and only one hyperbolic
component \[H\]. If the period of this orbit is \[p\] and the
multiplier is \[\lambda_p=\exp(2\pi im/p')\]
then the period of \[H\] is \[p\,p'\].}
\medskip}

By definition, the {\bit principal hyperbolic component\/} \[H\car\] is
the set of \[f_c\] having an attracting fixed point. If the multiplier is
\[\lambda_1\in D\] then a brief computation shows that
\[c=\lambda_1(2-\lambda_1)/4\], hence
$$	\lambda_1(c)~=~ 1-\sqrt{1-4c}~, $$
taking that branch of the square root which lies in the right half plane.
The boundary \[\partial H\car\] is the {\bit cardioid\/}, consisting of
all points which have the form
\[c=e^{i\theta}(2-e^{i\theta})/4\],\break so that \[\lambda_1(c)=
e^{i\theta}\].\smallskip
\eject

\midinsert
\insertRaster 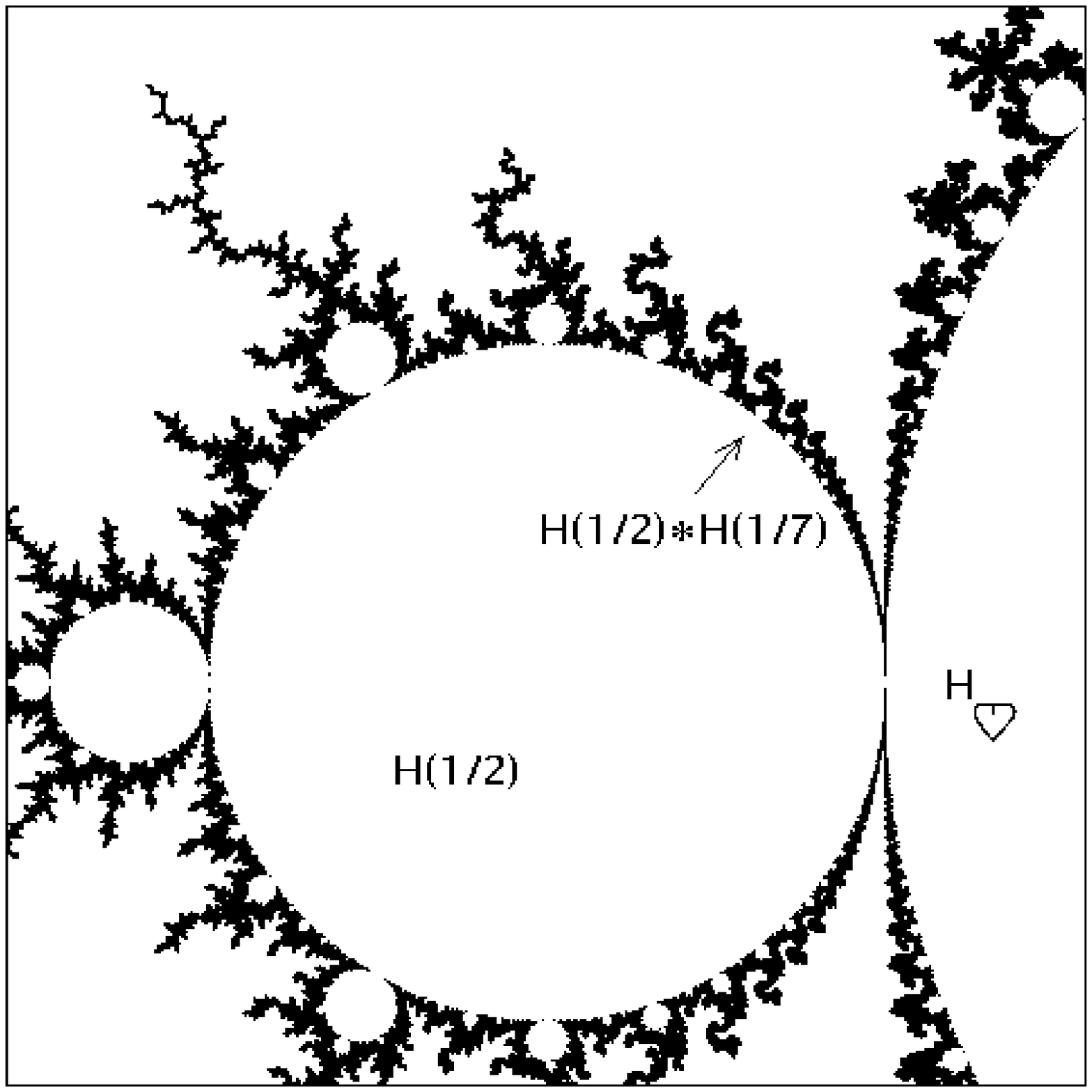 pixels 480 by 480 scaled 450
\smallskip
%{\QP\bit Figure 1. The period two component \[H(1/2)\] in the Mandelbrot set,
%with an arrow pointing to its satellite \[H(1/2)*H(1/7)\].\par}
\centerline{\bit Figure 20. The period two component \[H(1/2)\] in the
Mandelbrot}
\centerline{\bit set,
with an arrow pointing to its satellite \[H(1/2)*H(1/7)\].}

\endinsert
\vskip -.1in

{\bf Satellites.} Given any
hyperbolic component \[H\] of period \[p\ge 1\] and any root of
unity \[e^{2\pi im/p'}\ne 1\], it follows from  (2) that
there is a unique point \[r\in\partial H\] with\break \[\lambda_p(r)
=e^{2\pi im/p'}\]. According to (3), this \[r\] is the root point
of a new hyperbolic component \[H'\] of period \[p\,p'\]. We say
that \[H'\] is a {\bit satellite\/}, which is
{\bit attached\/} to \[H\] at {\bit internal angle\/} \[m/p'\]. In the special
case where \[H\] is the principal hyperbolic component \[H\car\], we will
use the notation \[H'=H(m/p')\] for this satellite, and in general we will
use the notation \[H'=H*H(m/p')\]. (Compare the discussion of ``tuning''
below.) As an example, taking \[m/p'=1/2\],
the point \[r=-3/4\] is the root point of the period two
component \[H(1/2)\], consisting of all \[c\] with \[|c+1|<1/4\].

\midinsert
\insertRaster 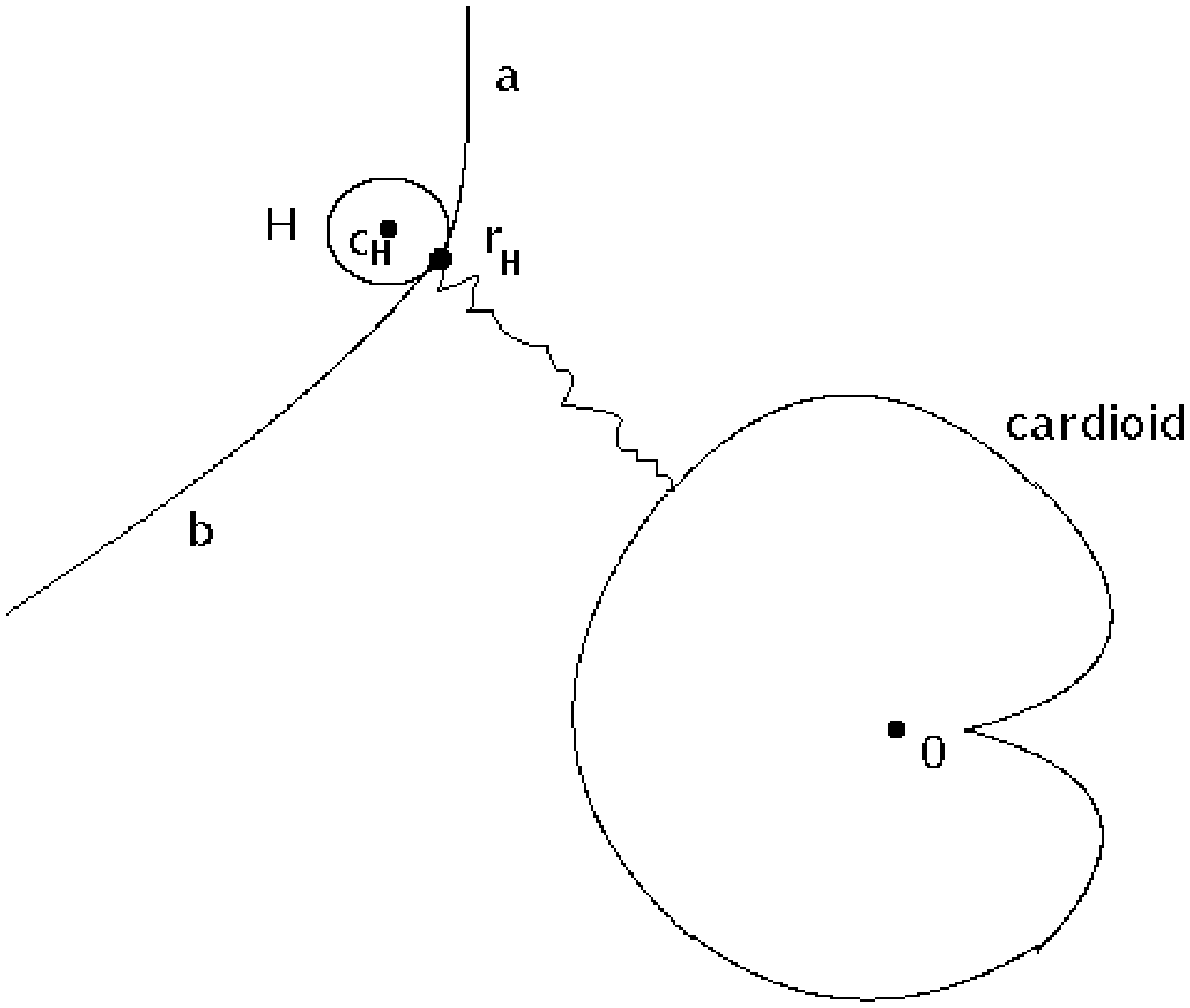  pixels 640 by 480 scaled 330
\vskip -.2in
\centerline{\bit Figure 21.
Schematic picture of external rays for}
\centerline{\bit a hyperbolic component
\[H\subset M\] in the \[c$-plane,}\smallskip

\insertRaster 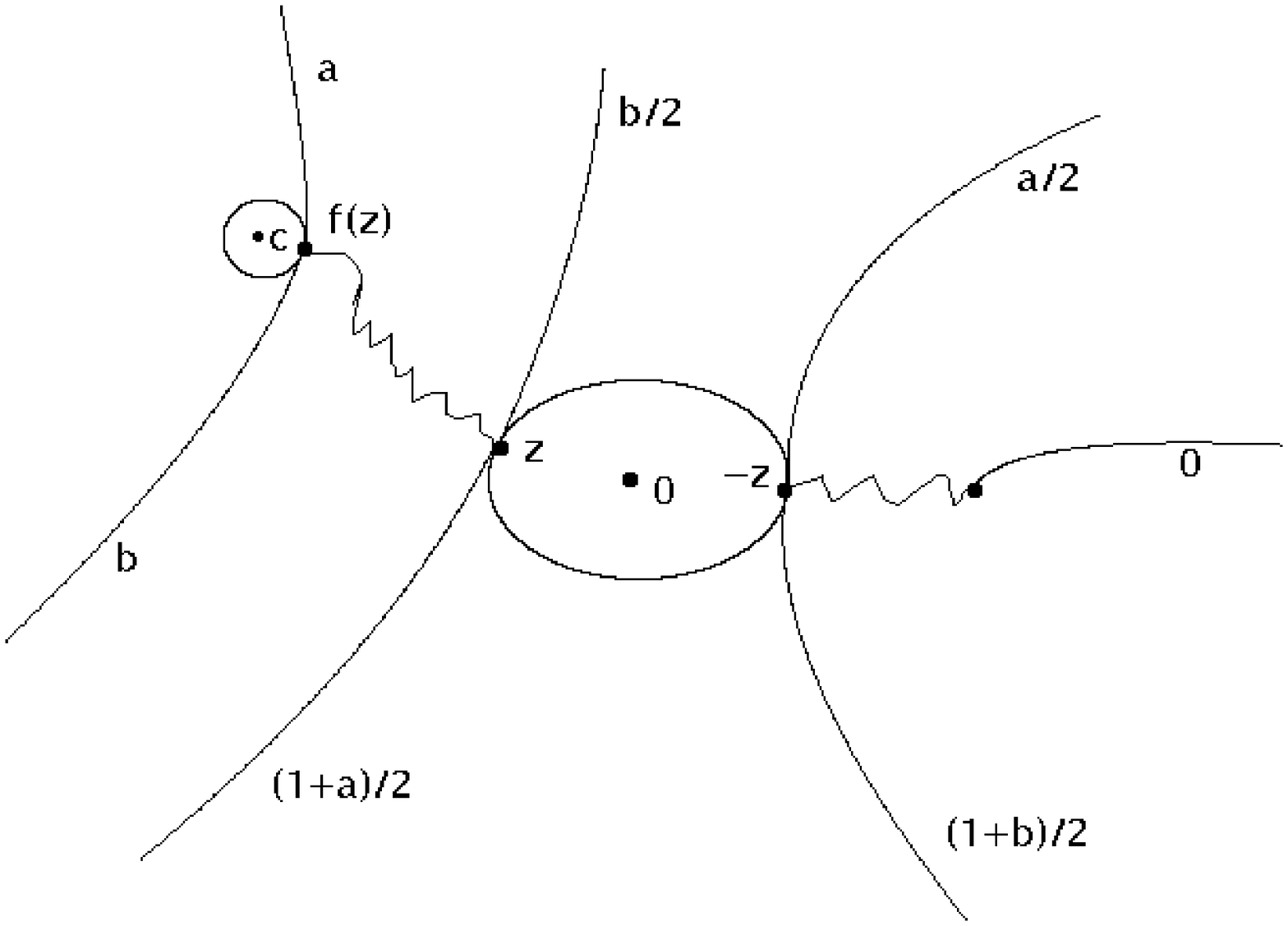 pixels 640 by 480 scaled 330
\centerline{\bit and a corresponding picture in the dynamic plane.}
%\vskip -.1in
\endinsert

{\bf External rays.} We consider not only external rays for the Julia set
in the \[z$-plane\break ($=$ {\bit dynamic plane\/}), but also external rays
for the Mandelbrot set in the \[c$-plane\break ($=$ {\bit parameter plane\/}).
Let \[H\subset M\] be a hyperbolic component of
period \[p>1\]. Then exactly two external
rays in \[\C\ssm M\] land at the
root point \[r_H\]. Let \[0<a<b<1\] be their angles. These angles have
period exactly \[p\] under doubling, and hence can be expressed as fractions
of the form \[n/(2^p-1)\]. As an example, for the satellite component
\[H(1/p)\] we have \[a=1/(2^p-1)\] and \[b=2/(2^p-1)\].

\midinsert
\insertRaster 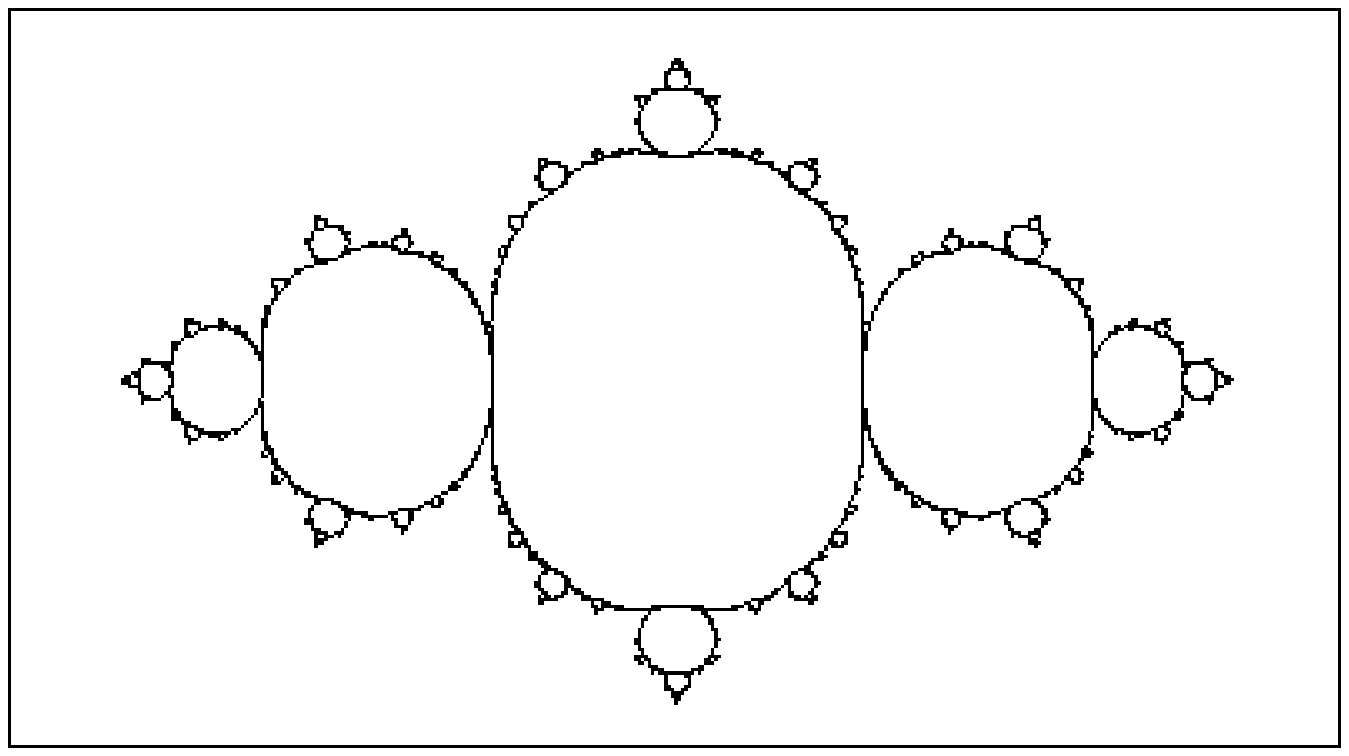 pixels 400 by 222 scaled 450
\smallskip
\centerline{\bit Figure 22. Julia set for the root point of \[H(1/2)\],}\medskip
\insertRaster 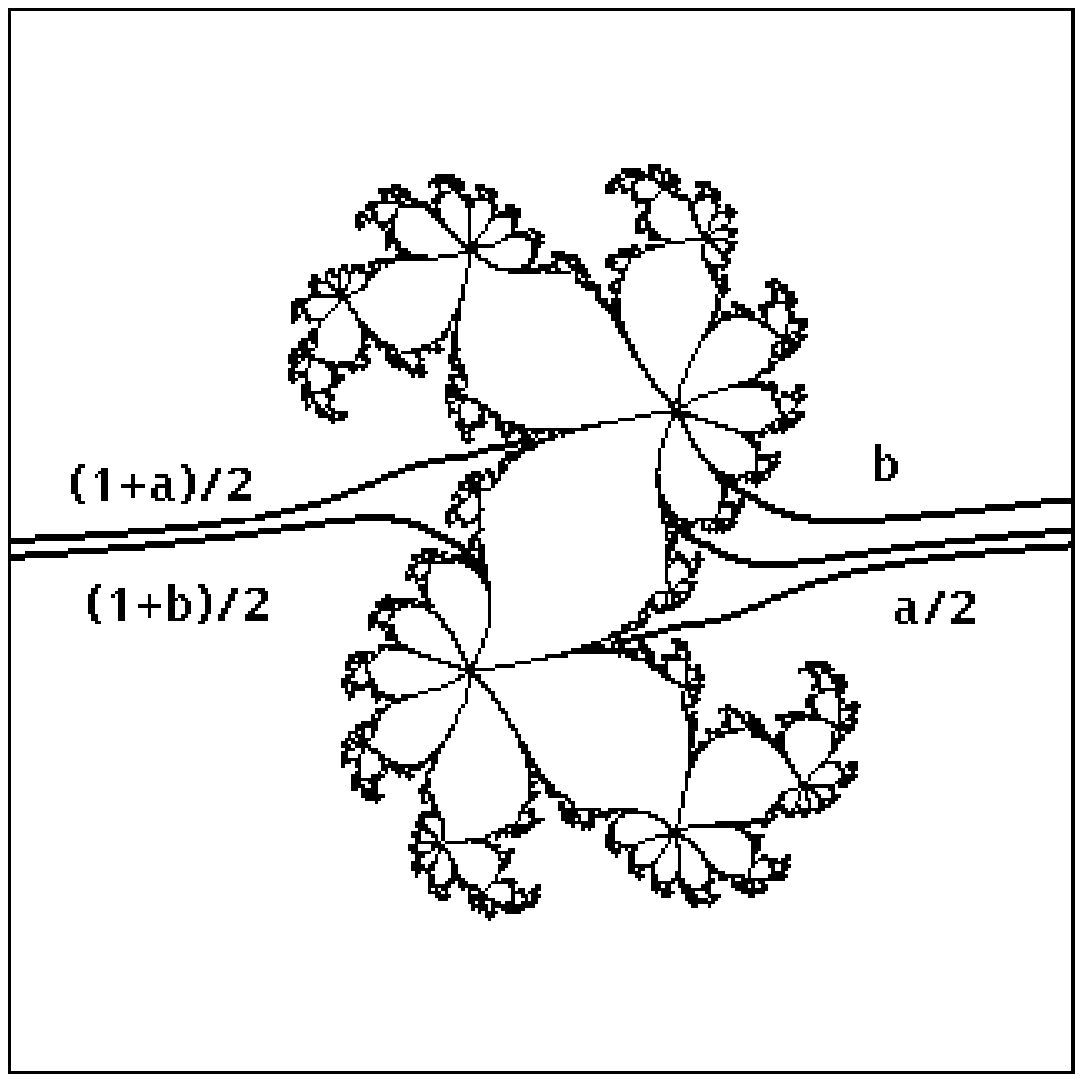 pixels 320 by 320 scaled 450
\smallskip
\centerline{\bit for the root point of \[H(1/7)\] (here the ray \[a=b/2\]
has not been labeled),}\medskip
\insertRaster 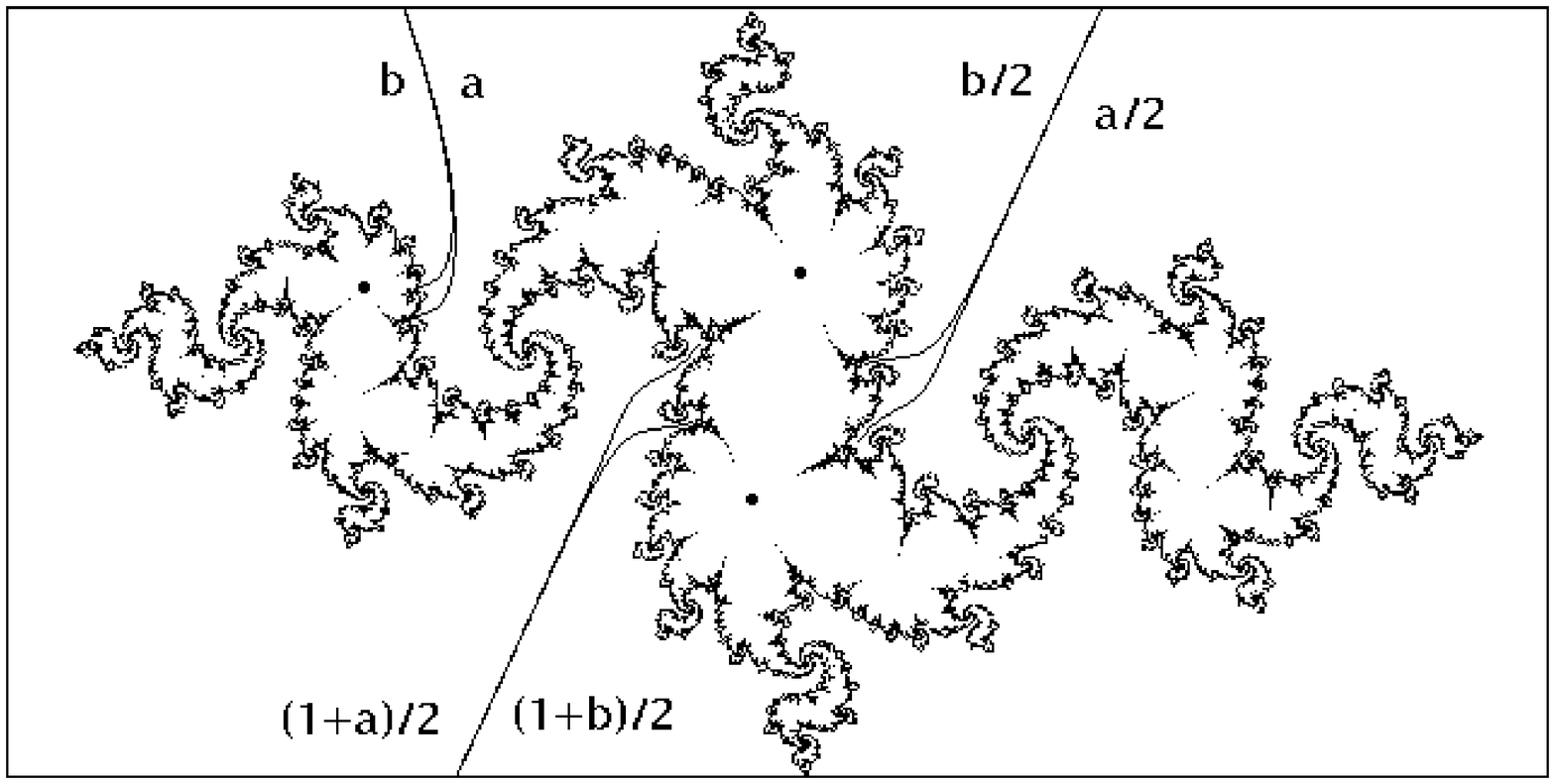 pixels 640 by 320 scaled 450
\centerline{\bit and for the root point of \[H(1/2)*H(1/7)\].}
\vskip -.05in
\endinsert

For any \[c\in H\cup
\{r_H\}\] the corresponding external rays \[R_a(c)\] and \[R_b(c)\]
in \[\C\ssm J(f_c)\]
land at a periodic point which can be described as the ``root point''
of the Fatou component containing \[c\]. See Figure 21 for a schematic
picture of the way these external rays are arranged in the parameter plane
and in the dynamic plane. In particular
$$	0\,<\,a/2\,<\,b/2\,\,~\le~\,\,a\,<\,b\,\,~\le\,\,~1+a/ 2\,<\,(1+b)/2\,<\,1~, $$
where the \[a/2$-ray and the \[b/2$-ray always land at distinct points:
exactly one of these two rays is periodic.\eject

{\bf Renormalization and Tuning.} A map \[f_c\in M\]
is {\bit renormalizable\/} if there is an integer
\[p\ge 2\] called the renormalization period,
and a closed topological disk \[\Delta\] in the \[z$-plane,
centrally symmetric about the origin, so that:
\vskip -.05in
{\QP (1) \[f^{\circ p-1}\] maps the image
\[f(\Delta)\] by a conformal isomorphism onto a disk \[f^{\circ p}(\Delta)\]
which contains \[\Delta\] in its interior, and\smallskip

\noindent (2) the entire orbit
of \[0\] under \[f^{\circ p}\] is contained in \[\Delta\].\smallskip}

\noindent
The set of all \[f_c\in M\] which are renormalizable of period \[p\] consists
of a finite number of small copies of \[M\] (sometimes
with the root point deleted).
Each of these small copies contains a unique hyperbolic component
of period \[p\]. Conversely, each hyperbolic component \[H\] of period \[p\ge
2\] determines a small copy of \[M\] which can be described as the image
of an associated mapping \[c\mapsto H*c\], which
embeds \[M\] homeomorphically onto a proper subset \[H*M\subset M\].
The elements of this form (with the possible exception of
the root point \[r_H\]), are
precisely the elements of \[M\] which are renormalizable of period \[p\].

This embedding \[c\mapsto H*c\] maps each hyperbolic component \[H'\]
of period \[p'\ge 1\] conformally onto a hyperbolic component \[H*H'\] of
period
\[pp'\] in such a way that the multiplier \[\lambda_{pp'}(H*c)\] is equal
to \[\lambda_{p'}(c)\].
In particular, this embedding maps center points to center points and root
points to root points. It carries
the principal hyperbolic component \[H\car\]
onto \[H\] itself, and maps the satellite \[H(m/p')\] of \[H\car\] onto the
satellite \[H*H(m/p')\] of \[H\].
This $*$-product operation between hyperbolic components
is associative, and has the
principal hyperbolic component \[H_\heartsuit\] as two-sided identity
element. Thus it makes the collection of hyperbolic
components into a monoid (=associative semigroup with identity)
which operates as a semigroup of embeddings of \[M\] into itself. This
monoid is free non-commutative.\smallskip

Intuitively, the Julia set for \[H*c'\] is obtained from the Julia set of
a map \[f_c\] belonging to the hyperbolic component
\[H\] by replacing each bounded component of \[\C\ssm J(f_c)\] by a copy
of \[J(f_{c'})\] and defining the dynamics appropriately so that only
the copy which is pasted in place of the critical component is mapped
non-homeomorphically. The result is described as \[f_c\] {\bit tuned by\/}
\[f_{c'}\].\smallskip

In terms of external angles, let \[a<b\] be the two external angles for
the root point of \[H\]. These angles have periodic binary
expansions of the form
 \[.a_1\cdots a_p\,a_1\cdots a_p\cdots\] and \[.b_1\cdots b_p\,
b_1\cdots b_p\cdots\], both with
period exactly equal to
\[p\]. {\bf Assertion:} If the point \[c'\in\partial M\] is the landing point
of a ray with angle \[t\in\R/\Z\], then the image
\[H*c'\] is the landing point of a ray whose binary expansion can be obtained
by inserting the \[p$-tuple \[a_1\cdots a_p\] in place of each zero in the
binary expansion of \[t\] and the \[p$-tuple \[b_1\cdots b_p\] in place of
each one. (See [D3].)\bigskip

\centerline{\bf A non locally connected Julia set.}\medskip

Consider a sequence of integers \[1<p_1<p_2<\cdots\].

{\QP{\bf Theorem 6.} \it If this sequence diverges to infinity sufficiently
rapidly, then the sequence of subsets
\[H(1/p_1)*H(1/p_2)*\cdots*H(1/p_k)*M\] intersects in a single point \[\omega
\in M\]
with the property that the Julia set \[J(f_\omega)\] is not locally connected.
\medskip}

The proof begins as follows. It will be convenient to use the abbreviation
\[r(k)\] for the root point \[r_{H(1/p_k)}\] on the cardioid.
Start with any \[p_1>1\] and consider the
embedding \[c\mapsto H(1/p_1)*c\] from \[M\] into itself, which carries
the period 1 root point \[r\car =1/4\] to the root point \[r(1)=
r_{H(1/p_1)}\]. Since this
embedding is continuous, we can place \[H(1/p_1)*c\] as close as we like
to \[r(1)\] by choosing \[c\] close to \[1/4\]. In particular, if \[p_2\]
is sufficiently large, then the entire \[1/p_2$-limb will be close to
\[1/4\] by the Yoccoz inequality, hence every point \[c\in
H(1/p_1)*H(1/p_2)*M\] will certainly
be close to \[r(1)\]. Under such a small perturbation,
note that the parabolic
fixed point of \[f_{r(1)}\] splits up
into a repelling fixed point for \[f_c\], together with a nearby period
\[p_1\] orbit. Thus, by choosing  \[p_2\] large, we can place this entire
period \[p_1\] orbit into an \[\epsilon$-neighborhood of
the parabolic fixed point of \[f_{r(1)}\].\smallskip

Similarly, choosing \[p_3\] even larger, for any
\[c\in H(p_1)*H(p_2)*H(p_3)*M\] we can guarantee
that the parabolic period \[p_1\] orbit for \[H(1/p_1)*r(2)\]
is replaced by a nearby
period \[p_1\,p_2\] orbit. In particular, we can guarantee that this
new orbit lies within the \[\epsilon+\epsilon/2\] neighborhood of
the original parabolic fixed point. Continuing inductively, we can guarantee
that all of the new orbits which are constructed lie within the \[2\epsilon\]
neighborhood of the original fixed point.
Furthermore, we can easily
guarantee that these successive copies of \[M\] have diameter shrinking to zero.

Thus, if \[\omega\] is the unique point in the intersection,
we know that the Julia set \[J(f_\omega)\]
contains orbits of period  \[1\,,\,p_1\,,
\,p_1\,p_2\,,\,\ldots\], all lying within the \[2\epsilon\] neighborhood of
the original parabolic fixed point. These various periodic points can
be described as the landing points of external rays of angles \[a_k<b_k\]
where \[0<a_1<a_2<\cdots<b_2<b_1\]. Furthermore, it is easy to check
that the difference \[b_k-a_k\] tends rapidly to zero, so that
\[~\lim_{k\to\infty} a_k=\lim_{k\to\infty} b_k\,\]. On the other hand, as noted
in the preceding section, the external rays of angle \[a_k/2\] and \[b_k/2\]
land at different points, which are negatives of each other. Each such landing
point is either close to the original parabolic fixed point, or close to
its negative. Thus we have a sequence of angles \[a_k/2\] and \[b_k/2\]
tending to a common limit, yet the landing points of the corresponding
external rays for \[J(f_\omega)\] do not tend to a common limit. According
to Caratheodory, this implies that \[J(f_\omega)\] is not locally connected.
\QED

More generally,
if \[H_1\,,\, H_2\,,\,\ldots\] is any sequence of hyperbolic components
which are not centered along the real axis, then it is again true that
\[~\lim_{k\to\infty} a_k=\lim_{k\to\infty} b_k\,\]. Let \[\omega\in M\]
be any element of the
%which converges sufficiently rapidly to \[1/4\], then a similar argument
intersection \[\bigcap_k H_1*H_2*\cdots*H_k*M\]. Whenever the landing
points in \[J(f_\omega)\] of the rays of angle \[a_k/2\] fail to converge
to the critical point, it follows that \[J(f_\omega)\] is not locally
connected. In particular, this is true whenever the sequence of components
\[H_k\] converges sufficiently rapidly to the root point \[r\car=1/4\].
\smallskip

This seems to be the only known obstruction to local connectivity in the
infinitely renormalizable case.
Thus whenever the landing points of the \[a_k/2\] rays do converge to the
critical point, we can ask whether \[J(f_\omega)\] is in fact locally
connected.
\bigskip\bigskip
\vfil\eject

%\input ../header.tex

%\footline={\hss\tenrm A-\folio\hss}
%\pageno=1

\def\md{{\rm mod}}
\def\ar{{\rm area}}
\def\L{{L}}
\def\T{{\bf T}}

\centerline{\bf Appendix. Length-Area-Modulus Inequalities.}\bigskip

The most basic length-area inequality is the following. Let \[I^2\subset\C\]
be the open unit square consisting of all \[z=x+iy\] with \[0<x<1\] and
\[0<y<1\]. By a {\bit conformal metric~} on \[I^2\] we mean a metric
of the form
$$	ds\= \rho(z)|dz| $$
where \[z\mapsto \rho(z)>0\] is any strictly positive continuous real
valued function on the open square. In terms of such a metric, the
{\bit length~} of a smooth curve \[\gamma:(a,b)\to I^2\] is defined to be the
integral
$$	{L}_\rho(\gamma)\= \int_a^b \rho(\gamma(t))|d\gamma(t)|\;, $$
and the {\bit area~} of a region \[U\subset I^2\] is defined to be
$$	\ar_\rho(U)\= \int\!\int_U\, \rho(x+iy)^2dx\,dy\;. $$
In the special case of the Euclidean metric \[ds=|dz|\], with \[\rho(z)\]
identically equal to 1, the subscript \[\rho\] will be omitted.

{\QP{\bf Theorem A.1.} \it If \[\ar_\rho(I^2)\] (the integral over the entire
square) is finite, then for
Lebesgue almost every \[y\in(0,1)\] the length \[\L_\rho(\gamma_y)\] of the
horizontal line \[\gamma_y:t\mapsto
(t,y)\] at height \[y\] is finite. Furthermore, there exists
\[y\] so that
$$	\L_\rho(\gamma_y)^2\;\le\;\ar_\rho(I^2) \;.\eqno (1) $$
In fact, the set consisting of all \[y\in(0,1)\] for which this inequality
is satisfied has positive Lebesgue measure.\smallskip}

{\bf Remark 1.} Evidently this inequality is best posible. For in the
case of the Euclidean metric \[ds=|dz|\] we have
$$	{L}(\gamma_y)^2\;=\;\ar(I^2)\;=\;1\;. $$

{\bf Remark 2.} It is essential here that we use a square, rather than
a rectangle. If we consider instead a rectangle \[R\] with base
\[\Delta x\] and height \[\Delta y\], then the corresponding inequality
would be
$$	\L_\rho(\gamma_y)^2\;\le\;{\Delta x\over\Delta y}\ar_\rho(R)\eqno (2) $$
for a set of \[y\] with positive measure.\smallskip

{\bf Proof of A.1.} 
We use the Schwarz inequality
$$	\left(\int_a^b f(x)g(x)\,dx\right)^2\;\le\;\left(\int_a^b f(x)^2\,
dx\right)\cdot\left(\int_a^b g(x)^2\,dx\right)\;, $$
which says (after taking a square root) that the inner product of any
two vectors in the Euclidean vector space
of square integrable real functions on an interval is less than or
equal to the product of their norms.
We may as well consider the more general case of a
rectangle\break \[R=(0,\,\Delta x)\times(0,\,\Delta y)\].
Taking \[f(x)\equiv 1\] and \[g(x)=
\rho(x,y)\] for some fixed \[y\], we obtain
$$	\bigg(\int_0^{\Delta x} \rho(x,y)\,dx\bigg)^2\;\le\; \Delta x
 \int_0^{\Delta x} \rho(x,y)^2\,dx\;, $$
or in other words
%for each \[y\in(0\,,\,\Delta y)\]. In other words
$$	L_\rho(\gamma_y)^2\;\le \; \Delta x\int_0^{\Delta x} \rho(x,y)^2dx\;, $$
%where \[\gamma_y\] is the horizontal curve \[x\mapsto(x,y)\].
for each constant height \[y\].
Integrating this inequality over the
interval \[0<y<{\Delta y}\] and then dividing by \[\Delta y\], we get
$$	{1\over\Delta y}\int_0^{\Delta y} L_\rho(\gamma_y)^2dy\;\le\;
	{\Delta x\over\Delta y}\; {\rm area}_\rho(A)\;.\eqno (3) $$
In other words, the {\bit average~} over all \[y\] in the interval
\[(0,\,\Delta y)\] of \[L_\rho(\gamma_y)^2\] is less than or equal
to \[{\Delta x\over\Delta y}\; {\rm area}_\rho(A)\,\]. Further details of the
proof are straightforward.\QED\smallskip

Now let us form a cylinder \[{\cal C}\] of circumference \[\Delta x\] and height
\[\Delta y\] by gluing the left and right edges of our
rectangle together. More precisely, let \[{\cal C}\] by the quotient space which
is obtained from the infinitely wide strip \[0<y<\Delta y\] in the \[z$-plane
by identifying each point \[z=x+iy\] with its translate \[z+\Delta x\].
Define the {\bit modulus\/} \[\md({\cal C})\] of such a cylinder to be the ratio
\[\Delta y/\Delta x\] of height to circumference.
By the {\bit winding number\/} of a closed curve \[\gamma\] in \[{\cal C}\] we
mean the integer
$$	w\={1\over\Delta x}\;\oint_\gamma dx\;. $$

{\QP{\bf Theorem A.2 (Length-Area Inequality for Cylinders).}
\it For any conformal metric
\[\rho(z)|dz|\] on the cylinder \[{\cal C}\] there exists some simple closed
curve \[\gamma\] with winding number \[+1\] whose length
\[\;\L_\rho(\gamma)=\oint_\gamma \rho(z)|dz|\;\] satisfies the inequality
$$	\L_\rho(\gamma)^2\;\le \; \ar_\rho(A)/\md(A)\;. \eqno (4) $$
Furthermore, this result is best possible: If we use the Euclidean metric
\[|dz|\] then
$$      \L(\gamma)^2\;\ge\;\ar(A)/\md(A) \eqno (5)
$$
for {\bit every} such curve \[\gamma\].\smallskip}

\noindent {\bf Proof.} Just as in the proof of A.1, we find a horizontal curve
\[\gamma_y\] with
$$ \L_\rho(\gamma_y)^2\;\le\;{\Delta x\over\Delta y}\ar_\rho({\cal C})\;=\;
{\ar_\rho({\cal C})\over\md({\cal C})}\;. $$
On the other hand in the Euclidean case, for any
closed curve \[\gamma\] of winding number one we have
$$ \L(\gamma)\;=\;\oint_\gamma |dz|\;\ge\;\oint_\gamma dx\;=\;\Delta x\;, $$
hence \[\L(\gamma)^2\;\ge\; (\Delta x)^2\;=\;\ar({\cal C})/\md({\cal C})\].\QED

{\bf Definitions.}
A Riemann surface \[A\] is said to be an {\bit annulus\/} if it is conformally
isomorphic to some cylinder. An embedded annulus \[A\subset {\cal C}\] is said
to be {\bit essentially embedded~} if it contains a curve which has winding
number one around \[{\cal C}\].

\noindent Here is an important consequence
of Theorem A.2.

{\QP{\bf Corollary A.3 (An Area-Modulus Inequality).} \it Let \[A\subset {\cal C}\]
be an essentially embedded annulus in the cylinder \[{\cal C}\], and suppose that
\[A\] is conformally isomorphic to a cylinder \[{\cal C}'\]. Then
$$	\md({\cal C}')\;\le\;{\ar(A)\over\ar({\cal C})}\,\md({\cal C})\;. \eqno (6)$$
In particular:
$$	\md({\cal C}')\;\le\;\md({\cal C})\;. \eqno (7) $$
\smallskip}

{\bf Proof.} Let \[\zeta\mapsto z\] be the embedding of \[{\cal C}'\] onto
\[A\subset {\cal C}\]. The Euclidean metric \[|dz|\] on \[{\cal C}\], restricted to \[A\],
pulls back to some conformal metric \[\rho(\zeta)|d\zeta|\]
on \[{\cal C}'\], where \[\rho(\zeta)=|dz/d\zeta|\]. According to A.2, there
exists a curve \[\gamma'\] with winding number \[1\] about \[{\cal C}'\] whose length
satisfies
$$	L_\rho(\gamma')^2\le\ar_\rho({\cal C}')/\md({\cal C}') \;. $$
%=\ar(A)/\md({\cal C}')\;. $$
This length coincides with the Euclidean length \[\L(\gamma)\] of the
corresponding curve \[\gamma\] in \[A\subset {\cal C}\], and \[\ar_\rho({\cal C}')\]
is equal to the Euclidean area
\[\ar(A)\], so we can write this inequality as
$$	L(\gamma)^2\le\ar(A)/\md({\cal C}') \;. $$
But according to (5) we have
$$	\ar({\cal C})/\md({\cal C})\;\le\;L(\gamma)^2\;. $$
Combining these two inequalities, we obtain
$$ \ar({\cal C})/\md({\cal C})\;\le\;\ar(A)/\md({\cal C}')\;, $$
which is equivalent to the required inequality (6).\QED

{\QP{\bf Corollary A.4.} \it The modulus of a cylinder is a well defined
conformal invariant.\smallskip}

%\noindent Hence the {\bit modulus~} of an annulus can be defined as the modulus
%of any conformally isomorphic cylinder.\smallskip

{\bf Proof.} If \[{\cal C}'\] is conformally isomorphic to \[{\cal C}\]
then (7) asserts
that\break \[\md({\cal C}')\le\md({\cal C})\], and similarly \[\md({\cal C})\le\md({\cal C}')\].\QED

It follows that the {\bit modulus~} of an annulus \[A\]
can be defined as the modulus of any conformally isomorphic cylinder.
Furthermore, if \[A\] is essentially embedded in some other annulus
\[A'\], then \[\md(A)\le\md(A')\].

{\QP{\bf Corollary A.5 (Gr\"otzsch Inequality).} \it Suppose that \[A'
\subset A\] and\break \[A''\subset A\] 
are two disjoint annuli, each essentially embedded in \[A\]. Then 
$$ \md(A')+\md(A'')\;\le\;\md(A)\;. $$
\smallskip}
\eject

{\bf Proof.} We may assume that \[A\] is a cylinder \[{\cal C}\].
According to (6) we have
$$	\md(A')\;\le \;{\ar(A')\over\ar({\cal C})}\,\md({\cal C})\;,\qquad
      \md(A'')\;\le \;{\ar(A'')\over\ar({\cal C})}\,\md({\cal C})\;. $$
where all areas are Euclidean. Using the inequality
$$	\ar(A')+\ar(A'')\le \ar({\cal C})\;, $$
the conclusion follows.\QED\smallskip

Up to this point, we have only considered cylinders or annuli of finite
modulus. If we take an arbitrary Riemann surface with free cyclic
fundamental group, then it is always conformally isomorphic to some cylinder,
providing that we allow also the possibility of a one-sided infinite or
two-sided infinite cylinder (that is, either the upper half-plane or the
full complex plane modulo the identification \[z\equiv z+1\]). By definition,
a Riemann surface conformally isomorphic to such an infinite cylinder
will be called an {\bit annulus of infinite modulus}.\smallskip

Now consider the following situation. Let \[U\subset\C\] be a bounded
simply connected open set, and let \[K\subset U\] be a compact subset
such that the difference \[A=U\ssm K\] is an annulus (which may have finite
or infinite modulus).
%. It is well
%known that such an annulus must be conformally isomorphic to a finite or
%infinite cylinder. (By definition an infinite cylinder, that is a cylinder of
%infinite height, has infinite modulus. In general
%an infinite cylinder may be either
%one-sided infinite --- conformally isomophic to a punctured disk, or two-sided
%infinite --- conformally isomorphic to the punctured plane. However,
%in the present
%case, since \[U\] is bounded, only the one-sided infinite case can occur.)

{\QP{\bf Corollary A.6.} \it Suppose that \[K\subset U\] as described
above. Then \[K\]
reduces to a single point if and only if the annulus \[A=U\ssm K\] has infinite
modulus. Furthermore, the diameter of \[K\] is bounded by the inequality
$$	4\;{\rm diam}(K)^2\;\le\;{\ar(A)\over\md(A)}\;\le\;{\ar(U)
 \over\md(A)}\;. \eqno (8)  $$
\smallskip}

{\bf Proof.} According to A.2, there exists a curve with winding number
one about \[A\] whose
length satisfies \[L^2\le \ar(A)/\md(A)\]. Since \[K\] is enclosed
within this curve, it follows easily that \[{\rm diam}(K)\le\L/2\],
and the inequality (8) follows. Conversely, if \[K\] is a single point then
using (7) we see easily that \[\md(A)=\infty\].\QED

{\QP{\bf Corollary A.7.
(Branner-Hubbard).} \it Let \[K_1\supset K_2\supset K_3
\supset\cdots\] be compact subsets of \[\C\] with each \[K_{n+1}\]
contained in the interior of \[K_n\]. Suppose further that each interior
\[K_n^o\] is simply connected, and that each difference \[A_n=
K_n^o\ssm K_{n+1}\] is an annulus. If \[\sum_1^\infty\md(A_n)\] is infinite,
then the intersection
\[\bigcap K_n\] reduces to a single point.\smallskip}

{\bf Proof.} It follows inductively from the Gr\"otzsch inequality that the
modulus\break \[\md(K_1^o\ssm K_n)\] tends to infinity as \[n\to\infty\].
Hence by (7) and A.7 the intersection of the \[K_n\] is a point.\QED\smallskip

{\QP{\bf Corollary A.8 (McMullen Inequality).} \it Again suppose that
\[K\subset U\subset\C\], and that \[A=U\ssm K\] is an annulus. Then
$$      \ar(K)\;\le\;{\ar(U)/ e^{4\pi\,\md(A)}}\;\;. \eqno (9)$$
\smallskip}

{\bf Proof.} (Compare [BH, II].) We will first prove the weaker inequality
$$      \ar(K)\;\le\;{\ar(U)\over 1+4\pi\,\md(A)}\;. \eqno (10) $$
The {\bit isoperimetric inequality~} asserts that the area enclosed by a plane
curve of length \[\L\] is at most \[\L^2/(4\pi)\], with equality if
and only if the curve is a round circle. (See for example [CR].) Combining
this with the proof of A.6, we see that
$$	\ar(K)\;\le\;{\L^2\over 4\pi}\;\le{\ar(A)\over 4\pi\,\md(A)}\,. $$
Writing this inequality
as \[\;4\pi\,\md(A)\le \ar(A)/\ar(K)\;\] and adding \[+1\]
to both sides we obtain the completely equivalent inequality
\[	1+4\pi\,\md(A)\le\ar(U)/\ar(K)\;. \] This
in turn is equivalent to (10).

To sharpen this inequality, we proceed as follows.
Cut the annulus \[A\] up into \[n\]
concentric annuli \[A_i\], each of modulus equal to \[\md(A)/n\]. Let \[K_i\]
be the bounded component of the complement of \[A_i\], and assume that
these annuli are nested so that\break \[A_i\cup K_i=K_{i+1}^o\] with \[K_1=K\].
Let \[K_{n+1}^o=A\cup K=U\]. Then $$\ar(K_{i+1})/
\ar(K_i)\,\ge\,1+4\pi\,\md(A)/n $$ by (10), hence
$$	\ar(U)/\ar(K)\;\ge\;(\,1\,+\,4\pi\,\md(A)/n\,)^n\;, $$
where the right hand side converges to \[\;e^{4\pi\,\md(A)}\;\]
as \[n\to\infty\].\QED
%\eject

As a final illustration of modulus-area inequalities,
consider a {\bit flat torus~} \[{\bf T}=\C/\Lambda\]. Here \[\Lambda
 \subset \C\] is to be a 2-dimensional {\bit lattice\/}, that is
an additive subgroup of the complex numbers, spanned by
two elements \[\lambda_1\] and \[\lambda_2\] where \[\lambda_1/\lambda_2
\not\in\R\]. Let \[A\subset\T\] be an embedded annulus.

By the ``{\bit winding number\/}'' of \[A\] in \[\T\] we will mean the lattice
element \[w\in\Lambda\] which is constructed as follows. Under the
universal covering map \[\C\to\T\], the central curve of \[A\] lifts
to a curve segment which joins some point \[z_0\in\C\] to a translate
\[z_0+w\] by the required lattice element. We say that \[A\subset\T\] is an
{\bit essentially embedded annulus~} if \[w \ne 0\].

{\QP{\bf Corollary A.9 (Bers Inequality).} \it If the annulus \[A\]
is embedded in the flat torus \[\T=\C/\Lambda\] with winding number
\[w \in\Lambda\], then
$$	\md(A)\;\le\;{\ar(T)\over|w |^2}\;. \eqno (11) $$
\smallskip}

\noindent Roughly speaking, if \[A\] winds many times around the torus,
so that \[|w |\] is large, then \[A\] must be very skinny. A slightly
sharper version of this inequality is given in Problem A-3 below.
\smallskip

{\bf Proof.} Choose a cylinder \[C'\] which is conformally isomorphic
to \[A\]. The Euclidean metric \[|dz|\] on \[A\subset\T\] corresponds to some
metric \[\rho(\zeta)|d\zeta|\] on \[C'\], with
$$ \ar_\rho(C')\= \ar(A)\;.$$
By A.2 we can choose a curve \[\gamma'\] of winding number one on \[C'\],
or a corresponding curve \[\gamma\] on \[A\subset\T\], with
$$	\L(\gamma)^2\=\L_\rho(\gamma')^2\;\le\; {\ar_\rho(C')\over\md(C')}
 \={\ar(A)\over\md(A)}\;\le\;{ \ar(T)\over\md(A)}\;. $$
Now if we lift \[\gamma\] to the universal covering space \[\C\] then
it will join some point \[z_0\] to \[z_0+w \]. Hence its Euclidean
length \[\L(\gamma)\] must satisfy \[\L(\gamma)\ge
|w|\]. Thus
$$ |w|^2\;\le\; {\ar(T)\over\md(A)}\;, $$
which is equivalent to the required inequality (11).\QED\smallskip
%\eject

\medskip\lln\medskip

{\bf Concluding Problems:}\smallskip

{\bf Problem A-1.} In the situation of Theorem
A.1, show that more than half of the
horizontal curves \[\gamma_y\] have length \[L_\rho(\gamma_y)\le
\sqrt{2\,\ar_\rho(I^2)}\]. (Here ``more than half'' is to be interpreted
in the sense of Lebesgue measure.)\smallskip

{\bf Problem A-2.}  Show that the
converse to A.7 is false: A nest \[K_1\supset K_2\supset K_3
\supset\cdots\] of compact subsets of \[\C\] may intersect in a point
even if \[\sum_1^\infty\md(K_n^o\ssm K_{n+1})\] is finite.
(For example, do this by showing that a closed
disk \[\overline{D}'\] of radius \[1/2\] can be embedded in the open unit
disk \[D\] so that the complementary annulus \[A=D\ssm\overline{D}'\] has
modulus arbitrarily close to zero.)\smallskip

{\bf Problem A-3 (Sharper Bers Inequality).} If the flat torus
\[\T=\C/\Lambda\] contains several disjoint
annuli \[A_i\], all with the same ``winding number'' \[w\in\Lambda \],
show that
$$	\sum\md(A_i)\;\le\; \ar(\T)/|w |^2\;. $$
If two essentially
embedded annuli are disjoint, show that they necessarily
have the same winding number.\bigskip
\vfil\eject
\end

\vfill\eject
%\magnification=1200

\def\oldQP{\narrower\smallskip\noindent}
\def\QP{\medskip\leftskip=.4in\rightskip=.4in\noindent}
\def\[{$\,}
\def\]{\,$}
\def \QED { \rlap{$\sqcup$}$\sqcap$\smallskip}
\def\mod{{\rm mod}}
\def\C{{\bf C}}
\def\inter{{\rm interior}}
\def\ref{\smallskip\hangindent=1pc \hangafter=1 \noindent}
\def\scd{{\rm scd}}
\font\bit=cmssi12 at 12truept
\font\tenmsy=msym10
\textfont8=\tenmsy
\mathchardef\ssm="7872
\abovedisplayskip=6pt
\belowdisplayskip=6pt
\parskip=0pt
\ifnum\pageno<2 \pageno=44 \fi
\input psfig

\centerline{\bf Renormalization and Tuning}\vskip .04in
\centerline{\bf ERRATA %}} % concerning Renormalization and Tuning}}
%\vskip .02in\centerline{\bf 
for}\vskip .02in
\centerline{{\bf ``Local Connectivity of Julia Sets: Expository Lectures''}}
%~ by J. Milnor}
%\centerline{\bf Errata for ``Local Connectivity of Julia Sets:
%Expository Lectures'' by J. Milnor}\smallskip

%\centerline{J. Milnor}
\bigskip

These notes make the following claim on p. 32:

{\QP\it Any map \[f_c(z)=z^2+c\]
in the Mandelbrot set \[M\] which is {\bit renormalizable\/} must belong to
one of the small embedded copies of \[M\] which are constructed by
Douady-Hubbard {\bit tuning\/.}\medskip} 

\noindent Curt McMullen points out
that this is false; here is a counter-example.
Let $$c=c_1\break\approx .41964338+.60629073\,i~, $$
%far out in the \[1/4$-limb of \[M\],
so that the polynomial \[f=f_c\] is post-critically finite, with critical orbit
$$0~\mapsto~ c_1~\mapsto~ c_2~\mapsto~ c_3~\leftrightarrow~ c_4=-c_2~. $$
(See Figure A. The number \[c\] can be described as the landing point of the
external ray of angle \[1/12\] for \[M\]. It is nearly the right-most point on
the \[1/4$-limb of \[M\].)
This map \[f\] is renormalizable of period \[2\]. In fact the union
$$	\Delta_0~=~\widehat P_3(0)\cup\widehat P_2(c_2)\cup\widehat
 P_2(-c_2) $$

\goodbreak
\midinsert
\insertRaster 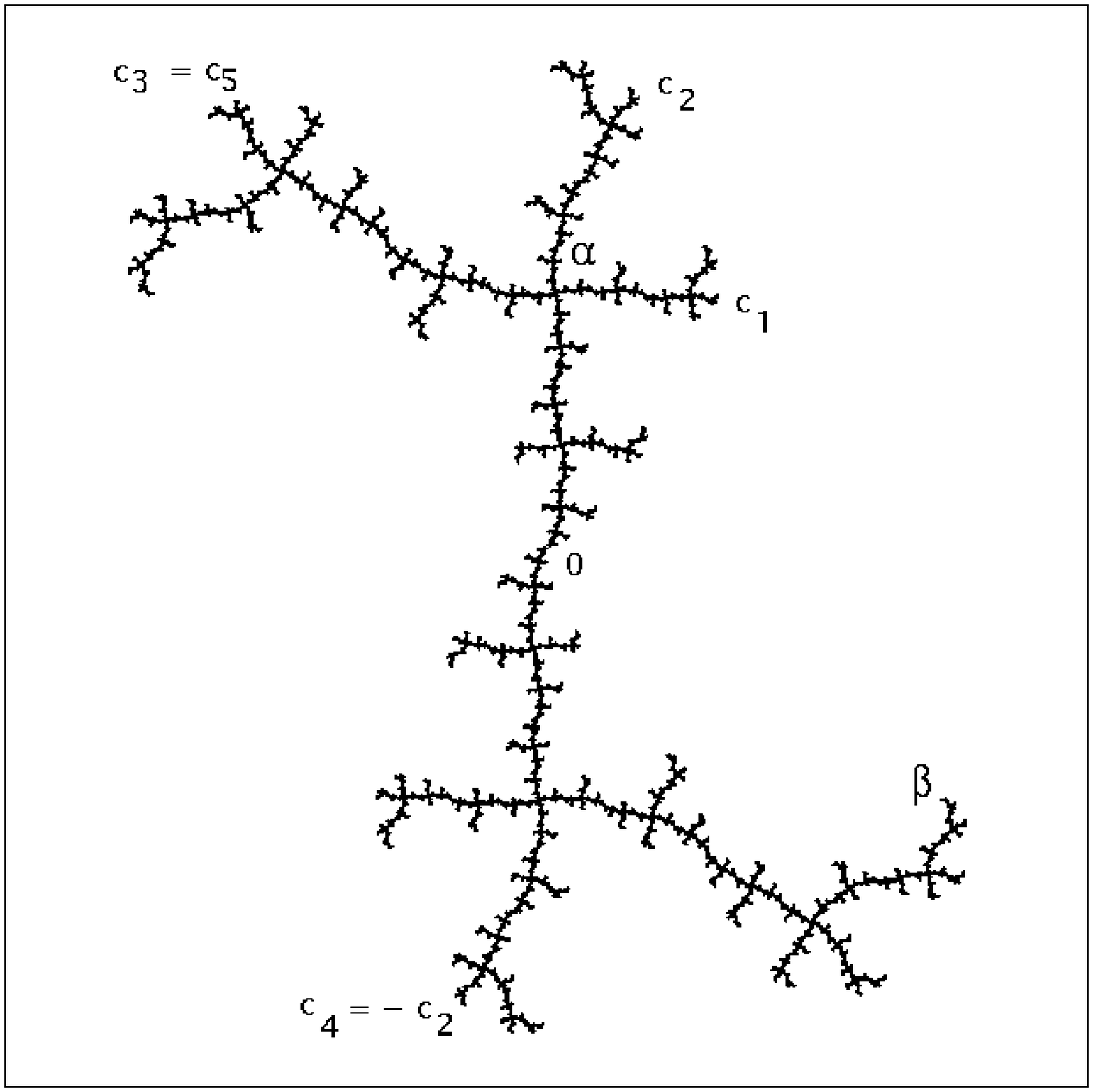 pixels 640 by 640 scaled 450
\smallskip
\centerline{\bit Figure A.
Julia set for \[~~f\,:\,z\,\mapsto\, z^2 +.41964338+.60629073\,i\].}
\endinsert

\noindent
is a topological disk, centrally symmetric about the origin, and the
composition \[f^{\circ 2}\] maps \[\Delta_0\] by a two-fold branched
covering onto a disk
$$	\Delta_2~=~\widehat P_0(0)\cup\widehat P_1(c_2) $$
which compactly contains it. Furthermore, the orbit of \[0\] under
\[f^{\circ 2}\] is completely\break contained in \[\Delta_0\].
(For definitions and notation,
see pp. 11, 12, 32.) The proof is not difficult.
In this example, the quadratic-like map \[f^{\circ 2}|\Delta_0\] is
hybrid-equivalent (in the sense of [DH3])
to the Tchebycheff map \[z\mapsto z^2-2\]. In particular, the filled Julia set
\[K_0= K(f^{\circ 2}|\Delta_0)\] is homeomorphic to the filled Julia set
for the Tchebycheff map, which is equal to the interval \[[-2,2]\].
\vskip -.1in
\midinsert
\insertRaster 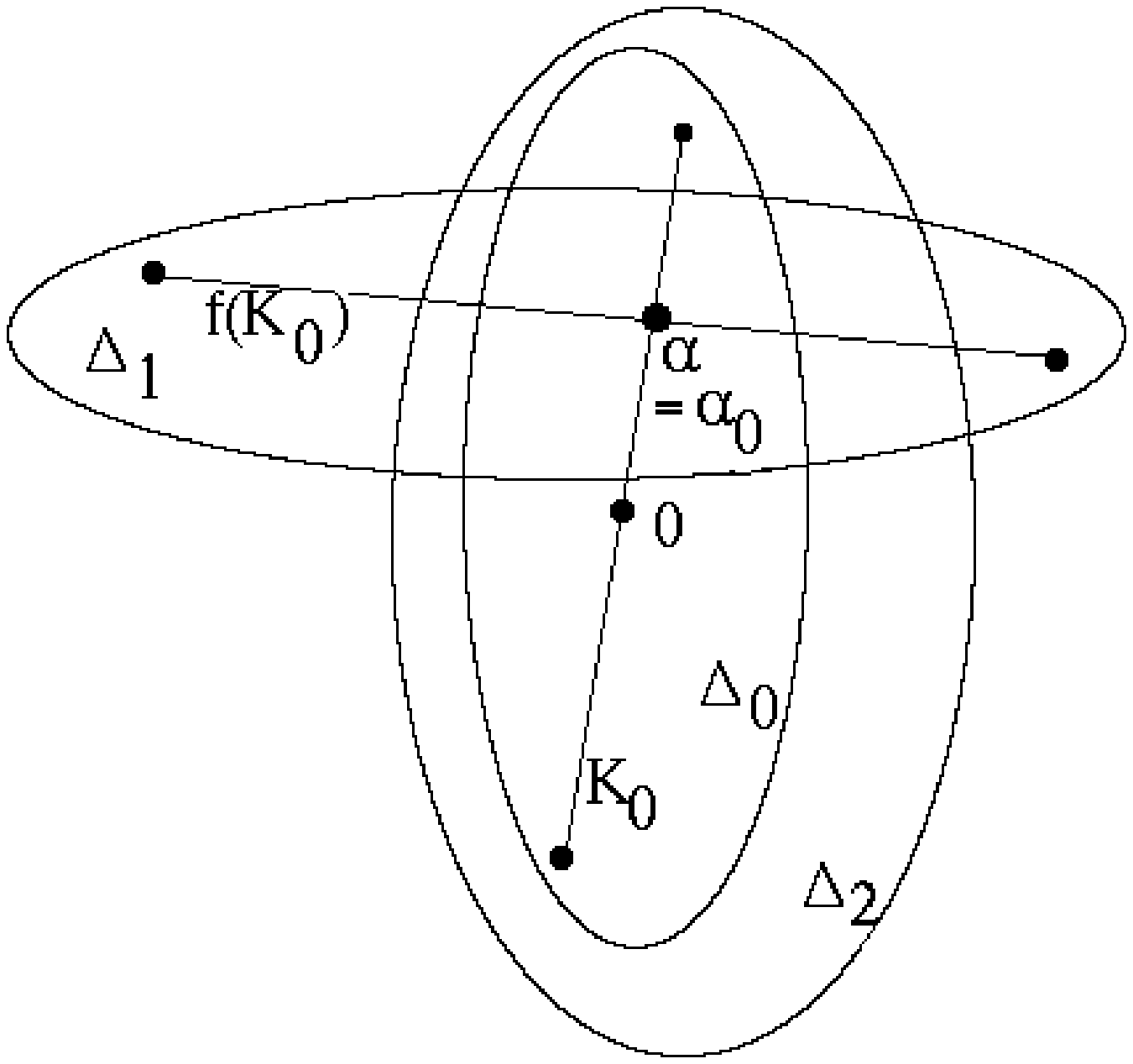 pixels 640 by 480 scaled 450
\vskip -.3in
\centerline{\bit Figure B. Schematic picture of this renormalization.}
\medskip
\endinsert

{\it Note that the critical tableau associated with this map
\[f\] is not periodic.\/}
In fact the orbit of the critical point \[0\] never returns to the puzzle
piece \[P_1(0)\]. It follows that \[f\] cannot be obtained by tuning.\smallskip

As a different example, let \[c\] be any point in a ``satellite''
hyperbolic component which is attached to the cardioid at internal angle
\[p/q\]. (See p. 31.) Then McMullen
shows that \[f_c\] is renormalizable, not only of period \[q\], but also of
period \[n\] for any divisor \[n\] of \[q\].\smallskip

The curious feature of all of these examples is that the disks \[\Delta_0\]
and \[\Delta_1=f(\Delta_0)\] cross over each other in a neighborhood of
the fixed point \[\alpha\]. In fact this fixed point belongs both to the
filled Julia set \[K_0=K(f^{\circ 2}|\Delta_0)\] and to its image \[f
(K_0)\]. (See Figure B.) Furthermore, \[\alpha\] separates each of these sets.
\smallskip

More generally, let \[g=f_{c_1}\] be any renormalizable quadratic map,
with\break
\[g^{\circ n}:\Delta_0\to \Delta_n\supset \Delta_0\]. The
associated filled Julia set
\[K_0=K(g^{\circ n}| \Delta_0)\] has two fixed points which we can label as
\[\alpha_0\] and \[\beta_0\]. Here \[\beta_0\] is
the non-separating fixed point of \[g^{\circ n}| K_0\],
corresponding to the end-point of the zero ray for the hybrid-equivalent
quadratic map. McMullen shows that a
forward image \[g^{\circ i}(K_0)\] with \[0< i<n\]
can intersect \[K_0\] at most in a single repelling fixed point of \[g^{\circ
n}\].\smallskip

{\bf Definition:} The renormalizable map \[g\] is
{\bit simply-renormalizable\/} if the sets \[K_0\] and \[g^{\circ i}
(K_0)\] intersect at most in the non-separating fixed point \[\beta_0\].
\smallskip

This concept will be studied in a forthcoming manuscript by McMullen.
\medskip
%\eject
\centerline{-----------------------------------------------------}\medskip

Here is a more detailed list of corrections and comments.
%, based on this definition.\bigskip
\medskip
\frenchspacing

\noindent {\bf p. 2, Theorem 1:} The statement can
be sharpened by substituting ``simply-renormaliz\-able'' in Hypothesis (3).
Compare the discussion of Lemmas 2 and 3 below.\medskip

\noindent {\bf p. 2, Hypothesis \[\bf (3')\]:} McMullen shows that a map
is infinitely renormalizable if and only if it is infinitely
simply-renormalizable.\medskip

\noindent{\bf p. 3:} For [HF] read [FH].\medskip

\noindent{\bf p. 10. Proof of (b):} For \[d'>n\] read \[d'>d\] and for
\[d'-3k\] read \[d'-2k\].
\medskip

\noindent{\bf p. 12, Lemma 2:}
The converse is false as stated. In fact McMullen shows
that the critical tableau is periodic if and only if \[f\] is
simply-renormalizable. (Here we assume as usual that the fixed point \[\alpha\]
is not parabolic, and also that the critical orbit does not
terminate on \[\alpha\] so that the critical tableau is
well defined.)\medskip

\noindent {\bf p. 12, Lemma 3:} Under the hypothesis of this lemma,
he shows that \[f\] is
simply-renormalizable of period \[q\]. It follows that the 
Corollary on p. 13 applies
whenever \[f\] is not simply-renormalizable.\medskip

\noindent{\bf p. 16, Problem 1-5:} The proof as outlined depends on the claim
that a generic map in the boundary of the Mandelbrot set is
not renormalizable, although the proof sketched in the footnote
shows only that it cannot be obtained by tuning.
McMullen and also J. Kwapisz point out that the easiest way to prove this
statement, as well as Problem 1-5, is by a quite different
argument which makes no reference to puzzles. Note first that the family
of polynomial maps \[c\mapsto f_c^{\circ n}(0)\] fails to be normal in a
neighborhood of \[c\] precisely
for \[c\] belonging to the boundary \[\partial M\]. Let \[U(p,\epsilon)\]
be the set of \[c\in\partial M\] such that the critical orbit of \[f_c\]
intersects the \[\epsilon$-neighborhood of every cycle of period \[p\].
This set is clearly open in \[\partial M\],
and using Montel's Theorem one can check that it is dense when \[p>1\].
Therefore, a generic \[f_c\] in \[\partial M\] belongs to the intersection of
the \[U(p,\epsilon)\]. It follows that
the critical orbit for a generic \[f_c\] is
dense in the Julia set \[J(f_c)\]. As a corollary, it follows that a generic
\[f_c\in\partial M\] is not renormalizable,
since a renormalizable map must have
critical orbit bounded away from its \[\beta\] fixed point.\medskip

\noindent{\bf p. 32, Renormalization and Tuning.} Replace lines 9--3 from the
bottom by:

\noindent It seems reasonable to
conjecture that a map \[f\in M\] is simply-renormalizable of period \[n\]
if and only if 

(1) it belongs to one of the small copies \[H*M\subset M\]
which are obtained by Douady-Hubbard tuning, where \[H\] ranges over
the hyperbolic components of period \[n\], and

(2) it is not the root point of a
satellite hyperbolic component of period \[n\].

\noindent Compare [DH3, p. 332].

\bigskip
\rightline{J. Milnor, Stony Brook, August 1992}

\rightline{Internet: jack@math.sunysb.edu\qquad}

\end

\vfill\eject
%\magnification 1200
\parskip=2pt
\ifnum\pageno<2 \pageno=45 \fi

\ref{\bf References}

\ref [A] L. Ahlfors, Conformal Invariants, McGraw-Hill 1973.

\ref [Be] A. Beardon, Iteration of Rational Functions, Grad. Texts Math.
{\bf 132}, Springer 1991.

\ref [Bl] P. Blanchard, Disconnected Julia sets, pp. 181-201 of ``Chaotic
Dynamics and Fractals'', ed. Barnsley and Demko, Academic Press 1986.

\ref [BDK] P. Blanchard, R. Devaney and L. Keen, The dynamics of complex
polynomials and automorphisms of the shift, Inv. Math {\bf 104} (1991) 545-580.

\ref [Br1] B. Branner, The parameter space for complex cubic polynomials,
pp. 169-179 of ``Chaotic
Dynamics and Fractals'', ed. Barnsley and Demko, Academic Press 1986.

\ref [Br2] B. Branner, The Mandelbrot set, pp. 75-105
of ``Chaos and Fractals'', edit.
Devaney and Keen, Proc. Symp. Applied Math. {\bf 39}, Amer. Math. Soc. 1989.

\ref [BH1] B. Branner and J. H. Hubbard, The iteration of cubic polynomials,
Part I: the global topology of parameter space, Acta Math. {\bf 160} (1988)
143-206;

\ref [BH2] B. Branner and J. H. Hubbard, The iteration of cubic polynomials,
Part II: patterns and parapatterns, Acta Math., to appear.

\ref [CR] R. Courant and H. Robbins, What is Mathematics?, Oxford U. Press
 1941.

\ref [D1] A. Douady, Syst\`emes dynamiques holomorphes,
S\'eminar Bourbaki, $35^e$ ann\'ee 1982\--83, ${\rm n}^o$ 599; Ast\'erisque
{\bf 105-106} (1983) 39-63.

\ref  [D2] A. Douady, Julia sets and the Mandelbrot  set, pp. 161-173
of ``The Beauty of Fractals'', edit. Peitgen and Richter, Springer 1986.

\ref [D3] A. Douady, Algorithms for computing angles in the Mandelbrot set,
pp. 155-168 of ``Chaotic Dynamics and Fractals'', ed. Barnsley \& Demko,
Acad. Press 1986.

\ref [D4] A. Douady, Disques de Siegel et anneaux de Hermann, S\'eminar
Bourbaki, $39^e$ ann\'ee 1986\--87, ${\rm n}^o$ 677; Ast\'erisque {\bf 152-153}
(1987-88) 151-172.
%\eject

\ref [DH1] A. Douady and J. H. Hubbard, It\'eration des polyn\^omes
quadratiques complexes, CRAS Paris {\bf 294} (1982) 123-126.

\ref [DH2] A. Douady and J. H. Hubbard, \'Etude dynamique des polyn\^omes
complexes I \& II, Publ. Math. Orsay (1984-85).

\ref [DH3] A. Douady and J. H. Hubbard, On the dynamics of polynomial-like
mappings, Ann. Sci. Ec. Norm. Sup. (Paris) {\bf 18} (1985), 287-343.

\ref [EL] A. Eremenko and M. Lyubich, The dynamics of analytic
transformations, Leningr. Math. J. {\bf 1} (1990) 563-634.

\ref [F] D. Faught, Local connectivity in a family of cubic polynomials,
thesis, Cornell 1992.

\ref [FH] D. Faught and J. H. Hubbard, Local connectivity of quadratic
Julia sets and the Mandelbrot set, in preparation.

\ref [GM] L. Goldberg and J. Milnor, Fixed point portraits of polynomial maps,
Stony Brook IMS preprint 1990/14.

\ref [He] M. Herman, Recent results and some open questions on Siegel's
linearization  theorem of germs of complex analytic diffeomorphisms of
\[C^n\] near a fixed point, pp. 138-198 of Proc \[8^{th}\] Int. Cong. Math.
Phys., World Sci. 1986.

\ref [Hu1] J. H. Hubbard, Puzzles and quadratic tableaux (according to
Yoccoz), preprint 1990.

\ref [Hu2] J. H. Hubbard, Parapuzzles and local connectivity in the
parameter plane (according to Yoccoz), preprint 1990.

\ref [K] J. Kahn, in preparation.

\ref [L1] M. Lyubich, The dynamics of rational transforms: the topological
picture, Russian Math. Surveys {\bf 41:4} (1986), 43-117.

\ref [L2] M. Lyubich, On the Lebesgue measure of the Julia set of a quadratic
polynomial, Stony Brook IMS preprint 1991/10.

\ref [LM] M. Lyubich and J. Milnor, The Fibonacci unimodal map, Stony Brook
IMS preprint 1991/15.

\ref [McM] C. McMullen, Automorphisms of rational maps, pp. 31-60 of
``Holomorphic Functions and Moduli'', edit. Drasin et al., Springer 1988.

\ref [M1] J. Milnor, Self-similarity and hairiness in the Mandelbrot set,
pp. 211-257 of ``Computers in Geometry and Topology'', edit. Tangora,
Lect. Notes Pure Appl. Math. {\bf 114}, Dekker 1989

\ref [M2] J. Milnor, Dynamics in One Complex Variable, Introductory
Lectures, Stony Brook IMS preprint 1990/5.

\ref [Pe] C. Petersen, On the Pommerenke-Levin-Yoccoz inequality,
preprint, IHES 1991.

\ref [Pr] F. Przytycki, Polynomials in the hyperbolic components, in
preparation.

\ref [Sh] M. Shishikura, The Hausdorff Dimension of the Boundary of the
Mandelbrot Set and Julia Sets, Stony Brook IMS preprint 1991/7.

\ref [S\o] D. E. K. S{\o}rensen, Local connectivity of quadratic
Julia sets, preprint, Tech. Univ. Denmark, Lyngby 1992.

\ref [Su] D. Sullivan, Conformal dynamical systems, pp. 725-752 of ``Geometric
Dynamics'', edit. Palis, Lecture Notes Math. {\bf 1007} Springer 1983.
\end